\documentclass[12pt,leqno]{article}

\usepackage{amsmath}
\usepackage{amsfonts}
\usepackage{url}
\usepackage{amsthm}
\usepackage{color}
\usepackage{xcolor}
\usepackage{graphicx}
\usepackage{bm}

\newtheorem{mydef}{Definition}
\newtheorem{mylem}{Lemma}

\newcommand{\bs}[1]{{\boldsymbol #1}}

\definecolor{hotcolor}{rgb}{1,0,0}
\definecolor{redcolor}{rgb}{0,0,0}
\definecolor{purplecolor}{rgb}{0,0,0}
\definecolor{greencolor}{rgb}{0,0,0}
\definecolor{bluecolor}{rgb}{0,0,0}

\newcommand{\purp}[1]{\textcolor{purplecolor}{#1}}
\newcommand{\revone}[1]{{\textcolor{greencolor}{#1}}}
\newcommand{\revtwo}[1]{{\textcolor{redcolor}{#1}}}
\newcommand{\revthree}[1]{{\textcolor{bluecolor}{#1}}}

\usepackage[margin=1.5cm]{geometry}
\begin{document}
\title{A sharp interface method for an immersed viscoelastic solid}
\author{Charles Puelz$^1$ and Boyce E. Griffith$^{2,3,4}$}
\date{%
    $^1$Courant Institute of Mathematical Sciences, New York University\\%
	$^2$Departments of Mathematics, Applied Physical Sciences, and Biomedical Engineering, University of North Carolina, Chapel Hill\\
	$^3$Carolina Center for Interdisciplinary Applied Mathematics, University of North Carolina, Chapel Hill\\
    $^4$McAllister Heart Institute, University of North Carolina, Chapel Hill\\
    [2ex]%    \today
}
\maketitle

\abstract{The immersed boundary--finite element method (IBFE) is an approach to describing the dynamics of an elastic structure immersed in an incompressible viscous fluid.
In this formulation, there are discontinuities in the pressure and viscous stress at \purp{fluid--structure interfaces.}  
The standard immersed boundary approach, \purp{which connects the Lagrangian and Eulerian variables via integral transforms with regularized Dirac delta function kernels}, smooths out these discontinuities, \purp{which generally} leads to low order accuracy.  
\purp{This paper describes an approach to accurately resolve} pressure discontinuities for these types of formulations, in which the solid may undergo large deformations.
Our strategy is to decompose the physical pressure field into a sum of two pressure--like fields, one defined on the entire \purp{computational} domain, \purp{which includes both the fluid and solid subregions}, and one defined \purp{only} on the solid \purp{subregion}. 
Each of these fields is continuous on its domain of definition, \purp{which enables high accuracy via standard discretization methods without sacrificing} sharp resolution of the pressure discontinuity.  
Numerical tests demonstrate that this method improves rates of convergence for displacements, velocities, stresses, and pressures, as compared to the conventional IBFE method. 
Further, it produces much smaller errors at reasonable numbers of degrees of freedom.  
The performance of this method is tested on several cases with analytic solutions, a nontrivial benchmark problem of incompressible solid mechanics, \purp{and an example involving a thick, actively contracting torus.}}

%\vspace{0.5cm}
%\noindent
%(Reviewer 1 revisions in \revone{green}, Reviewer 2 revisions in \revtwo{red}, and Reviever 3 revisions in \revthree{blue}.)

\section{Introduction}

The immersed boundary (IB) method is a general approach to modeling fluid--structure interaction that was introduced by Peskin \cite{Peskin72, Peskin77} for heart valve dynamics.  
Its strength lies in the representation of the fluid in Eulerian form, which enables approximation on a fixed, Cartesian mesh, and the representation of the solid in Lagrangian form.
\purp{This approach is appealing because it does not require discretizations that conform to the fluid--structure interface, but} instead relies on integral transforms with \purp{Dirac} delta function kernels to describe interactions between the Eulerian and Lagrangian frames.  
In practice, conventional IB methods use discretizations of these transforms with regularized delta function kernels. \purp{This has the effect of smoothing discontinuities} in the pressure and viscous stress that generically appear at the interface between the fluid and solid.

\revthree{Many techniques have been developed to improve the accuracy in numerical treatments of discontinuities arising in \revthree{elliptic problems} and fluid--structure interaction formulations. One set of approaches focuses on immersed interfaces, i.e., structures that have codimension 1 with respect to the ambient space.
\revthree{LeVeque and Li considered elliptic interface problems and altered finite difference stencils to properly incorporate jump conditions in the discretized equations \cite{Leveque94}.  We also note the work of Bedrossian et al. in which interface conditions are enforced using a Lagrange multiplier approach \cite{Bedrossian10}.}}

  With respect to codimension 1 interfaces immersed in an incompressible viscous fluid, relevant work includes that of Ye et al. \cite{Ye99} for stationary structures and Udaykumar et al. \cite{Udaykumar01}, for moving structures.
In these methods, the equations of motion are approximated using finite volume methods, and modifications to the discrete operators are developed that account for jump conditions at the immersed interface.  
Related work by Seo et al. \cite{Seo11} uses a cut--cell approach to obtain higher volume accuracy. 
Local finite--volume approximations are used in cells that are cut by the interface to modify the equations of motion to account for the imposed jump or boundary conditions.  
\purp{Lee and LeVeque also considered this problem, and imposed jump conditions directly in the Poisson solve for the pressure field used in the discrete solution to the incompressible Navier--Stokes equations \cite{Lee03}.}
\revthree{The book by Li and Ito contains a review of some numerical approaches for these types of problems \cite{Li06}.}

The focus of this paper is on a sharp interface method for solids with codimension 0.  
Recent work in high order methods for stationary solids with codimension 0 includes the immersed boundary smooth extension method (IBSE) of Stein et al. \cite{Stein16}.  
Their approach involves the construction of smooth extensions of solutions from the fluid domain to an {\em extension} domain that typically overlaps with the solid region.  
In the IBSE method, these extensions are spectrally computed from a sequence of harmonic problems and are used to modify forcing functions in the problem formulation. 
Our approach is similar to the IBSE method in that we solve a harmonic problem and modify the forcing function so the pressure field is continuous.  Unlike the approach of Stein et al., however, our approach is able to handle solids undergoing large deformations. 
We also employ standard finite--volume methods for the Eulerian equations and finite--element methods for the Lagrangian equations.   

This study uses the immersed boundary--finite element method (IBFE), in which solid displacements and forces are approximated using finite element discretizations \cite{Griffith17}. 
The first version of such a method appeared in the work of Wang and Liu \cite{Wang04}, and was expanded upon in the work of Zhang et al.~\cite{Zhang04} and Liu et al.~\cite{Liu06}. 
Our paper uses a version of the immersed boundary finite element method based on one introduced by Boffi et al.~\cite{Boffi08}.  
Boffi et al. use finite element discretizations for the Eulerian and Lagrangian variables, treat the delta function kernel in a variational way, and obtain suboptimal convergence rates for test problems with respect to the approximation spaces used.  
In our work, we use the model described by Boffi et al. \cite{Boffi08}, a finite element approximation for the discretization of the solid, and a second-order staggered--grid finite volume scheme for the discretization of the fluid.  
This approach was detailed by Griffith and Luo \cite{Griffith17}.

The method described herein follows the method of Griffith and Luo, but avoids regularizing pressure discontinuities at fluid--structure interface by splitting the physical pressure into the sum of two pressure--like fields that are continuous on their domains of definition. 
Specifically, we use a global pressure--like field that is defined over the entire Eulerian computational domain along with a second field that is restricted to the Lagrangian domain. 
Boundary conditions for the Lagrangian pressure--like field ensure that the total pressure recovers the correct jump conditions.
This allows us to obtain higher--order accuracy while using only standard discretization methods. \revone{Our approach requires one additional harmonic or diffusion equation solve per timestep, compared to the conventional formulation of the IBFE method.} 

\section{Equations of motion and jump conditions}
\label{sec:motion}

\purp{We consider a coupled fluid--structure interaction problem in a domain $\Omega = \Omega_t^\text{s} \cup \Omega_t^\text{f}$, with the solid and fluid domains at time $t$ defined to be $\Omega_t^\text{s}$ and $\Omega_t^\text{f}$ respectively.
The reference configuration of the solid is denoted by $U$. 
The motion map for the solid is defined to be $\bm{\chi}(\cdot, t)$ so that ${\bs x} = \bm{\chi}({\bs X},t)$ is the current position at time $t$ of the reference coordinate ${\bs X} \in U$.
The image $\bm{\chi}(U,t)$ is the current configuration of the solid at time $t$, so that $\Omega_t^\text{s} = \bm{\chi}(U,t)$ and $\Omega_t^\text{f} =\Omega \sim  \bm{\chi}(U,t)$.}
  
The fluid--structure system that we consider in this paper is described by the total Cauchy stress tensor,
\begin{align*}
\mathbb{\sigma}(\bm{x},t) = -\mathbb{I} p(\bm{x},t) + \mu ( \nabla {\bs u}(\bm{x},t) + \nabla {\bs u}(\bm{x},t)^T ) + 
\begin{cases}
{\bs \sigma}^\text{e}(\bm{x},t), \quad &\bm{x} \in \Omega_t^\text{s}, \\
{\bs 0}, \quad &\bm{x} \in \Omega_t^\text{f},
\end{cases}
\end{align*}
in which $\bs \sigma^\text{e}$ is the elastic Cauchy stress tensor of the solid, $\mu$ is the dynamic viscosity, $p$ is the pressure, and ${\bs u}$ is the velocity.  
Denote the deformation gradient $\mathbb{F} = \partial {\bs \chi}/{\partial {\bs X}}$ and its determinant $J =\text{det}(\mathbb{F})$. 
The first Piola--Kirchoff elastic stress corresponding to ${\bs \sigma}^\text{e}$ is
\begin{align*}
\mathbb{P}^\text{e}({\bs X},t) = J({\bs X},t)\, {\bs \sigma}^\text{e}({\bs \chi}({\bs X},t),t)\, \mathbb{F}^{-T}({\bs X},t). 
\end{align*}
As shown by Boffi et al. \cite{Boffi08}, the equations of motion for this fluid--structure system are:
\begin{align}
\label{eq:em1}
\rho \left(\frac{\partial \bm{u}}{\partial t}(\bm{x},t) + \bm{u} (\bm{x},t)   \cdot \nabla \bm{u}(\bm{x},t) \right) &= -\nabla p(\bm{x},t)  
+ \mu \nabla^2 \bm{u}(\bm{x},t) + \bm{f}(\bm{x},t), && \revthree{\bm{x} \in \Omega,}\\
\label{eq:em2}
\nabla \cdot \bm{u}(\bm{x},t) &= 0, && \revthree{\bm{x} \in \Omega,} \\
\label{eq:em3}
\bm{f}(\bm{x},t)  &= \int_{U} \nabla_{\bm X} \cdot \mathbb{P}^{\text e}(\bm{X},t)\, \delta(\bm{x} - \bm{\chi}(\bm{X},t))\,d\bm{X},  \\
\nonumber
&\quad  -\int_{\partial U} \mathbb{P}^{\text e}(\bm{X},t)\bm{N}(\bm{X})\, \delta(\bm{x} - \bm{\chi}(\bm{X},t))\,dA, && \revthree{\bm{x} \in \Omega,} \\
\label{eq:em5}
\frac{\partial \bm{\chi}}{\partial t}(\bm{X},t) &= \int_\Omega \bm{u}(\bm{x},t) \, \delta(\bm{x} - \bm{\chi}(\bm{X},t))\, d\bm{X}, && \revthree{\bm{X} \in U,}
\end{align}
in which $\rho$ is the density of the fluid.  
Equation \eqref{eq:em1} expresses balance of momentum for the fluid and solid.  Incompressibility is imposed in equation \eqref{eq:em2}.  
Equation \eqref{eq:em3} explicity defines the force density exerted from the solid onto the fluid in terms of the first Piola--Kirchoff stress.  
This term includes a volumetric force density and a surface force density, the latter of which generates discontinuities in the pressure, derivatives of the pressure, and derivatives of the velocity at the fluid--structure interface.  
Equation \eqref{eq:em5} requires the velocity of the solid to equal the velocity of the background fluid, the so called no--penetration and no--slip condition.

In this paper, we assume the solid is hyperelastic.  This implies that the stress is determined by a strain energy density $W$ via
\begin{align*}
\mathbb{P}^\text{e} = \frac{\partial W}{\partial \mathbb{F}}.
\end{align*}
This assumption is not a requirement of the IBFE approach or our sharp interface method.

We recall some results for deriving a jump condition in the pressure.  Sketches of proofs can be found in the work of Lai and Li \cite{Lai01} and Peskin and Printz \cite{Peskin93}, and we recall them here in our notation for completeness. 
Let ${\bs n} = {\bs n}({\bs x},t)$ and ${\bs N} = {\bs N}({\bs X})$ denote the outward unit normal vectors to the solid region in the current and reference configurations respectively. 
The following notation is useful in discussing jumps of variables at the fluid--structure interface in the current configuration, defined as $\Gamma_t^\text{fs} = \partial \Omega_t^\text{s}$.  

\begin{mydef}
\label{def:jump}
The jump of a scalar valued function across the fluid--structure interface is defined as:
\begin{align*}
[g({\bs x})] = \lim_{\varepsilon \rightarrow 0} g({\bs x} + \varepsilon {\bs n}) - \lim_{\varepsilon \rightarrow 0} g({\bs x} - \varepsilon {\bs n}), \quad {\bs x} \in \Gamma_t^\emph{fs}.
\end{align*}
The jump of vector or tensor valued variables is defined in the same way, componentwise.
\end{mydef}
When we use the square bracket notation for the jump described above, we implicitly assume this jump is evaluated at the fluid--structure interface $\Gamma_t^\text{fs}$.  The following results are necessary in deriving a jump condition for the pressure under the  assumption the velocity field is continuous.
\begin{mylem}
\label{lem:tanvel}
Let ${\bs t}$ and ${\bs b}$ be the unit tangent vectors at the fluid--structure interface.  The tangential derivatives of the velocity are continuous, i.e.
\begin{align*}
[(\nabla {\bs u})\, {\bs t}] = [(\nabla {\bs u})\, {\bs b}] = 0. 
\end{align*}
\begin{proof}
Consider a parametrized curve ${\bs \beta} = {\bs \beta}(s)$ defined on $\Gamma_t^\text{fs}$ which contains the point at which we consider the jump.  This curve is constructed so its tangent vector $d {\bs \beta}/ds$ is equal to ${\bs t}$ at this point.  Because the velocity field is continuous, we can consider the velocity evaluated along this curve.  The derivative of a component of the velocity field $u_i$ along this curve is the tangential derivative,
\begin{align*}
\frac{d}{ds} u_i({\bs \beta}(s)) = \frac{d {\bs \beta}}{ds} \cdot \nabla u_i ({\bs \beta}(s)) = {\bs t} \cdot \nabla u_i ({\bs \beta}(s)). 
\end{align*}
This calculation shows the tangential derivative is defined along the boundary and must be continuous.  A similar argument is applied for the derivative in the direction of ${\bs b}$ and for the other components of the velocity field.
\end{proof}
\end{mylem}
\begin{mylem}
\label{lem:vel}
$[{\bs n} \cdot (\nabla {\bs u})\,{\bs n}] = [{\bs n} \cdot (\nabla {\bs u})^T\,{\bs n}] = 0.$
\begin{proof}
We follow the argument in \cite{Lai01}.  
To establish notation, let \purp{${\bs x}^\star$} be the point in the current configuration on $\Gamma_t^\text{fs}$ where we consider the jump, let ${\bs n} = (n_x,n_y,n_z)$ be the unit normal vector at \purp{${\bs x}^\star$}, and let the orthonormal pair ${\bs t} = (t_x,t_y,t_z)$ and ${\bs b} = (b_x,b_y,b_z)$ span the tangent space.  
Given an arbitrary point ${\bs x}$, consider a linear transformation to a new point $\hat{\bs x}$ defined by
\begin{align*}
\hat{\bs x} = 
\begin{pmatrix}
n_x & n_y & n_z \\
t_x & t_y & t_z \\
b_x & b_y & b_z
\end{pmatrix}
({\bs x} - \purp{{\bs x}^\star}) := {\bs T}\,({\bs x} - \purp{{\bs x}^\star}),
\end{align*}  
in which $\hat{\bs x}$ are local coordinates for ${\bs x}$ expressed in the basis $\{{\bs n}, {\bs t}, {\bs b}\}$, and ${\bs T}$ is the matrix containing these vectors in its rows.  
The notation $\nabla$ refers to derivatives with respect to the physical Cartesian coordinates. 
At this point, we invoke continuity of tangential derivatives of $\nabla {\bs u}$ from Lemma \ref{lem:tanvel}. By the chain rule, on $\Gamma_t^\text{fs}$ we have
\begin{align*}
0 = [(\nabla {\bs u})\, {\bs t}] = [(\nabla_{\hat{\bs x}} {\bs u})\, {\bs T}\, {\bs t}].
\end{align*}
With ${\bs u} = (u_1, u_2, u_3)$ and $\hat{\bs x} = (\hat{x}, \hat{y}, \hat{z})$, by orthogonality of the first and last rows of ${\bs T}$ with respect to ${\bs t}$, the above statement reads
\begin{align}
\label{eq:j1}
\left[\frac{\partial u_1}{\partial \hat{y}}\right] = \left[\frac{\partial u_2}{\partial \hat{y}}\right] = \left[\frac{\partial u_3}{\partial \hat{y}}\right] = 0.
\end{align}
The same argument with ${\bs b}$ gives
\begin{align}
\label{eq:j2}
\left[\frac{\partial u_1}{\partial \hat{z}}\right] = \left[\frac{\partial u_2}{\partial \hat{z}}\right] = \left[\frac{\partial u_3}{\partial \hat{z}}\right] = 0.
\end{align}
The incompressibility condition, written in terms of $\hat{\bs x}$, along with \eqref{eq:j1}--\eqref{eq:j2} implies
\begin{align*}
0 = [\nabla \cdot {\bs u}] = \left[\frac{\partial u_1}{\partial \hat{x}}\right]\frac{\partial \hat{x}}{\partial x} + \left[\frac{\partial u_2}{\partial \hat{x}}\right]\frac{\partial \hat{x}}{\partial y} + \left[\frac{\partial u_3}{\partial \hat{x}}\right]\frac{\partial \hat{x}}{\partial z}.
\end{align*}
Finally, we examine
\begin{align*}
{\bs n} \cdot [\nabla {\bs u}]\, {\bs n} = {\bs n} \cdot [\nabla_{\hat{\bs x}} {\bs u}]\, {\bs T}\, {\bs n} = \left[\frac{\partial u_1}{\partial \hat{x}}\right]\frac{\partial \hat{x}}{\partial x} + \left[\frac{\partial u_2}{\partial \hat{x}}\right]\frac{\partial \hat{x}}{\partial y} + \left[\frac{\partial u_3}{\partial \hat{x}}\right]\frac{\partial \hat{x}}{\partial z} = 0.
\end{align*}
\end{proof}
\end{mylem}
The discontinuity in the pressure field is derived in the following lemma.
\begin{mylem}
\label{lem:press}
The pressure satisfies the following jump condition at the fluid--structure interface:
\begin{align*}
[p] = -{\bs n} \cdot {\bs \sigma}^\text{e} \,{\bs n}. 
\end{align*}
\begin{proof}
Continuity of the traction vector $[{\bs \sigma} \, {\bs n}] = 0$ on $\Gamma_t^\text{fs}$ implies:
\begin{align*}
-[p] {\bs n} + \mu ( [\nabla {\bs u}] + [\nabla {\bs u}^T])\, {\bs n} - {\bs \sigma}^\text{e}\, {\bs n} = 0. 
\end{align*}
Taking the inner product of this equation with ${\bs n}$ and using Lemma \ref{lem:vel}, we obtain the jump condition for the pressure.
\end{proof}
\end{mylem}

\section{Description of the method}
\label{sec:method}

%We first rewrite the pressure jump condition in Lemma \ref{lem:press} in the reference configuration.  Using the definition of the first Piola--Kirchoff stress:
%\begin{align*}
%{\bs \sigma}^\text{e}\,{\bs n}\,da = \mathbb{P}^\text{e}\,{\bs N}\,dA
%\end{align*}

The method introduced here splits the physical pressure field $p$ into a sum of two components, $\varphi$ and $\pi$.  
The field $\varphi$ is defined on the solid domain $\Omega_t^\text{s}$ as a solution to a harmonic problem, to be specified below.  
The field $\pi$ is defined on $\Omega$ and, by choosing appropriate boundary conditions for $\varphi$, is continuous at the fluid--structure interface.  
Once $\varphi$ and $\pi$ are known, the physical pressure is recovered by $p = \pi + \varphi$.  
Explicitly, the physical pressure is defined as:
\begin{align*}
p({\bs x},t) = \pi({\bs x},t) + 
\begin{cases}
\varphi({\bs x},t), \quad &\text{if } {\bs x} \in \Omega_t^\text{s}, \\
0, \quad &\text{otherwise.}
\end{cases}
\end{align*}
This splitting is achieved by modifying the elastic stress ${\bs \sigma}^\text{e}$ so the normal component of the modified traction vanishes on $\Gamma_t^\text{fs}$.  
We take the field $\varphi$ to modify the first Piola--Kirchoff stress in the following way:
\begin{align*}
\tilde{\mathbb{P}}^\text{e} = \mathbb{P}^\text{e} - J\, \varphi\, \mathbb{F}^{-T}.
\end{align*}
Note that $\varphi$ has units of pressure, and we define a modified elastic Cauchy stress $\tilde{\bs \sigma}^\text{e}$ as:
\begin{align*}
\tilde{{\bs \sigma}}^\text{e} = J^{-1} \,\tilde{\mathbb{P}}^\text{e}\, \mathbb{F}^T = {\bs \sigma}^\text{e} - \varphi\, \mathbb{I}.
\end{align*}
To encode the pressure discontinuity in the boundary condition for $\varphi$, we require that it satisfy the interface condition:
\begin{align*}
\varphi({\bs x},t) = {\bs n} \cdot {\bs \sigma}^\text{e}\,{\bs n}, \quad {\bs x} \in \Gamma_t^\text{fs}.
\end{align*}
This condition removes the normal component of the modified traction at the fluid--structure interface, rendering the pressure--like field continuous; refer to Lemma \ref{lem:press}. 
The pressure--like field $\pi$ is defined to be the pressure solution to the equations of motion when solved with this modified stress. Notice that $\pi$ remains a Lagrange multiplier for $\nabla \cdot {\bs u} = 0$, and its value is completely determined (at least up to an additive constant, depending on boundary conditions) by imposing this constraint. 

It is convenient to compute $\varphi$ in the reference configuration of the solid. To express the interface condition for $\varphi$ in the reference configuration, we use Nanson's relation and the definition of the first Piola--Kirchoff stress,
\begin{align*}
{\bs n} \, da = J\, dA\, \mathbb{F}^{-T}{\bs N} \quad \text{and} \quad {\bs \sigma}^\text{e}{\bs n}\, da = \mathbb{P}^\text{e} {\bs N} \,dA. 
\end{align*} 
These equations imply 
\begin{align*}
{\bs n} = \frac{\mathbb{F}^{-T} {\bs N}}{\|\mathbb{F}^{-T} {\bs N}\|} \quad \text{and} \quad {\bs \sigma}^\text{e}\,{\bs n} = J^{-1}\frac{\mathbb{P}^\text{e}\,{\bs N}}{\|\mathbb{F}^{-T}{\bs N}\|},
\end{align*}
which results in an interface condition for $\varphi$ formulated in the reference configuration:
\begin{align*}
\varphi({\bs X},t) = J^{-1}\frac{\mathbb{F}^{-T} {\bs N}}{\|\mathbb{F}^{-T} {\bs N}\|^2} \cdot \mathbb{P}^\text{e}\,{\bs N}, \quad {\bs X} \in \partial U.
\end{align*}
We consider two approaches for computing $\varphi$. 
The first approach is to compute $\varphi$ as a solution to a {\em steady state} harmonic problem at each timestep:
\begin{align}
\label{eq:phi1}
-\nabla^2_{\bs X}\, \varphi({\bs X}, t) &= 0, \\
\label{eq:phi2}
\varphi({\bs X},t) &= J^{-1}\frac{\mathbb{F}^{-T} {\bs N}}{\|\mathbb{F}^{-T} {\bs N}\|^2} \cdot \mathbb{P}^\text{e}\,{\bs N}, \quad {\bs X} \in \partial U.
\end{align}
The second approach requires $\varphi$ to solve a {\em diffusion equation} each timestep, for some diffusion constant $\gamma$:
\begin{align}
\label{eq:dtphi1}
\frac{\partial}{\partial t} \varphi({\bs X}, t) - \gamma\,\nabla^2_{\bs X}\, \varphi({\bs X}, t) &= 0, \\
\label{eq:dtphi2}
\varphi({\bs X},t) &= J^{-1}\frac{\mathbb{F}^{-T} {\bs N}}{\|\mathbb{F}^{-T} {\bs N}\|^2} \cdot \mathbb{P}^\text{e}\,{\bs N}, \quad {\bs X} \in \partial U, \\
\label{eq:dtphi3}
\varphi({\bs X},0) &= \varphi_0({\bs X}).
\end{align}
The initial condition $\varphi_0$ is taken to be the solution to the harmonic problem at time $t = 0$.  
%The diffusion equation is discretized with Crank-Nicolson in time, resulting in a linear system with milder condition number compared to the one obtained from the steady state harmonic equation \eqref{eq:phi1}--\eqref{eq:phi2}.  
The parameter $\gamma$ can help to control the condition number for the linear system arising from spatial discretization.
\revone{We remark that the solutions to either \eqref{eq:phi1}--\eqref{eq:phi2} or \eqref{eq:dtphi1}--\eqref{eq:dtphi3} are generally different, resulting in a different splitting of the physical pressure $p = \pi + \varphi$ on the current configuration of the solid.} 
%%%%%%
% some settings for 0.38 0 figures I am documenting here:
% legend position: 0.38 0.125
% xscale = 150%, yscale = 70%
% font height for legend: 0.06
%%%%%%

The equations of motion for this sharp interface formulation of the IBFE scheme, using the modified stress and the steady state harmonic equation for $\varphi$, can be stated as:

\begin{align}
\label{eq:modem1}
\rho \left(\frac{\partial \bm{u}}{\partial t}(\bm{x},t) + \bm{u} (\bm{x},t)   \cdot \nabla \bm{u}(\bm{x},t) \right) &= -\nabla \pi(\bm{x},t)  
+ \mu\nabla^2 \bm{u}(\bm{x},t) + \tilde{\bm{f}}(\bm{x},t), && \revthree{\bm{x} \in \Omega,} \\
\label{eq:modem2}
\nabla \cdot \bm{u}(\bm{x},t) &= 0, && \revthree{\bm{x} \in \Omega,}\\
\label{eq:modem3}
\tilde{\bm{f}}(\bm{x},t)  &= \int_{U} \nabla_{\bm X} \cdot \tilde{\mathbb{P}}^{\text e}(\bm{X},t)\, \delta(\bm{x} - \bm{\chi}(\bm{X},t))\,d\bm{X},  \\
\nonumber
&\quad  -\int_{\partial U} \tilde{\mathbb{P}}^{\text e}(\bm{X},t)\bm{N}(\bm{X})\, \delta(\bm{x} - \bm{\chi}(\bm{X},t))\,dA, && \revthree{\bm{x} \in \Omega,} \\
\label{eq:modem5}
\frac{\partial \bm{\chi}}{\partial t}(\bm{X},t) &= \int_\Omega \bm{u}(\bm{x},t) \, \delta(\bm{x} - \bm{\chi}(\bm{X},t))\, d\bm{X}, && \revthree{\bm{X} \in U,} \\
\label{eq:modem6}
\nabla^2_{\bs X}\, \varphi({\bs X}, t) &= 0, && \revthree{\bm{X} \in U,}\\
\label{eq:modem7}
\varphi({\bs X},t) &= J^{-1}\frac{\mathbb{F}^{-T} {\bs N}}{\|\mathbb{F}^{-T} {\bs N}\|^2} \cdot \mathbb{P}^\text{e}\,{\bs N}, && \revthree{\bm{X} \in \partial U,} \\
\label{eq:modem8}
\tilde{\mathbb{P}}^\text{e}({\bs X},t) &= \mathbb{P}^\text{e}({\bs X},t) - J\, \varphi({\bs X},t)\, \mathbb{F}^{-T}, && \revthree{\bm{X} \in U.}
\end{align}
The equations are similar for the case when $\varphi$  is computed via the diffusion equation and are omitted for brevity. 
%Numerical methods for these equations are discussed in the next section.
%where we incorporated the first approach \eqref{eq:phi1}--\eqref{eq:phi2} for computing $\varphi$.  The second approach is added to the equations of motion in an identical way. 

\section{Numerical approximation}
\label{sec:numerics}

The numerical approximation for the equations of motion follows \cite{Griffith17}.  
The time step is denoted $\Delta t$.  The Cartesian grid spacing parameter is denoted $h$.  
The mesh factor $M_\text{fac}$, referenced in some of the results below, corresponds to the ratio between the approximate edge length in the solid finite element mesh and the Cartesian grid spacing $h$. 
Solid displacements and forces are approximated using finite elements via the ``unified weak formulation'' as described by Griffith and Luo \cite{Griffith17}, for which we seek an approximate volumetric force density ${\bs G} = {\bs G}({\bs X},t)$ that is variationally equivalent to the sum of the volumetric ($\nabla_{\bs X} \cdot \mathbb{P}^\text{e}$) and surface ($\mathbb{P}^\text{e}\, {\bs N}$) Lagrangian force densities in \eqref{eq:em3}. 
More explicitly, given some finite element space $\mathcal{F}_h$, the approximate volumetric force density satisfies:
\begin{align*}
\int_U {\bs G}({\bs X},t) \cdot {\bs V}_h({\bs X})\,d {\bs X} = -\int_U \mathbb{P}^\text{e}({\bs X}) : \nabla_{\bs X}{\bs V}_h({\bs X})\,d {\bs X}, \quad \text{for all } {\bs V}_h \in \mathcal{F}_h. 
\end{align*}
%This force density is spread onto the grid using a delta function kernel, and the unified weak formulation is stated as follows:
%\begin{align}
%\label{eq:uw1}
%\rho \left(\frac{\partial \bm{u}}{\partial t}(\bm{x},t) + \bm{u} (\bm{x},t)   \cdot \nabla \bm{u}(\bm{x},t) \right) &= -\nabla p(\bm{x},t)  
%+ \mu\nabla^2 \bm{u}(\bm{x},t) + \bm{g}(\bm{x},t), \quad && {\bs x} \in \Omega, \\
%\label{eq:uw2}
%\nabla \cdot \bm{u}(\bm{x},t) &= 0, \quad && {\bs x} \in \Omega, \\
%\label{eq:uw3}
%\bm{g}(\bm{x},t)  &= \int_{U} {\bs G}(\bs{X},t)\, \delta(\bm{x} - \bm{\chi}(\bm{X},t))\,d\bm{X}, \quad && {\bs X} \in U, \\
%\label{eq:uw4}
%\frac{\partial \bm{\chi}}{\partial t}(\bm{X},t) &= \int_\Omega \bm{u}(\bm{x},t) \, \delta(\bm{x} - \bm{\chi}(\bm{X},t))\, d\bm{X}, \quad && {\bs X} \in U.
%\end{align}
Nodal Lagrange finite elements are used for the finite element space $\mathcal{F}_h$ defined on the solid mesh, and in all simulations, we use bilinear ($Q^1$) elements.  

The equations for $\varphi$, either the steady state harmonic equation \eqref{eq:phi1}--\eqref{eq:phi2} or the diffusion equation \eqref{eq:dtphi1}--\eqref{eq:dtphi3}, are discretized in space using a standard finite element discretization.  Bilinear elements are also used for the finite element space, and Dirichlet boundary conditions are imposed via the penalty method.  The diffusion equation is discretized in time using Crank--Nicolson.  The resulting linear systems are solved with GMRES, using an incomplete $LU$ preconditioner.  For all numerical experiments, relative linear solver tolerances are set to $1\times 10^{-12}$.

The fluid equations are discretized using a second--order accurate staggered grid method, in which the pressure is approximated at the cell center and the velocity components are approximated at the edges (in two spatial dimensions) or faces (in three spatial dimensions) of the Cartesian grid cell.  
Lagrangian--Eulerian interaction through the delta function kernels is discretized by designing approximations to these operators that are discrete adjoints.  Details are provided in prior work \cite{Griffith17}.

\section{Results}
\label{sec:results}

To test this method, we consider four examples.  
The first two involve a thick ring that is either pre--stressed in its reference configuration or inflated to a steady--state final configuration with a fluid source.  
Both cases involve substantial pressure discontinuities along the fluid--structure interface, which we demonstrate are well resolved by the method introduced herein.  
Comparisons to analytic solutions reveal the overall errors in velocities, displacements, stresses and pressures are much smaller with our method, as compared to the conventional IBFE formulation, and converge at faster rates.  
The third example describes compression of a neo--Hookean block of material. 
This problem is more challenging because it contains discontinuities in the surface forces that generate solid displacements on the the top and bottom of the block.  
Further, the material model for the compressed block includes a volumetric energy useful for penalizing compressible deformations \cite{Vadala18}, but this term prominently contributes to discontinuities at the fluid--structure interface that present additional numerical challenges.  
\revtwo{The fourth example, inspired by the work of McQueen and Peskin \cite{McQueen89}, applies the method to an actively contracting thick torus.  The contraction produces a transient velocity field and also large pressures in the toroidal wall which leads to pressure discontinuities at the fluid--structure interface.}

\subsection{Thick orthotropic static ring}

For our first example, we consider a two--dimensional thick ring introduced in by Boffi et al.~\cite{Boffi08}.  
The Lagrangian \purp{curvilinear} coordinates of the ring are ${\bs s} = (s_1, s_2) \in U = [0,2\pi R]\times[0,w]$, and its initial configuration is defined by the motion map at $t = 0$:
\begin{align*}
{\bs \chi}({\bs s}, 0) = \big(\cos(s_1/R)(R+s_2) + 0.5,\, \sin(s_1/R)(R+s_2) + 0.5 \big),
\end{align*}
with $R = 0.25$ mm and $w = 0.0625$ mm.  
The constitutive model is taken to be
\begin{align*}
\mathbb{P}^e = \frac{\mu_e}{w} \mathbb{F},
\end{align*}
in which $\mu_e$ is the elastic stiffness coefficient.
Note that in this model, the Lagrangian coordinates are not the same as the reference coordinates.  
This model is initially at equilibrium, and based on the definition of the motion map from curvilinear Lagrangian coordinates to reference coordinates, leads to a discontinuous pressure field at the fluid--structure interface.  
An analytic solution is available for this problem \cite{Boffi08}, and we use it to compute convergence rates for the pressure and velocity.  
Under the condition that the mean of the pressure field is zero, and with $p_0 = \frac{\pi \mu_e}{3 w}\left( 3wR + R^2 - \frac{(R+w)^3}{R} \right)$ and $r = \|{\bs x} - (0.5, 0.5)\|$, the exact pressure field is:
\begin{align*}
p({\bs x}, t) = 
\begin{cases}
 p_0 + \mu_e \left( \frac{1}{R} - \frac{1}{R+w} \right) \quad  &r \leq R, \\
  p_0 + -\frac{\mu_e}{w} \left( \frac{1}{R}(R+w - r) + \frac{R}{R + w} \right) \quad &R \leq r \leq R+w, \\
p_0 \quad &R+w < r. 
\end{cases}
\end{align*}
%
%\begin{figure}[h!]
%\begin{center}
%\includegraphics[scale=0.125,trim=0 0 -100 0]{IMAGES/trimmed_cart_phi_static_ring.png}
%\includegraphics[scale=0.125,trim=-100 0 0 0]{IMAGES/trimmed_cart_pi_static_ring.png}
%\caption{The pressure--like fields for the splitting approach, using the steady state formulation for $\varphi$.  The field $\varphi$ is on the left and $\pi$ is on the right.}
%\end{center}
%\end{figure}
%

\begin{figure}[h!]
\begin{center}
\includegraphics[scale=0.11,trim=0 0 -250 0]{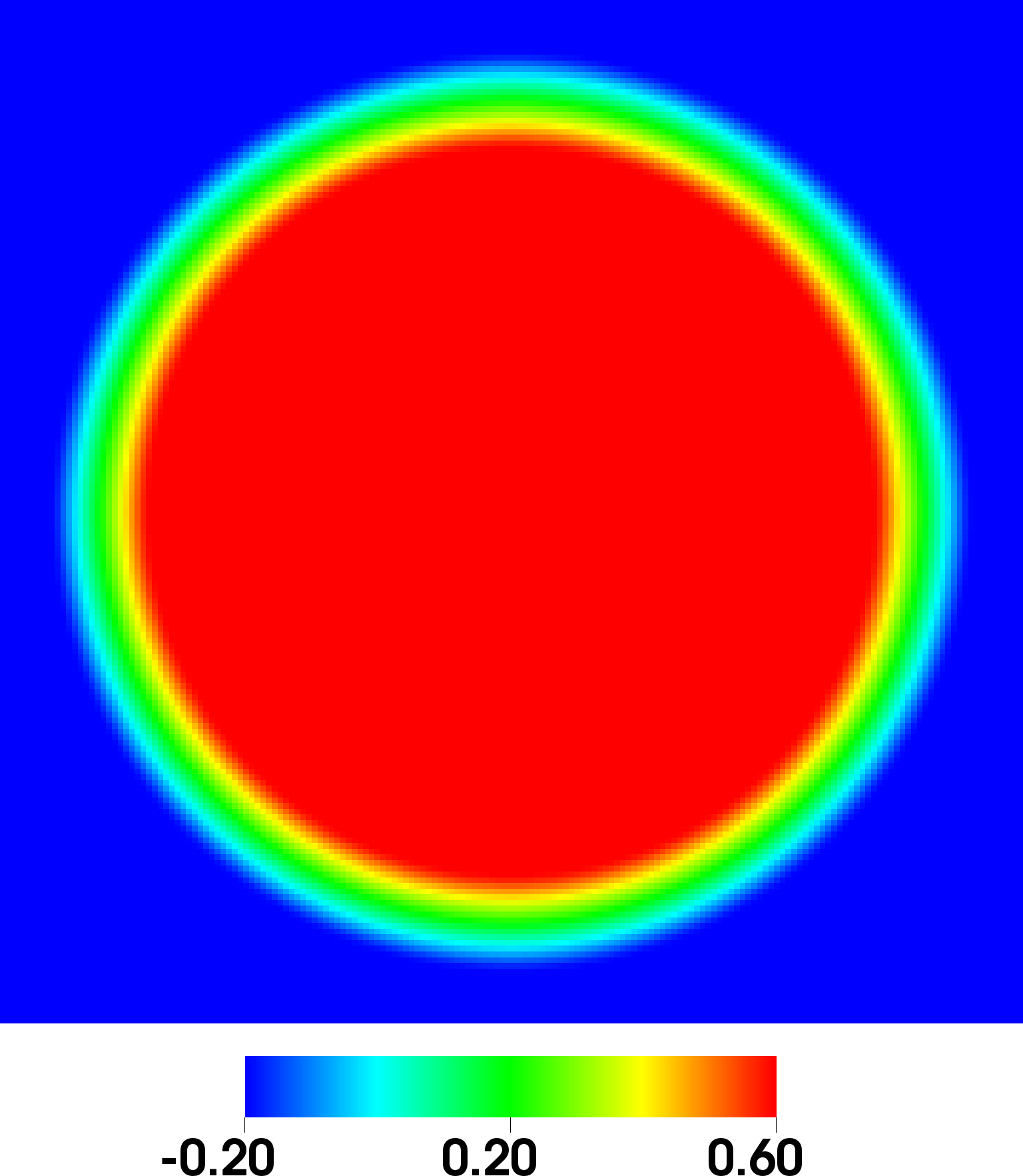}
\includegraphics[scale=0.11,trim=-250 0 0 0]{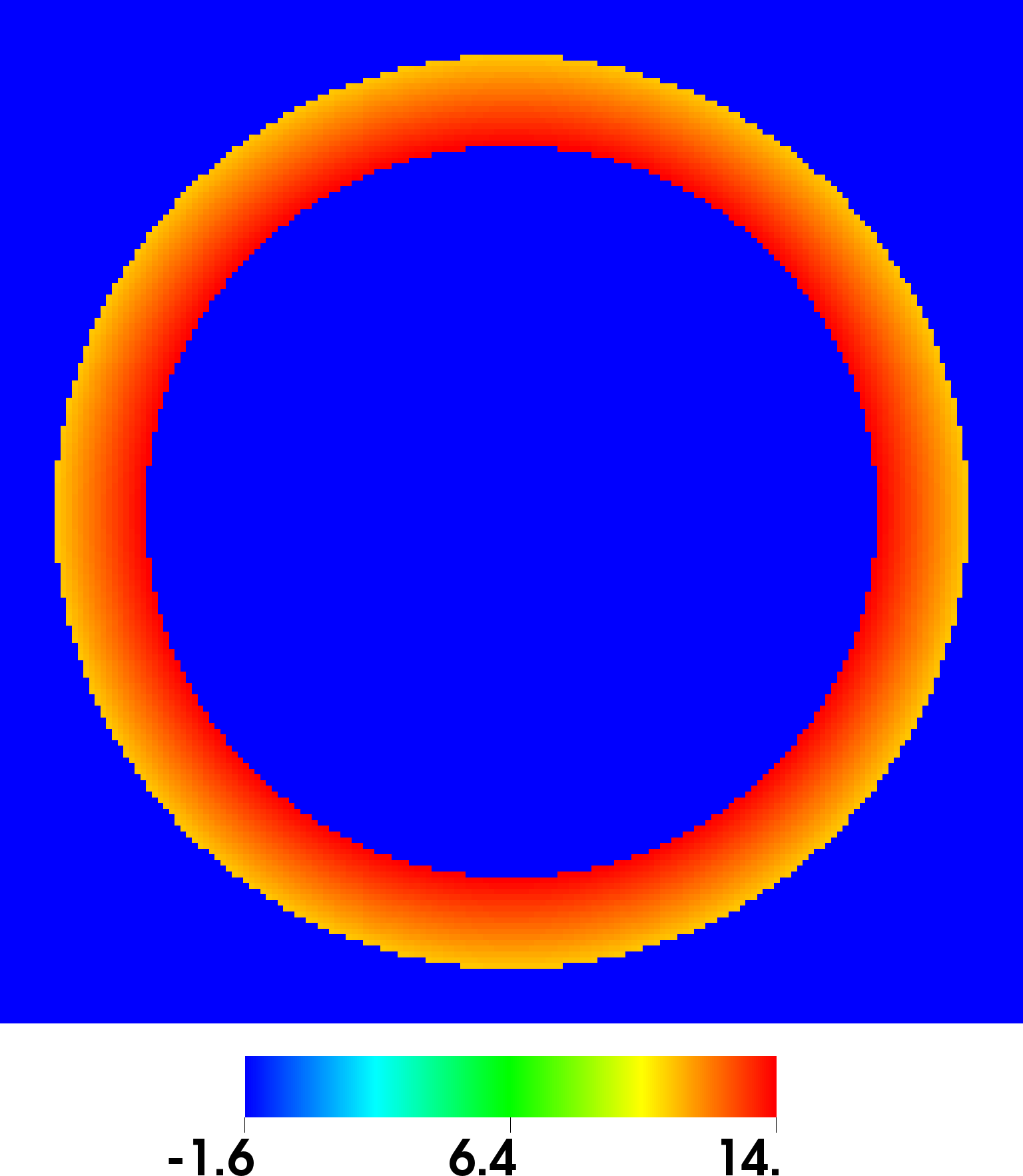}
\caption{The pressure--like fields for the sharp interface method, with the Cartesian grid composed of $128^2$ cells. The $\pi$ field is on the left and the $\varphi$ field is on the right.  The $\varphi$ field is computed using the steady state formulation.  The $\pi$ field is continuous at the fluid--structure interface.  Also, note that the $\varphi$ field is computed on the solid mesh, and is being projected onto the background Cartesian grid for output purposes and to compute errors.}
\label{fig:static1}
\end{center}
\end{figure}

For the numerical simulations, the fluid domain is taken to be $[0,L]^2$ with $L = 1$ mm. 
In this test we set $\mu_e = 1$ N/mm, $\rho = $1 kg/mm$^3$ and $\mu = $1 N$\cdot$s/mm$^2$.  
The velocity is set to zero at the boundary of the fluid domain, and the pressure field $p$ and pressure--like fields $\varphi$ and $\pi$ are normalized to have zero mean.  
\revthree{We use the spatial discretization parameter $h = L/N$ with $N = 2^{m}$ for $m = 5,6,7,8,$ and $9$.} The time step size is $\Delta t = 0.25 \times h$, we set $M_\text{fac} = 2$, and we perform simulations to a final time of $T = 0.01$ s, at which time we compute the errors on the Cartesian grid.  
Cartesian representations of the pressure--like fields $\pi$ and $\varphi$ computed using the sharp interface method are shown in Figure \ref{fig:static1}.  
Notice the continuity of $\pi$ at the fluid--structure interface, because the discontinuity in the physical pressure field $p$ is accounted for in the boundary condition for $\varphi$. 

\begin{figure}[h!]
\begin{center}
\includegraphics[scale=0.11,trim=0 0 -250 0]{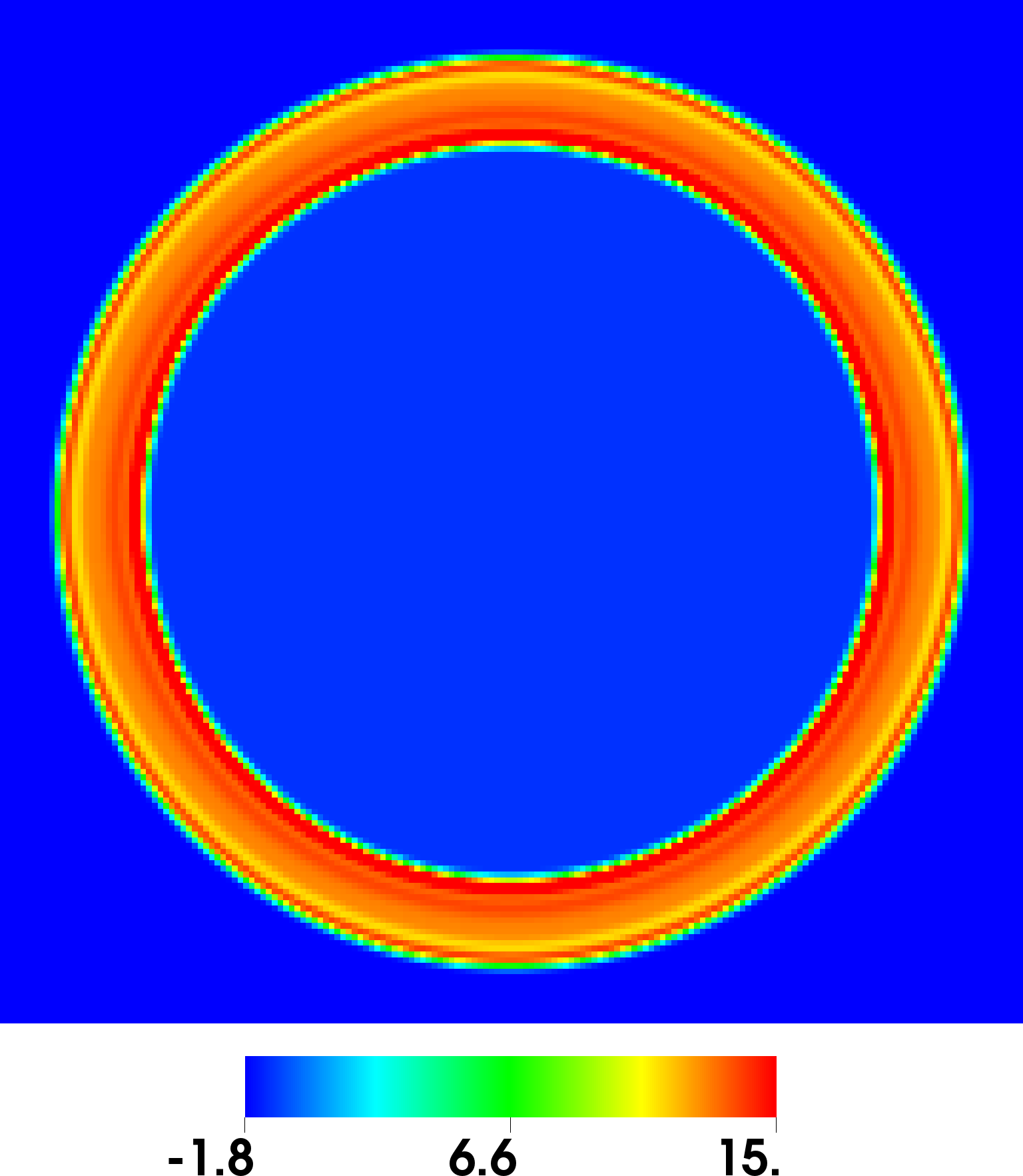}
\includegraphics[scale=0.11,trim=-250 0 0 0]{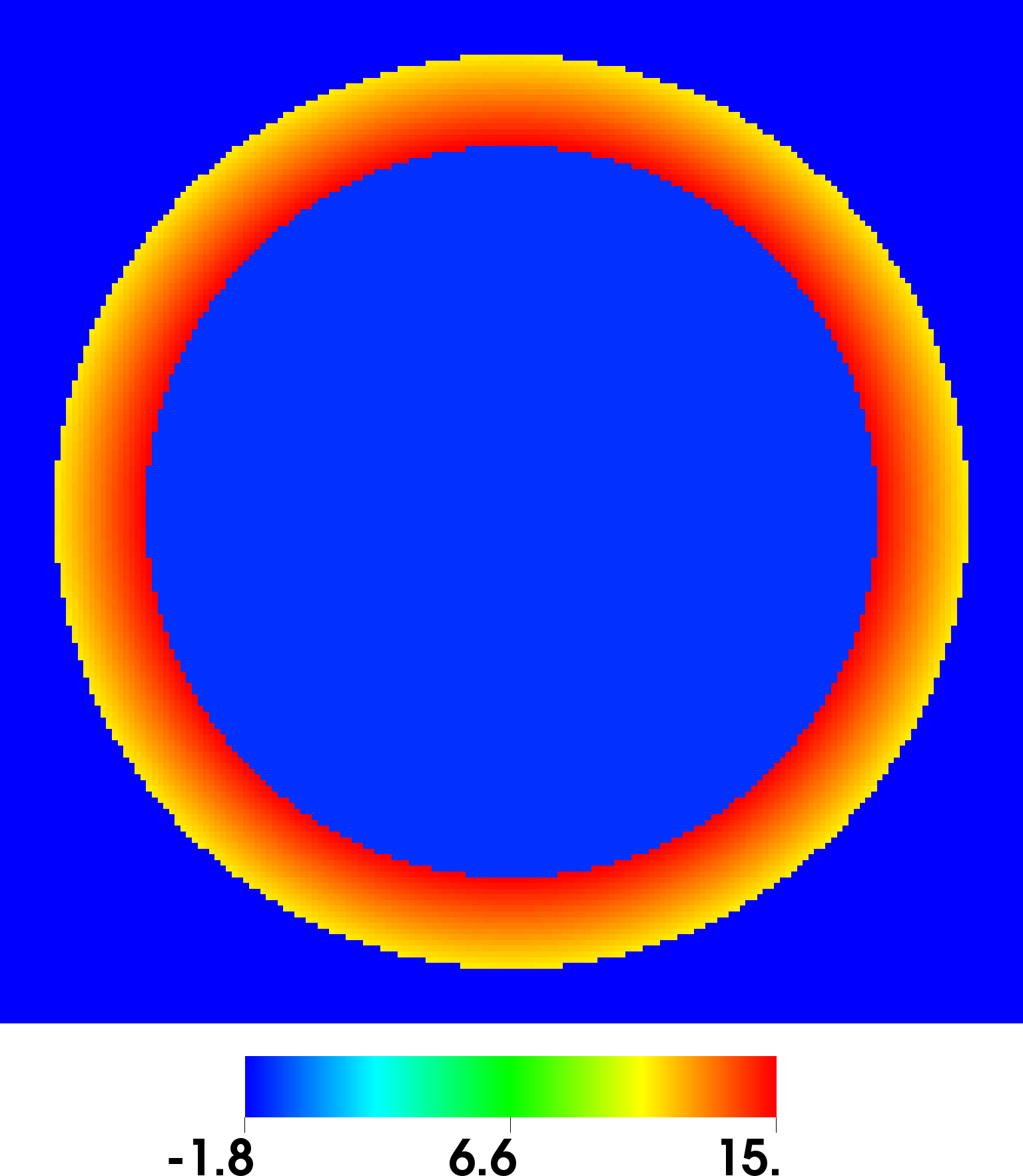}
\caption{On the left is the pressure field $p$ computed using the original IBFE method. On the right is the pressure field $p = \varphi + \pi$, with $\varphi$ computed using the steady state formulation.  The Cartesian grid contains $128^2$ cells.  There are artifacts at the fluid--structure interface in the results from the original IBFE method.  These artifacts appear to be eliminated in the results from the sharp interface method.}
\label{fig:static2}
\end{center}
\end{figure}

The pressure field computed with the original IBFE method, along with physical pressure field $p = \varphi + \pi$ from the sharp interface method, are shown in Figure \ref{fig:static2} on the left and right respectively.  
The pressure field determined by the original IBFE method contains artifacts at the fluid--structure interface, whereas the sharp interface approach cleanly resolves the pressure discontinuities.

\begin{figure}[h!]
\begin{center}
\includegraphics[scale=0.5,trim=0 0 -20 0]{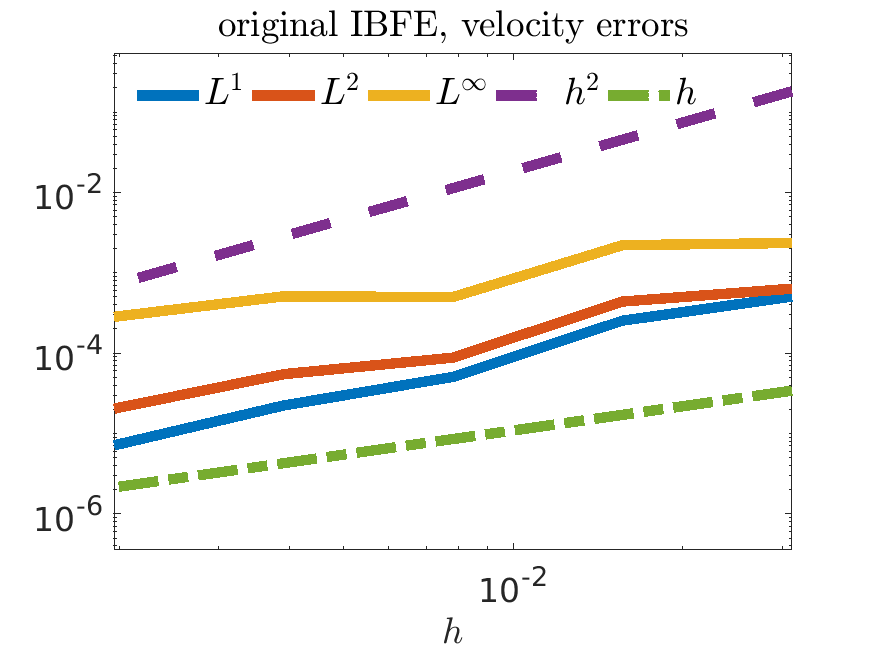}
\includegraphics[scale=0.5,trim=0 0 0 0]{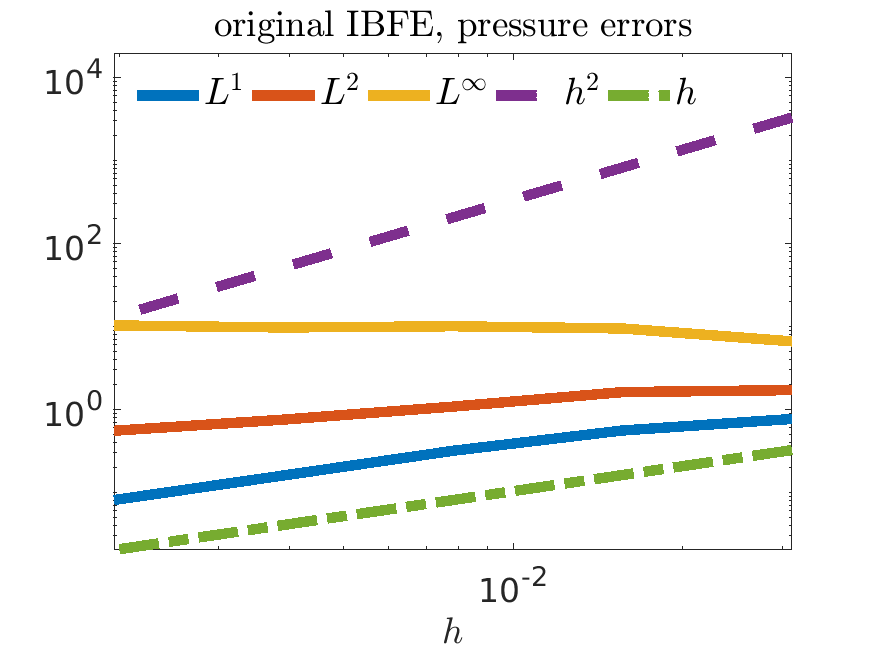}
\caption{Velocity and pressure absolute errors for static ring test with the original IBFE method. The pressure does not converge pointwise.}
\label{fig:static_errors0}
\end{center}
\end{figure}

The impact of the sharp interface method in this example can be quantified by considering the velocity and pressure errors.
The errors are computed on the Cartesian grid by projecting $\varphi$ onto the grid.  
This representation of $\varphi$ is visualized on the right in Figure \ref{fig:static1}. 
We remark that the errors displayed are absolute, to highlight the differences in absolute error between the original IBFE method and the sharp interface method.  
Figure \ref{fig:static_errors0} shows the errors for the original IBFE method.  
For velocity, we obtain approximately a rate of 2 for the $L^1$ error, a rate of 1.5 for the $L^2$ error, and a rate of 1 for the $L^\infty$ error.  For the pressure, we get a rate of 1 for the $L^1$ error, a rate of 0.5 for the $L^2$ error, and non--convergence for the $L^\infty$ error.  This is consistent with the results reported in \cite{Griffith17}.

\begin{figure}[h!]
\begin{center}
\includegraphics[scale=0.5,trim=0 0 -20 0]{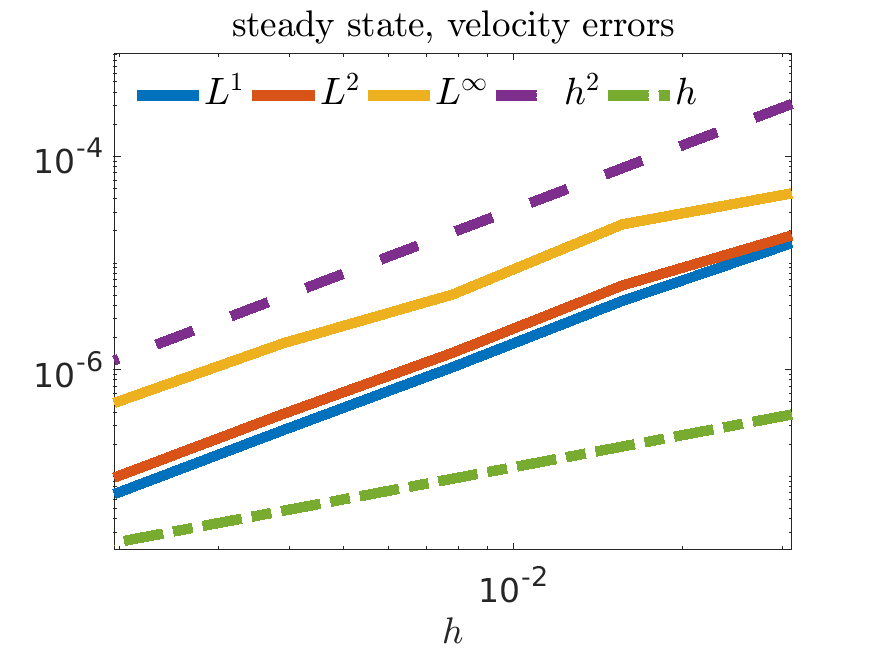}
\includegraphics[scale=0.5,trim=0 0 0 0]{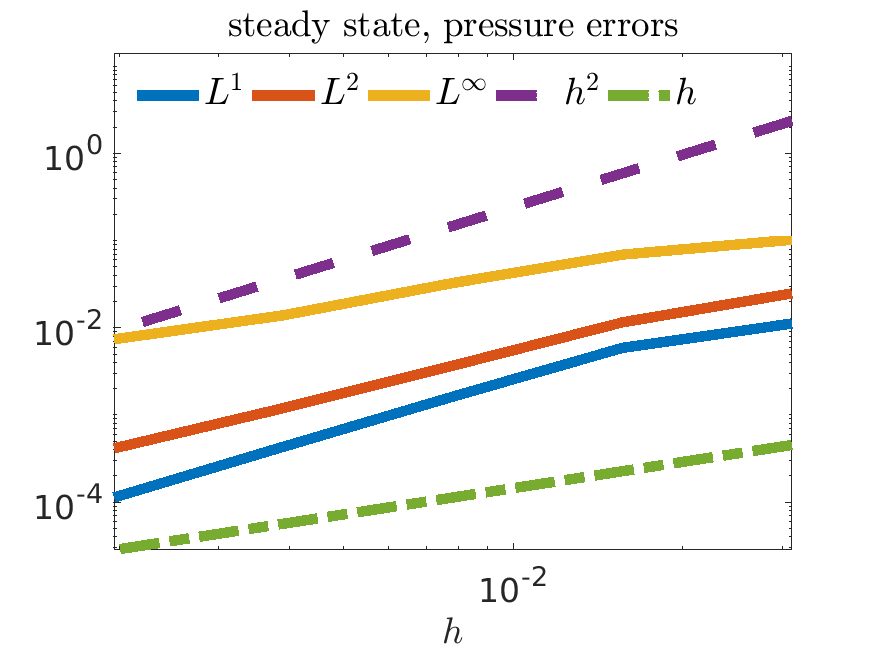}
\caption{Velocity and pressure absolute errors for the static ring test using the sharp interface method and the steady state formulation to determine $\varphi$.  The overall errors are several orders of magnitude smaller, compared to the results from the original IBFE method.  Also, the errors converge \purp{more rapidly under grid refinement}, and in particular the pressure converges pointwise.}
\label{fig:static_errors1} 
\end{center}
\end{figure}

\begin{figure}[h!]
\begin{center}
\includegraphics[scale=0.5,trim=0 0 -20 0]{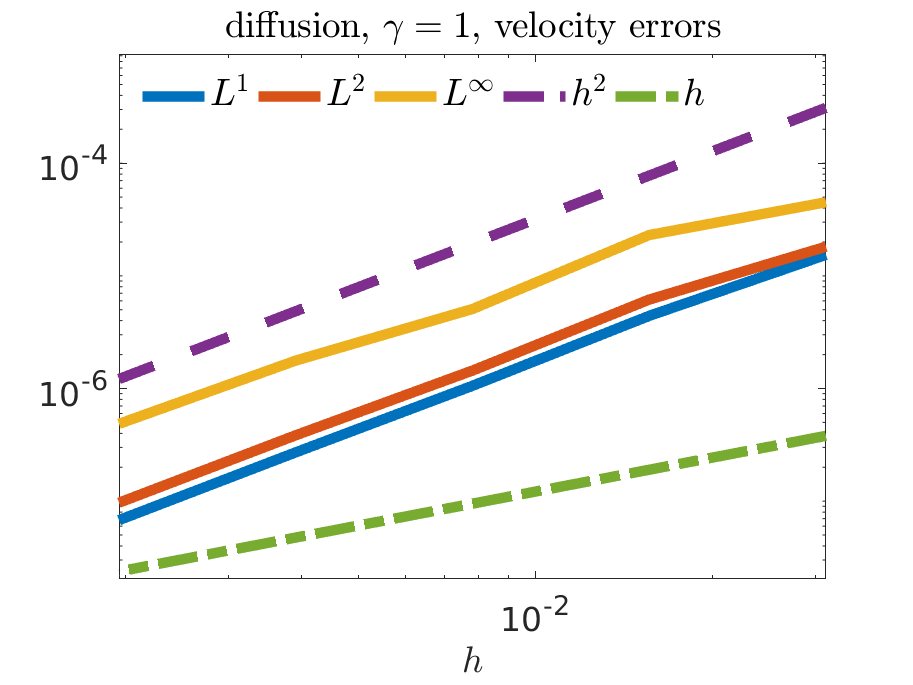}
\includegraphics[scale=0.5,trim=0 0 0 0]{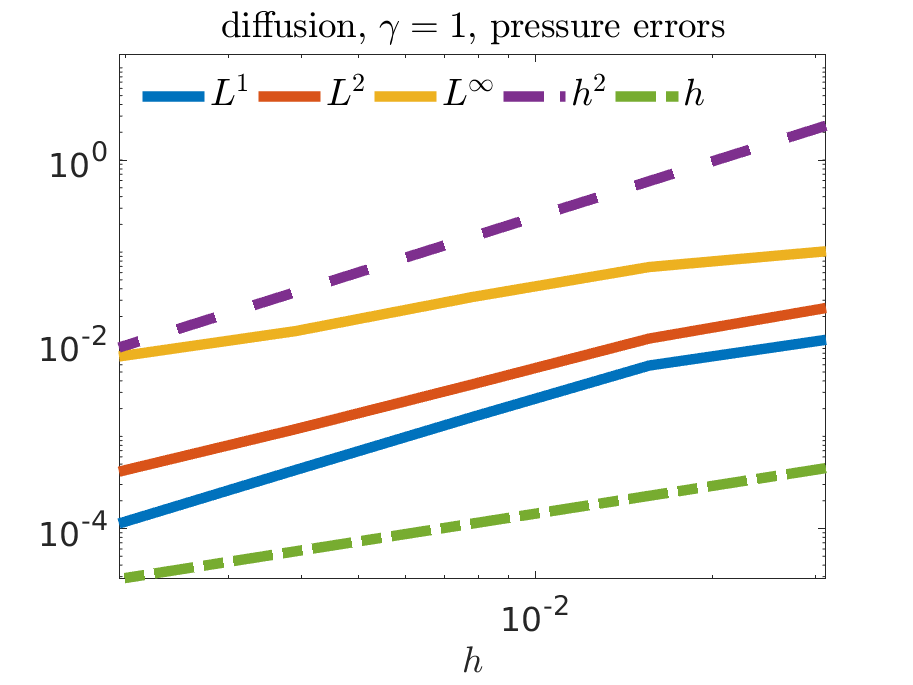}
\caption{Velocity and pressure absolute errors for static ring test using the sharp interface method and the diffusion formulation with $\gamma = 1$ to determine $\varphi$.  The errors are several orders of magnitude smaller, and error rates are improved, compared to the results from the original IBFE method.}
\label{fig:static_errors2}
\end{center}
\end{figure}

\begin{figure}[h!]
\begin{center}
\includegraphics[scale=0.5,trim=0 0 -20 0]{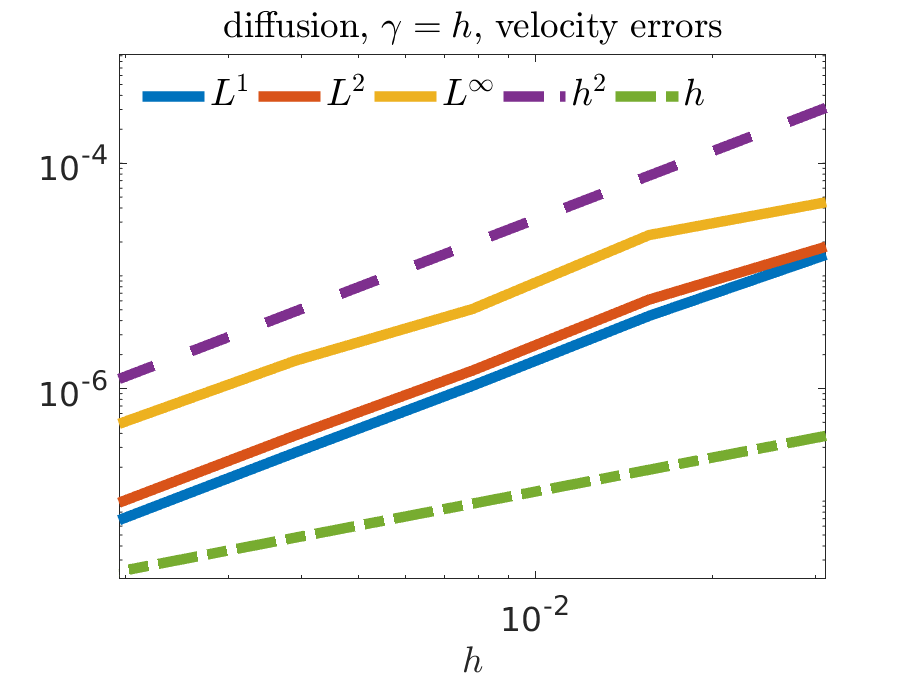}
\includegraphics[scale=0.5,trim=0 0 0 0]{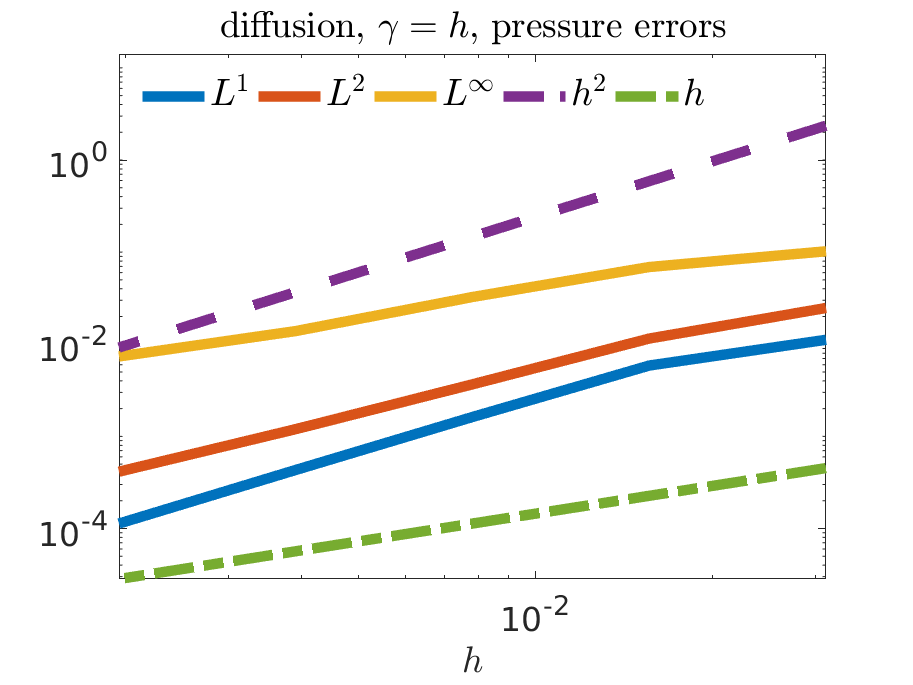}
\caption{Velocity and pressure absolute errors for static ring test using the sharp interface method and the diffusion formulation with $\gamma = h$ to determine $\varphi$. The errors are several orders of magnitude smaller, and error rates are improved, compared to the results from the original IBFE method.}
\label{fig:static_errors3}
\end{center}
\end{figure}

Figure \ref{fig:static_errors1}--\ref{fig:static_errors3} show the corresponding errors for the different formulations of the sharp interface method. 
Improved convergence rates are clearly apparent.  
Each of these figures corresponds to a different method for computing $\varphi$; Figure \ref{fig:static_errors1} uses the steady state formulation and Figures \ref{fig:static_errors2} and \ref{fig:static_errors3} use the diffusion formulation with $\gamma = 1$ and $\gamma = h$ respectively.  
In all cases, the velocity error converges to zero at a second order rate in each norm.  The $L^1$ error for the pressure converges with a rate of 2, the $L^2$ error has a rate of 1.5, and the $L^\infty$ error has a rate of 1.  
Perhaps most significantly, the sharp interface method yields an improvement in the errors of several orders of magnitude.

\begin{figure}[h!]
\begin{center}
\includegraphics[scale=0.5,trim=0 0 0 0]{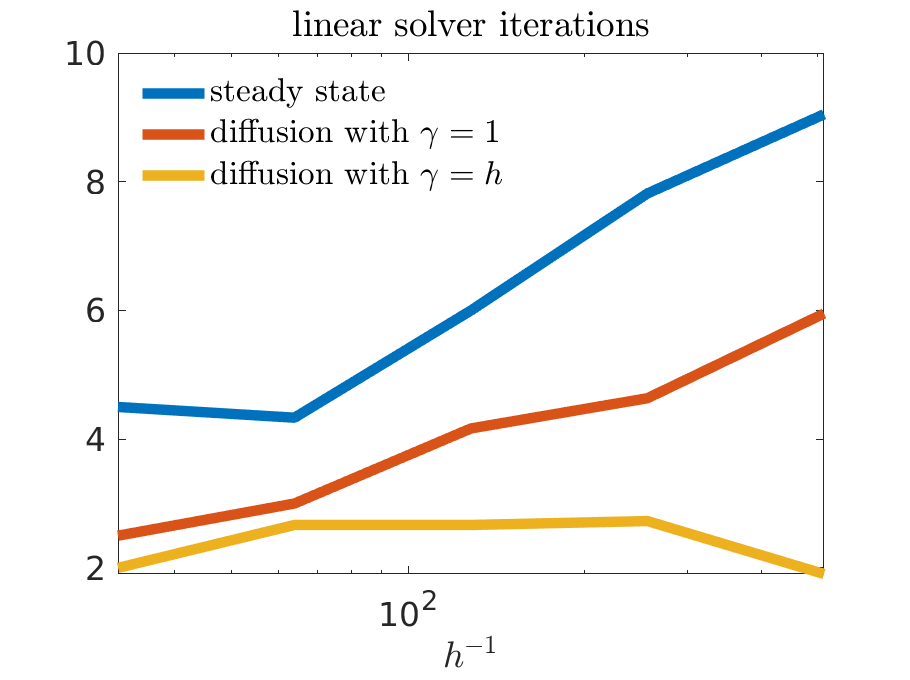}
\caption{Number of linear solver iterations for computing $\varphi$ in various ways for the static ring test.  The diffusion formulation is discretized in time using Crank--Nicolson, and the parameter $\gamma$ can be used to control the condition number of the resulting linear system.  The diffusion formulation for $\varphi$ leads to a milder scaling in the number of iterations as the number of degrees of freedom increase.  When $\gamma = h$, the number of iterations appears to be grid independent.}
\label{fig:num_iter}
\end{center}
\end{figure}

Figure \ref{fig:num_iter} plots the average number of linear solver iterations required to compute $\varphi$, as we vary $h$.  
We consider the steady state formulation for $\varphi$ as well as the diffusion formulation, with diffusion constant $\gamma$ set equal to either $1$ or $h$.  Notice that the number of iterations scale more mildly for cases that determine $\varphi$ with the diffusion formulation. If $\gamma = h$, we observe that the number of iterations is essentially independent of $h$.

\clearpage
\subsection{Thick inflating ring}

For this test, we consider the same ring as in the first example, but with a different constitutive model.  
In this case, the reference configuration is circular, i.e. 
\begin{align*}
(R,\Theta) \in U = [R_\text{in}, R_\text{out}] \times [0, 2\pi].
\end{align*}
\purp{In particular, we note the reference coordinates are the initial coordinates.} The material model is defined with respect to the polar deformation gradient $\mathbb{F}_\text{polar}$ and polar first Piola--Kirchoff stress $\mathbb{P}^\text{e}_\text{polar}$ as follows:
\begin{align*}
\mathbb{P}_\text{polar}^\text{e} = \mu_e \left(\mathbb{F}_\text{polar} - \mathbb{F}_\text{polar}^{-T} \right).
\end{align*}
A volume of fluid $A_\text{add}$ is injected gradually in the center of the ring, and the pressure is held at zero outside of the ring. 
At steady state, the exact solution for the pressure field is
\begin{align*}
p(r) = 
\begin{cases}
 \frac{\mu_e A_\text{add}}{2\pi}\left( \frac{1}{r_\text{in}^2} - \frac{1}{r_\text{out}^2} \right) , \quad  r \leq r_\text{in} \\
  -\frac{\mu_e A_\text{add}}{2\pi}\left( \frac{1}{r^2} + \frac{1}{r_\text{out}^2} \right) , \quad r_\text{in} \leq r \leq r_\text{out} \\
0, \quad r > r_\text{out} 
\end{cases}
\end{align*}
where $r_\text{in} = r(R_\text{in})$ and $r_\text{out} = r(R_\text{out})$ are given as
\begin{align*}
r(R) = \left(R^2 + \frac{A_\text{add}}{\pi} \right)^{1/2}, \quad R_\text{in} \leq R \leq R_\text{out}.
\end{align*} 
A derivation can be found in the appendix, along with an explanation of the notation.  
We remark that the pressure discontinuities can be made to be very large in this example, by increasing the value of $\mu_e$.  
The pressure in the center of the ring is takes a positive value, but within the solid, the pressure takes large negative values for a large choices of $\mu_e$.  
In this test we set $\mu_e = 10^4$ N/mm$^2$, $\rho = $1 kg/mm$^3$ and $\mu = $1 N$\cdot$s/mm$^2$. 

For our simulations, the fluid volume $A_\text{add}$ is added to the ring center over the first 0.1 s, and we run the model to a final time of 1 s to reach an approximate steady state. 
A combination of normal traction and tangential no--slip velocity boundary conditions are applied to the computational domain $\Omega$ to hold the pressure at zero outside of the ring.  
The computational domain $\Omega$ is $[-L,L]^2$ with $L = 1$ mm.  
The initial inner and outer radii of the ring are $R_\text{in} = 0.25$ mm and $R_\text{out} = 0.3125$ mm respectively.
\revthree{The spatial discretization parameter is $h = L/N$ with $N = 2^{m}$ for $m = 5,6,7,8,$ and $9$.}
The time step size is $\Delta t = 0.025 \times h$, and we set $M_\text{fac} = 1$.

\begin{figure}[h!]
\begin{center}
\includegraphics[scale=0.11,trim=100 0 -100 0]{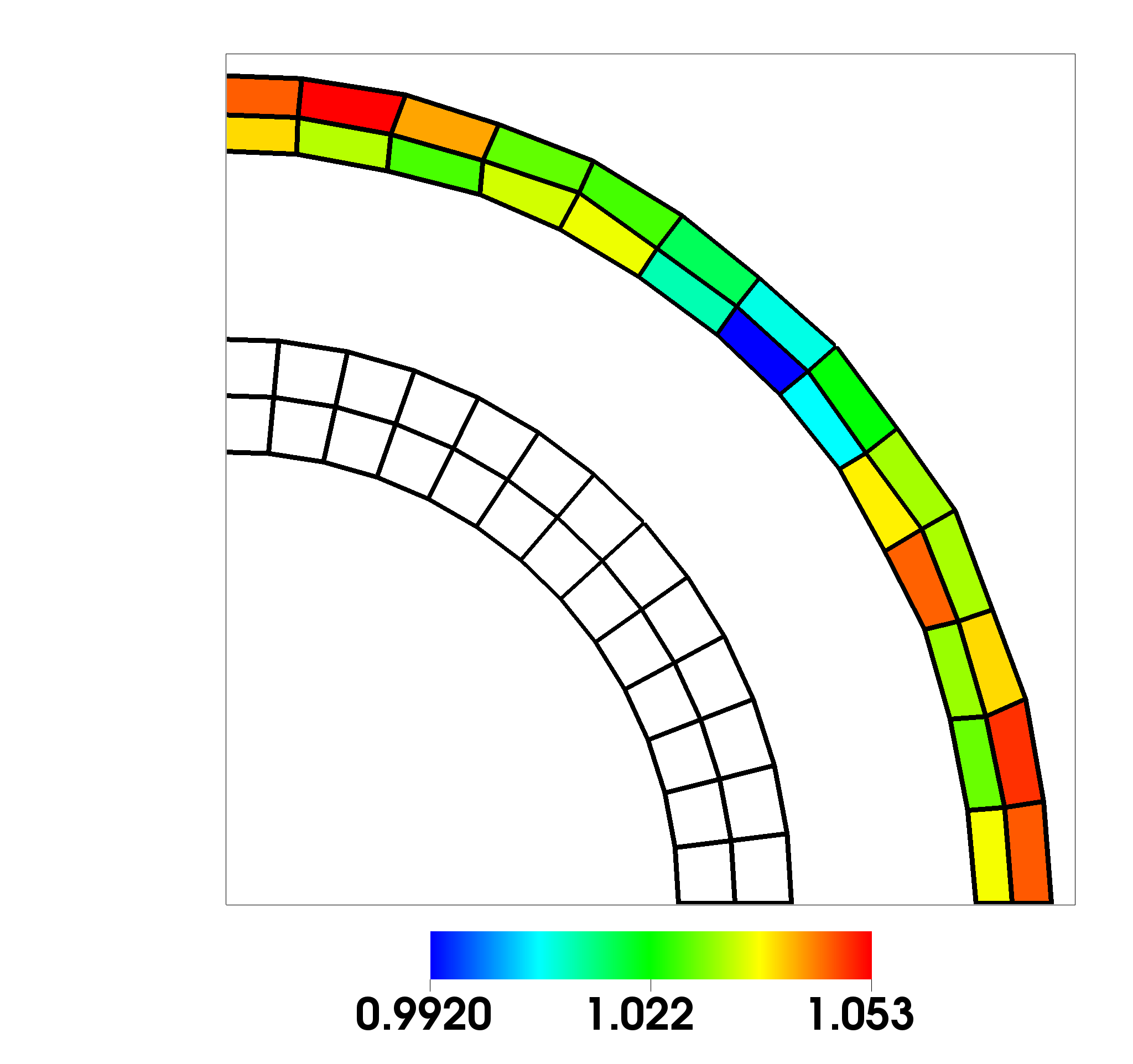}
\includegraphics[scale=0.11,trim=-100 0 0 0]{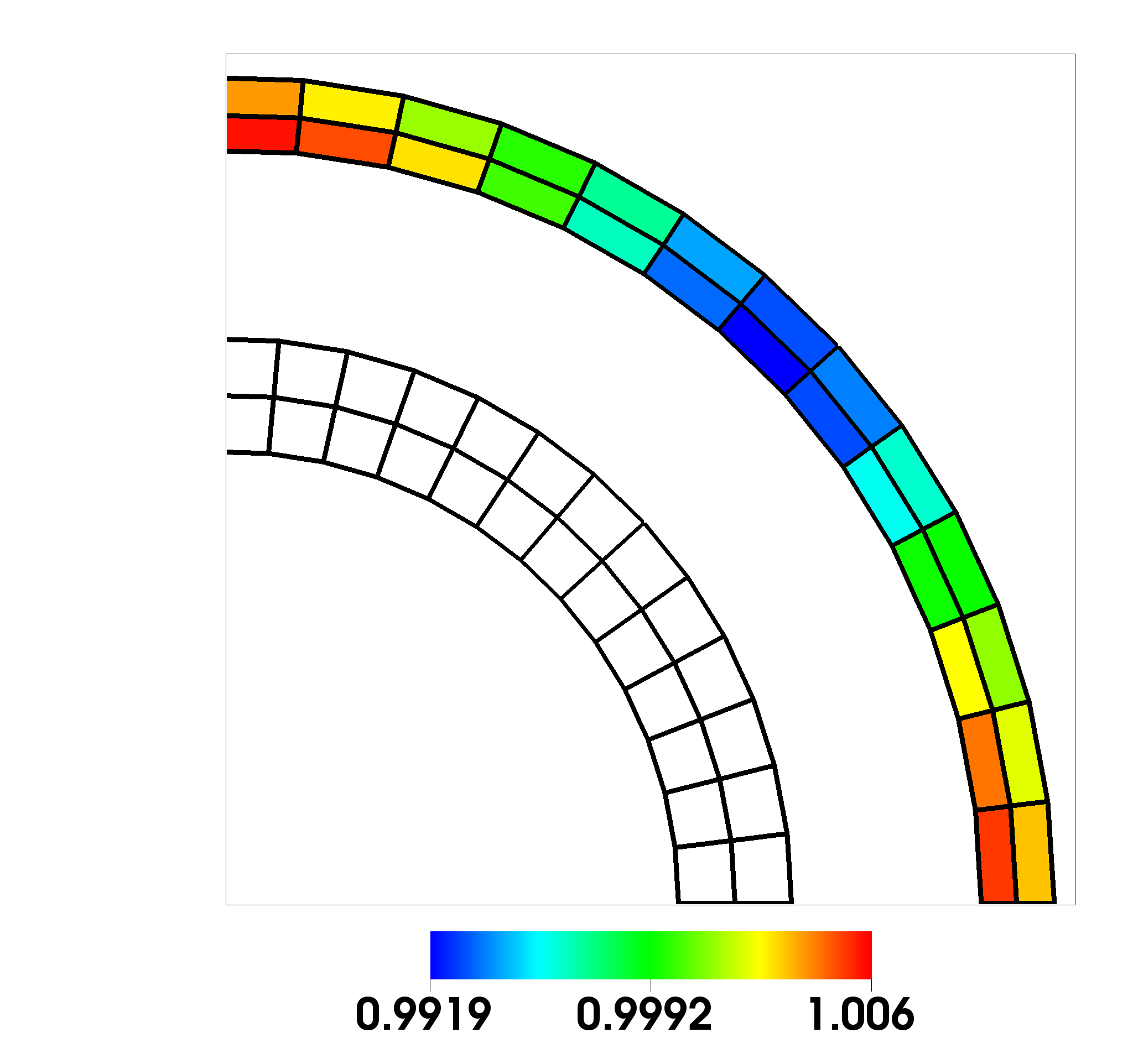}
\caption{A visualization of the initial and final finite element mesh configurations for the inflating ring test, with the original IBFE method on the left and the sharp interface method on the right. The steady state formulation is used to determine $\varphi$.  
The coloring in the final configuration corresponds to the deteminant $J = \det(\mathbb{F})$. \purp{Note that the color bars are different between the panels.} The Cartesian grid is composed of $32^2$ cells.
The deformations in the original IBFE method are nonuniform and less incompressible, compared to the results from the sharp interface method. }
\label{fig:inflate1}
\end{center}
\end{figure}

Figure \ref{fig:inflate1} shows the nontrivial displacement of the solid finite element from its initial to final configuration, with the original IBFE method on the left and the sharp interface method on the right.
In this figure, we use the steady state formulation for $\varphi$.  
Because the problem is axially symmetric, we show only one quarter of the ring. 
The final configuration is colored, indicating the elementwise value of $J = \det(\mathbb{F})$.   
In this case, the Eulerian grid is very coarse and contains $32^2$ fluid cells; because $M_\text{fac} = 1$, the initial finite element mesh indicates the approximate size of the fluid cells.  
Inflation of the ring results in a deformed configuration in which the ring thickness is close to the size of a single Cartesian grid cell.  The deformations in the original IBFE method deviate from incompressible deformations (i.e., corresponding to $J = 1$), perhaps not surprisingly because the Cartesian grid is coarse with respect to the thickness of the ring in the final configuration.  
In contrast, the deformations determined by the sharp interface method are much closer to being incompressible, highlighting its utility in numerically conserving volume, at least for coarse finite element meshes.

%\begin{figure}[h!]
%\begin{center}
%\includegraphics[scale=0.225,trim=0 100 150 0]{IMAGES/inflate_shell_initial}
%\includegraphics[scale=0.225,trim=0 100 0 0]{IMAGES/inflate_shell_final}
%\caption{Initial (left, blue) and final (right, red) configurations for the inflating thick shell mesh.  The Cartesian fluid mesh is also displayed in the background, in grey.}
%\end{center}
%\end{figure}

\begin{figure}[h!]
\begin{center}
\includegraphics[scale=0.11,trim=0 0 -250 0]{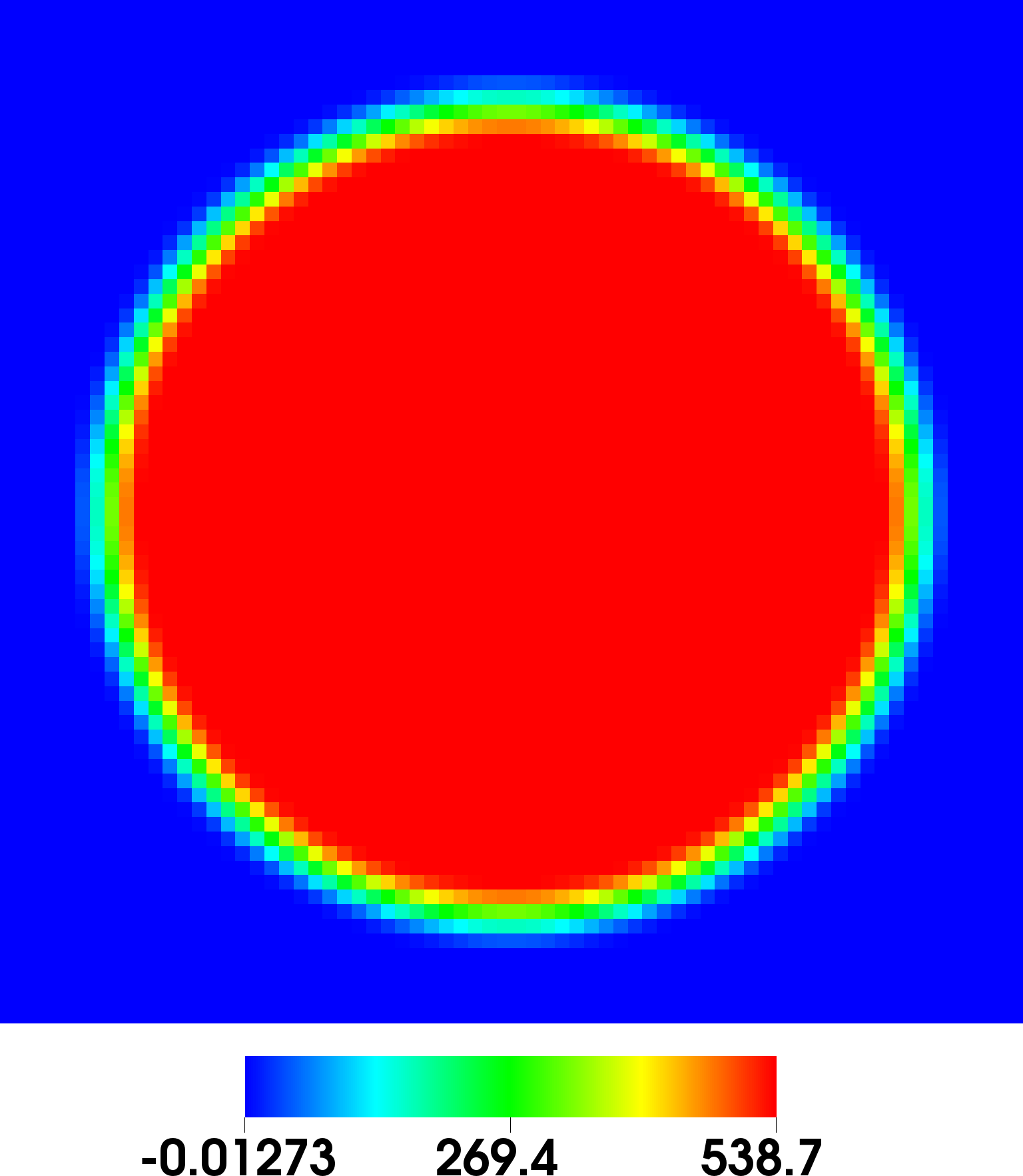}
\includegraphics[scale=0.11,trim=-250 0 0 0]{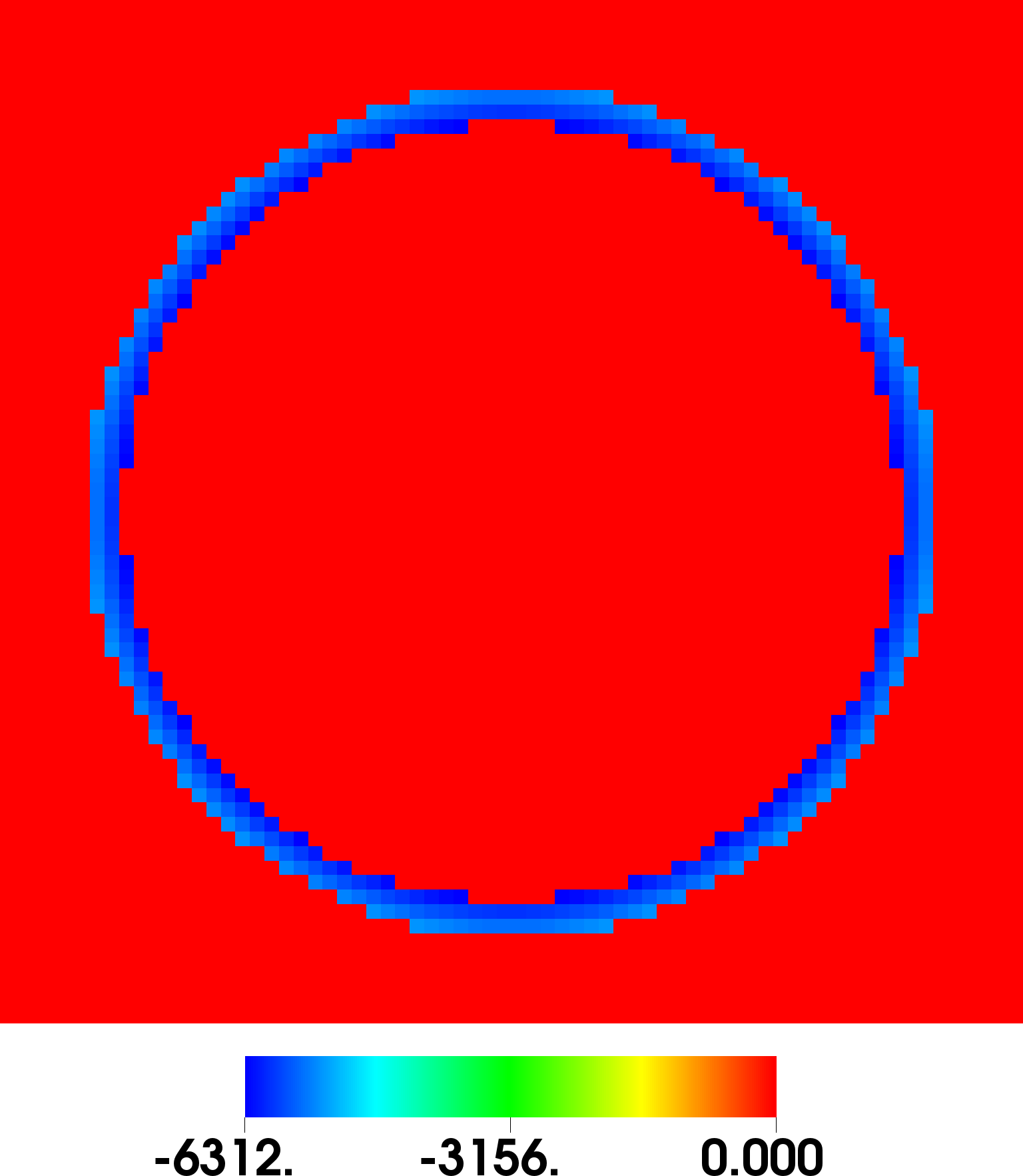}
\caption{On the left is the $\pi$ field and on the right is the $\varphi$ field, computed with the steady state formulation for the inflating ring test.  The $\pi$ field is continuous at the fluid--structure interface, and the $\varphi$ field contains the jump in the physical pressure field.  
The Cartesian grid composed of $128^2$ cells.}
\label{fig:inflate2}
\end{center}
\end{figure}

\begin{figure}[h!]
\begin{center}
\includegraphics[scale=0.11,trim=0 0 -250 0]{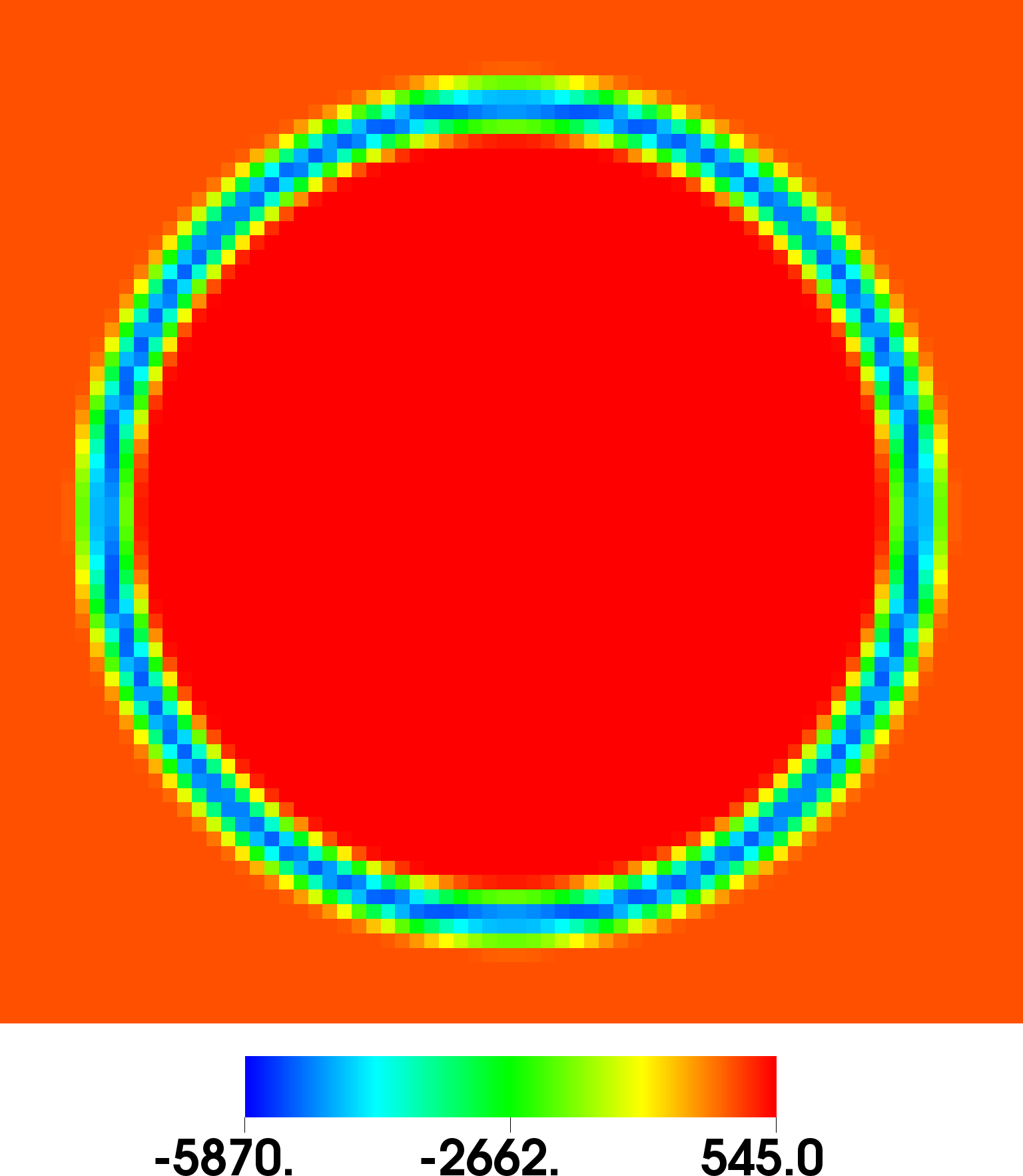}
\includegraphics[scale=0.11,trim=-250 0 0 0]{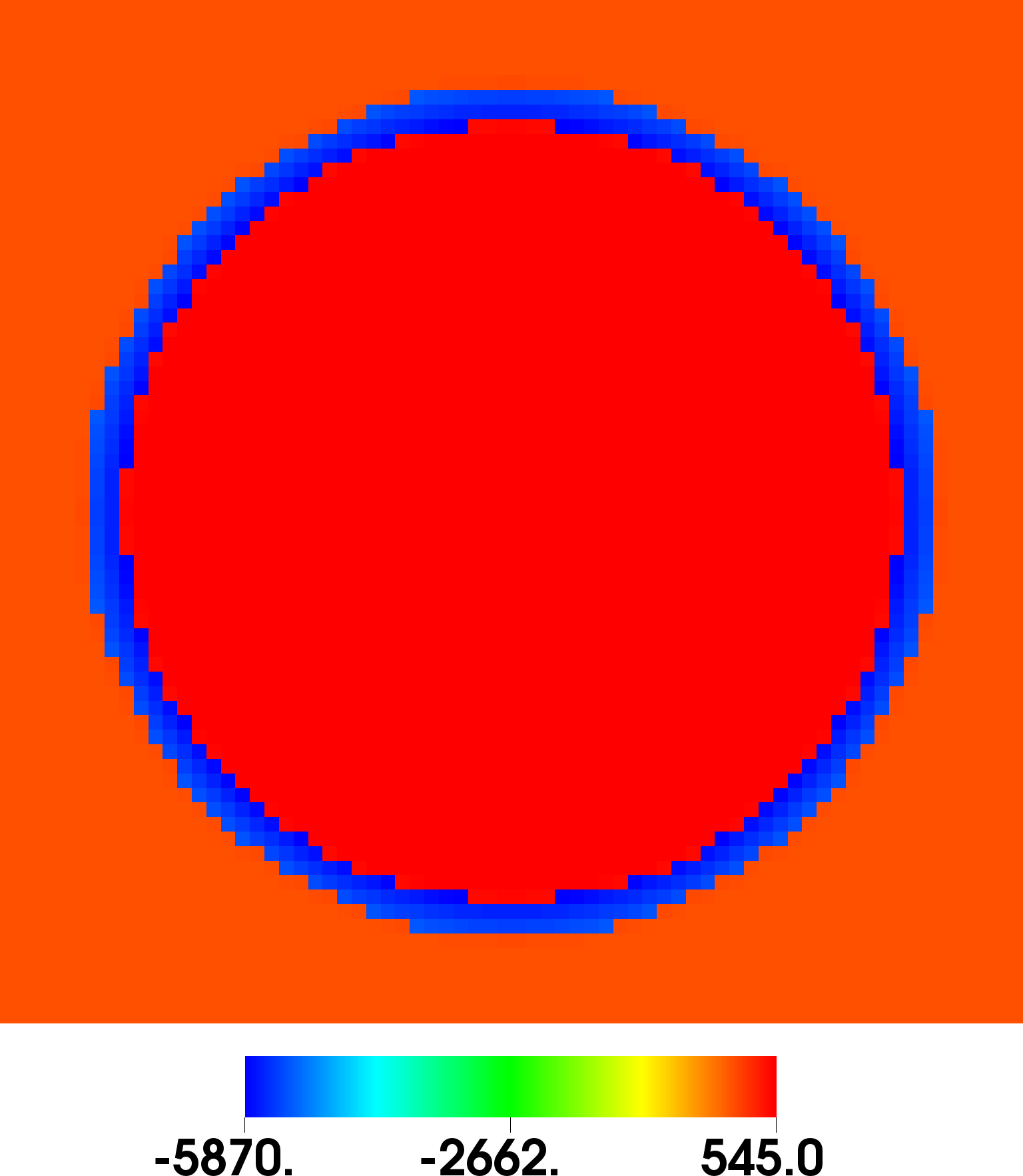}
\caption{On the left is the pressure field computed for the inflating ring test with the original IBFE method, and on the right is the pressure field $p = \pi + \varphi$ computed with the sharp interface method. The steady state formulation is used to solve for $\varphi$, and the Cartesian grid contains $128^2$ cells.  The results from the original IBFE method have substantial artifacts at the fluid--structure interface.  These artifacts appear to be eliminated by the sharp interface method.}
\label{fig:inflate3}
\end{center}
\end{figure}

Figure \ref{fig:inflate2} shows the pressure--like fields $\pi$ and $\varphi$ on the Cartesian grid.  
As in the static ring case, $\pi$ is continuous at the fluid--structure interface, and the pressure discontinuity is completely accounted for by the boundary conditions imposed in the $\varphi$ field.  
Figure \ref{fig:inflate3} shows the physical pressure $p$. 
Results from the original IBFE method are shown on the left and the physical pressure $p = \pi + \varphi$ determined by the sharp interface method are shown on the right.  
The original IBFE method clearly smears out the pressure discontinuity, whereas the discontinuity is cleanly resolved in the sharp interface method.

\begin{figure}[h!]
\begin{center}
\includegraphics[scale=0.5,trim=0 0 -50 0]{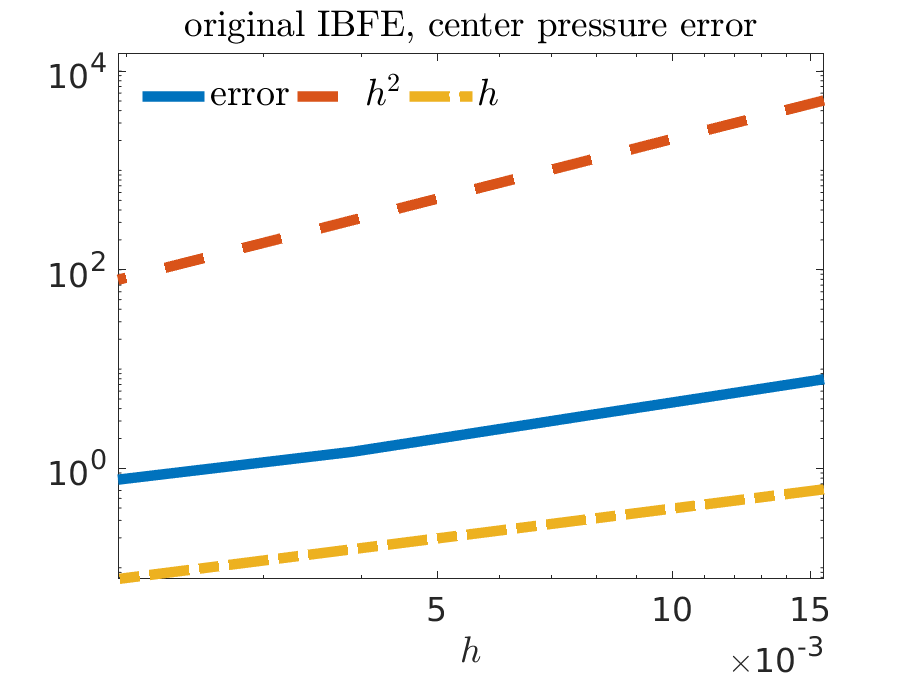}
\includegraphics[scale=0.5]{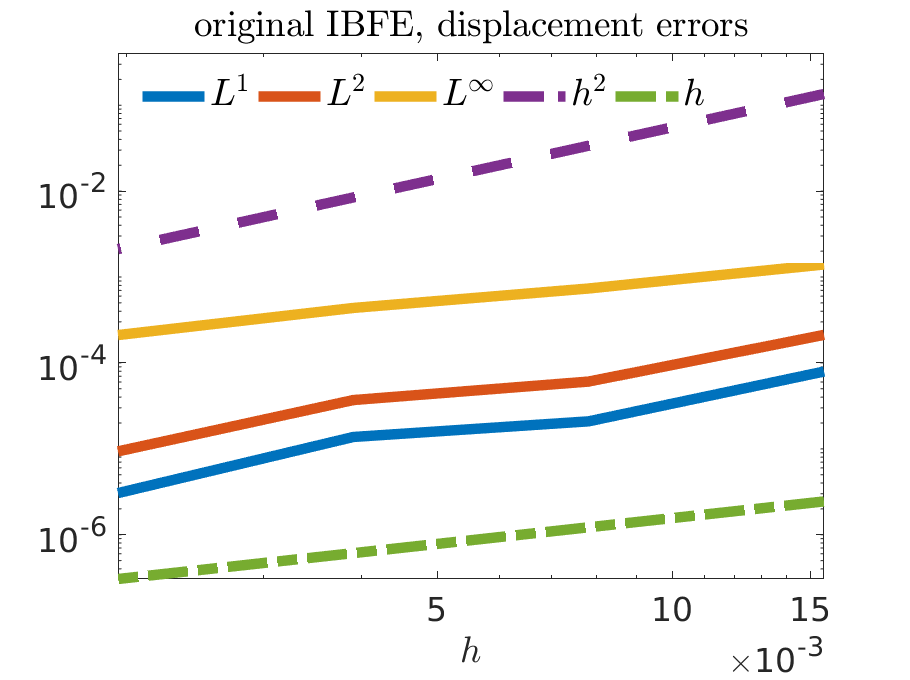} \\
\includegraphics[scale=0.5,trim=0 -30 -50 -30]{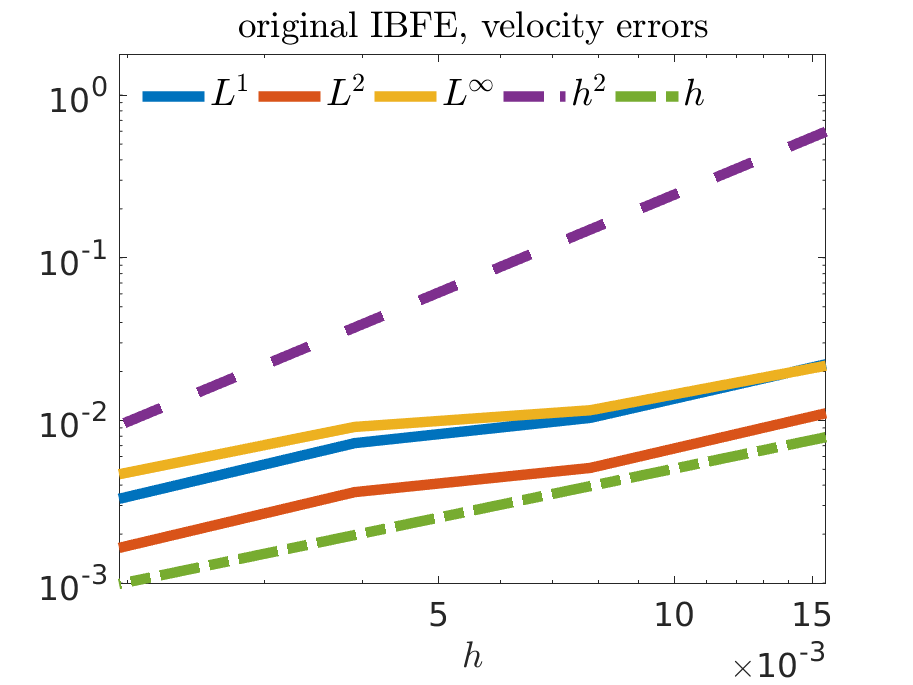}
\includegraphics[scale=0.5,trim=0 -30 0 -30]{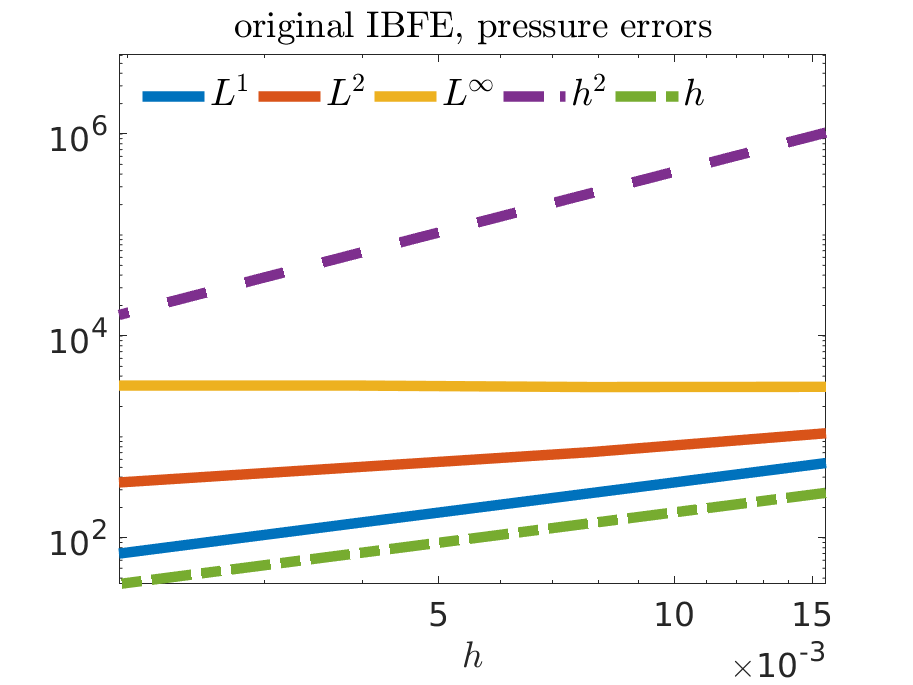} \\
\includegraphics[scale=0.5]{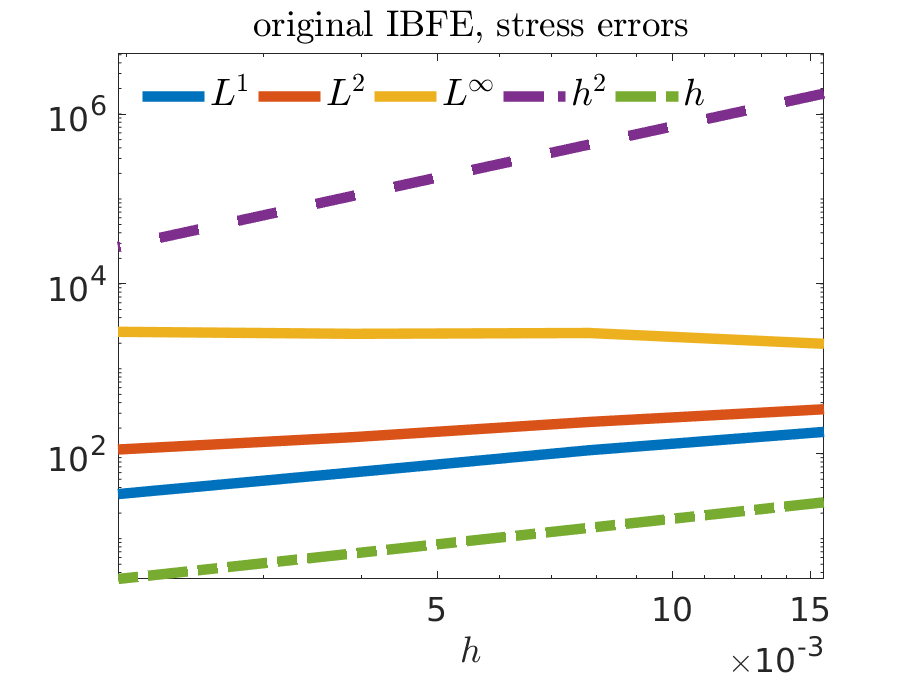}
\caption{Velocity, pressure, displacement, and elastic stress absolute errors for the inflating ring test with the original IBFE method.  Results are similar to the static ring example.  There is no pointwise convergence of the pressure and elastic stress.}
\label{fig:inflate_errors1}
\end{center}
\end{figure}

\begin{figure}[h!]
\begin{center}
\includegraphics[scale=0.5,trim=0 0 -50 0]{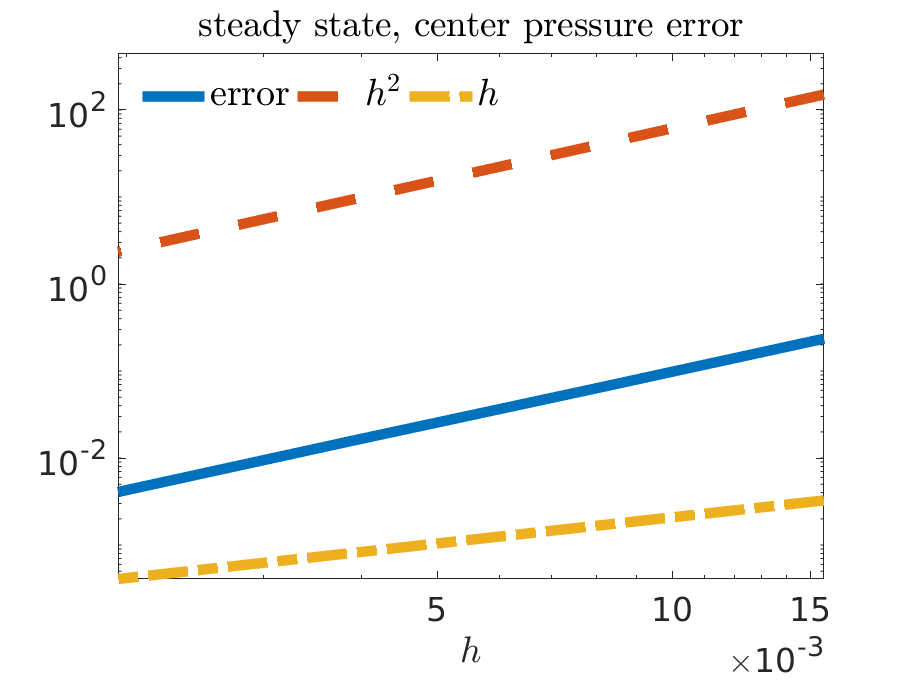}
\includegraphics[scale=0.5]{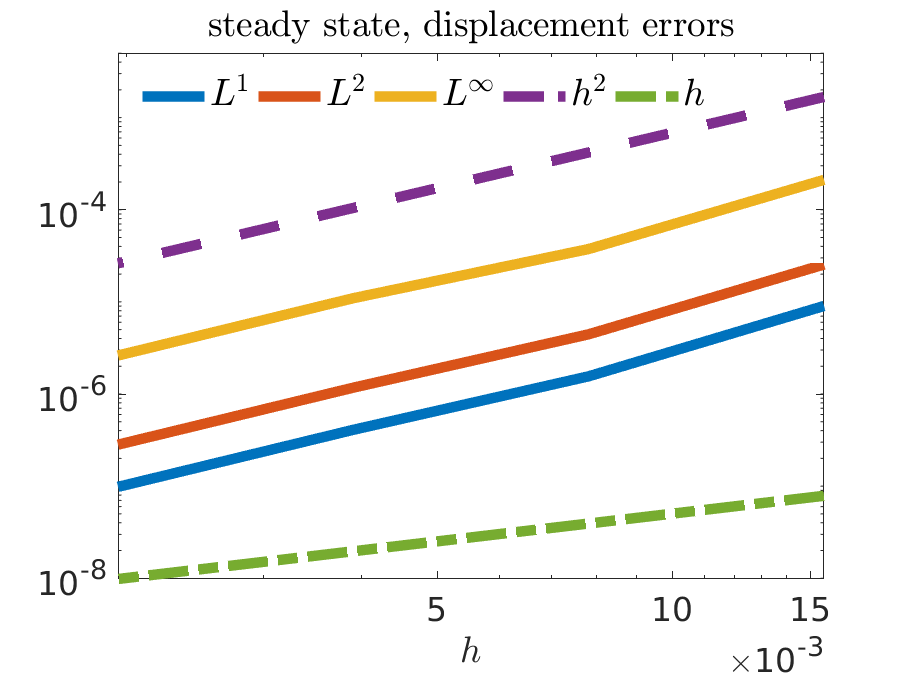} \\
\includegraphics[scale=0.5,trim=0 -30 -50 -30]{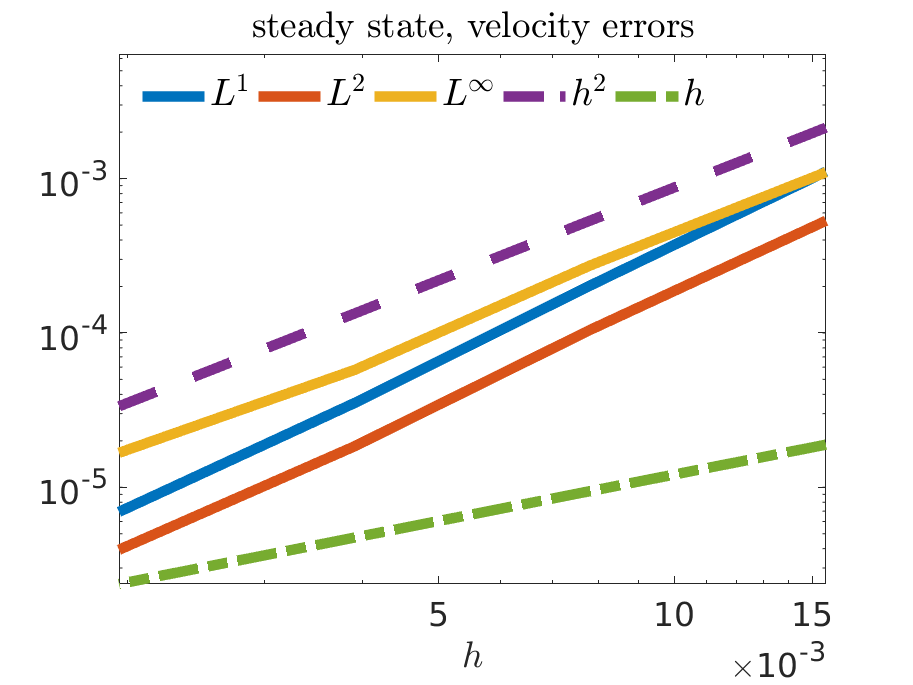}
\includegraphics[scale=0.5,trim=0 -30 0 -30]{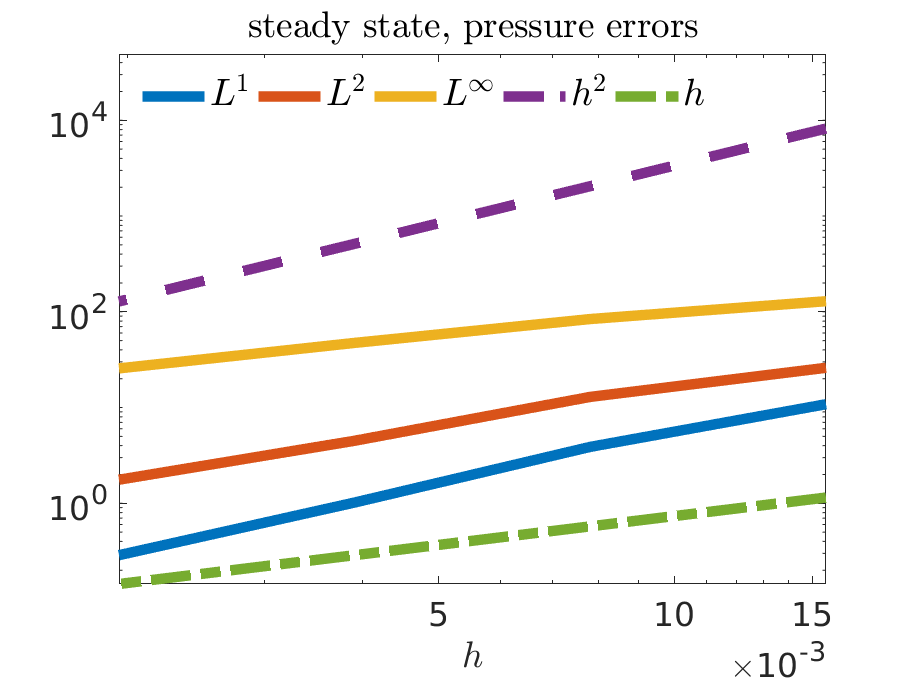} \\
\includegraphics[scale=0.5]{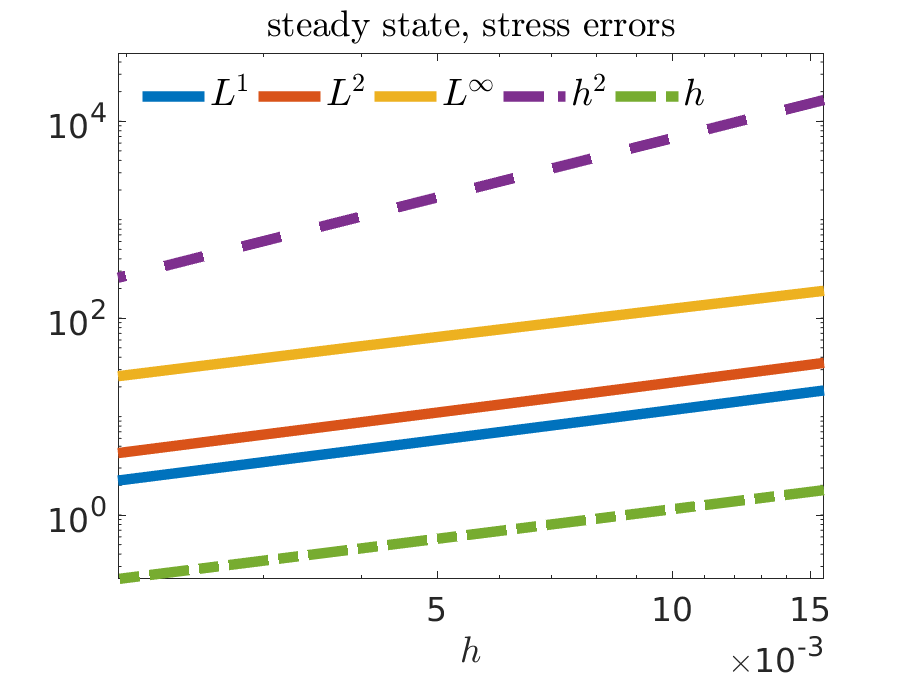}
\caption{Velocity, pressure, displacement, and elastic stress absolute errors for the inflating ring test with the sharp interface method, using the steady state formulation to determine $\varphi$.  Absolute errors and rates are much improved, compared to the results from the original IBFE method.  In particular, the pressure and elastic stress converge pointwise.}
\label{fig:inflate_errors2}
\end{center}
\end{figure}

\revtwo{Figures \ref{fig:inflate_errors1} and \ref{fig:inflate_errors2}} examine the absolute errors in the pressure, velocity, displacement, and elastic stress.  
The pressure at the center of the ring is sampled using a cosine kernel with radius equal to 0.1 mm. 
Figure \ref{fig:inflate_errors1} reports results for the original IBFE method. 
Results for the sharp interface method are shown in \revtwo{Figure \ref{fig:inflate_errors2}, with the steady state formulation used compute $\varphi$.}

In Figure \ref{fig:inflate_errors1}, for the original IBFE method, we observe the same convergence rates for the velocity and pressure as in the static ring test. 
The velocity errors converge with a rate of 2 in $L^1$, a rate of 1.5 in $L^2$, and a rate of 1 in $L^\infty$, and the pressure errors converge with a rate of 1 in $L^1$ and a rate of 0.5 in $L^2$.  
With the conventional method, there is no pointwise convergence of the overall pressure field, and the pressure sampled at the center of the ring converges with a rate of 1.  The displacement errors converge with the same rates as the corresponding velocity errors.  
The elastic stress errors are effectively the errors in the gradient of the displacement, which yields no pointwise convergence of the elastic stress error.

%\begin{figure}[h!]
%\begin{center}
%\includegraphics[scale=0.5,trim=0 0 -50 0]{IMAGES/ipdg_inflate_ring_center-pressure-errors}
%\includegraphics[scale=0.5]{IMAGES/ipdg_inflate_ring_displacement-errors} \\
%\includegraphics[scale=0.5,trim=0 -30 -50 -30]{IMAGES/ipdg_inflate_ring_velocity-errors}
%\includegraphics[scale=0.5,trim=0 -30 0 -30]{IMAGES/ipdg_inflate_ring_pressure-errors} \\
%\includegraphics[scale=0.5]{IMAGES/ipdg_inflate_ring_stress-errors}
%\caption{Velocity, pressure, and elastic stress errors for the inflating ring and the scheme with splitting, using IPDG.}
%\end{center}
%\end{figure}

\revtwo{Figures \ref{fig:inflate_errors2} indicates} that the sharp interface method gives improved rates of convergence as well as much smaller absolute errors.   All velocity and displacement errors converge with a rate of 2,  and correspondingly the elastic stress errors converge with a rate of 1.  
The error rates for the overall pressure field are improved: the $L^1$ error converges with a rate of 2, the $L^2$ error converges with a rate of 1.5, and the $L^\infty$ error converges with a rate of 1. 
The pressure sampled at the center of the ring converges with a rate of 2.  \revtwo{The results for the case where $\varphi$ is computed with the diffusion formulation with $\gamma = 1$ are the same.  When $\gamma = h$, the results are essentially the same except for the velocity errors.} 
In this case, the velocity errors appear to converge with a rate of 1 in all norms.

\clearpage 
\subsection{Compressed block}

Our third example is related to a standard solid mechanics test problem: compression of a two--dimensional block of isotropic material \cite{Reese99}.  
The strain energy density is neo--Hookean and is split into isochoric $W^\text{iso}$ and dilational $W^\text{dil}$ parts.  
This decomposition is achieved by using the first invariant $I_1 = \text{tr}(\mathbb{F}^T\mathbb{F})$, which is modified by scaling the deformation gradient so that its determinant is 1:
\begin{align*}
\bar{I}_1 =I_1(\bar{\mathbb{F}}^T \bar{\mathbb{F}}), \quad \bar{\mathbb{F}} = J^{-1/3} \mathbb{F}.
\end{align*}
Then the strain energy density is defined by
\begin{align*}
W = \frac{\mu_e}{2}(\bar{I}_1 - 3) + \frac{\lambda}{2} (\log J)^2 = W^\text{iso}(\bar{I}_1) + W^\text{dil}(J),
\end{align*}  
for shear modulus $\mu_e$ and bulk modulus $\lambda = \lambda(\nu)$, taken to be a function of a numerical Poisson ratio $\nu$ which is detailed in \cite{Vadala18}.  
The dilational energy $W^\text{dil}$ vanishes at the continuous level because $J = 1$.  
In the discretized equations, this term is not necessarily equal to zero, and it can be interpreted as a stabilization that enforces discrete incompressibility \cite{Vadala18}.  
The first Piola--Kirchoff stress takes the form:
\begin{align*}
\mathbb{P}^e = \mu_e J^{-2/3}\left(\mathbb{F} - \frac{I_1}{3} \, \mathbb{F}^{-T} \right) + \lambda(\nu) \log(J)\, \mathbb{F}^{-T},
\end{align*}
with $\lambda$ related to $\nu$ by:
\begin{align*}
\lambda(\nu) = \frac{2\mu_e(1+\nu)}{3(1-2\nu)}.
\end{align*}
As suggested in Ref. \cite{Vadala18}, we consider values for the numerical Poisson ratio $\nu =-1, 0, 0.4$.

In this test, we include surface force densities on the boundary of the block to weakly impose Dirichlet boundary conditions for the displacement.  First, define the following surface force density which tethers points to their reference configuration:
\begin{align*}
{\bs F}_\text{surface}({\bs X},t) = \kappa \left( {\bs X} - {\bs \chi}({\bs X},t) \right).
\end{align*}
Let ${\bs e}_1$ and ${\bs e}_2$ denote the standard Cartesian unit vectors.  
The reference configuration of the block is setup so its sides are parallel to the Cartesian axes. The surface force density on the top is then defined as 
\begin{align*}
{\bs F}_\text{top}({\bs X},t) = {\bs F}_\text{surface}({\bs X},t) - \left({\bs F}_\text{surface}({\bs X},t) \cdot {\bs e}_2\right){\bs e}_2,
\end{align*}
and on the bottom is defined as
\begin{align*}
{\bs F}_\text{bottom}({\bs X},t) = {\bs F}_\text{surface}({\bs X},t) - \left({\bs F}_\text{surface}({\bs X},t) \cdot {\bs e}_1\right){\bs e}_1.
\end{align*}
In particular, the surface force density on the bottom of the block allows the bottom boundary to slide in the horizontal direction.  
The parameter used to determine the strength of the force is:
\begin{align*}
\kappa = 0.1 \times h \,\Delta t^{-2}.  
\end{align*}  
For this example, we consider both discontinuous and continuous loading pressures to compress the block. 
A discontinuous loading pressure was part of the original benchmark presented in Ref. \cite{Reese99}; our sharp interface method improves results in this case, but it also encounters some difficulties. 
In this light, we also consider a smoothed loading pressure that is defined to be the discontinuous pressure multiplied by a mollifier.   

\subsubsection{Discontinuous loading pressure}

Define ${\bs X} = (X_1, X_2)$.  
The discontinuous loading pressure takes the form
\begin{align*}
\tilde{P}_\text{load}(X_1,t) = 
\begin{cases}
P_\text{ramp}(t) & \text{if } \tilde{a} < X_1 < \tilde{b}, \\
0 & \text{otherwise},
\end{cases}
\end{align*}
with 
\begin{align*}
P_\text{ramp}(t) = 
\begin{cases}
\frac{t}{t_\text{load}} \times P_\text{max} & \text{if } t < t_\text{load},\\
P_\text{max} & \text{otherwise}.
\end{cases}
\end{align*}
In the reference configuration, the block geometry spans the horizontal direction from $X_1 = 0$ mm to $X_1 = 20$ mm and the vertical direction from $X_2 = 0$ mm to $X_2 = 10$ mm. 
We set $\tilde{a} = 5$ mm, $\tilde{b} = 15$ mm, and the loading pressure is applied to the top of the block at $X_2 = 10$ mm.  
The maximum loading pressure is taken to be $P_\text{max} = 200$ N/mm$^2$, and we set $\mu_e = 80.194$ N/mm$^2$.  
The fluid density is $\rho = $1 kg/mm$^3$ and the viscosity is $\mu = $0.16 N$\cdot$s/mm$^2$.

The loading pressure is ramped up over the first 10 s of the simulation, so $t_\text{load} = $ 10 s.  
The model is run to a final time of 50 s to reach an approximate steady state.  
As in the inflating ring test, a combination of normal traction and tangential no--slip velocity boundary conditions are applied to the computational domain $\Omega$ to hold the pressure at zero.  
The computational domain $\Omega$ is $[0,L]^2$ with $L = 30$ mm, and the block is positioned in its center. 
\revthree{The spatial discretization parameter is $h = L/N$ with $N = 2^m$ for $m = 4,5,6,$ and $7$.}  
The time step size is $\Delta t = 0.005 \times h$ and we set $M_\text{fac} = 1$. 

\begin{figure}[h!]
\begin{center}
\includegraphics[scale=0.1,trim=0 0 0 0]{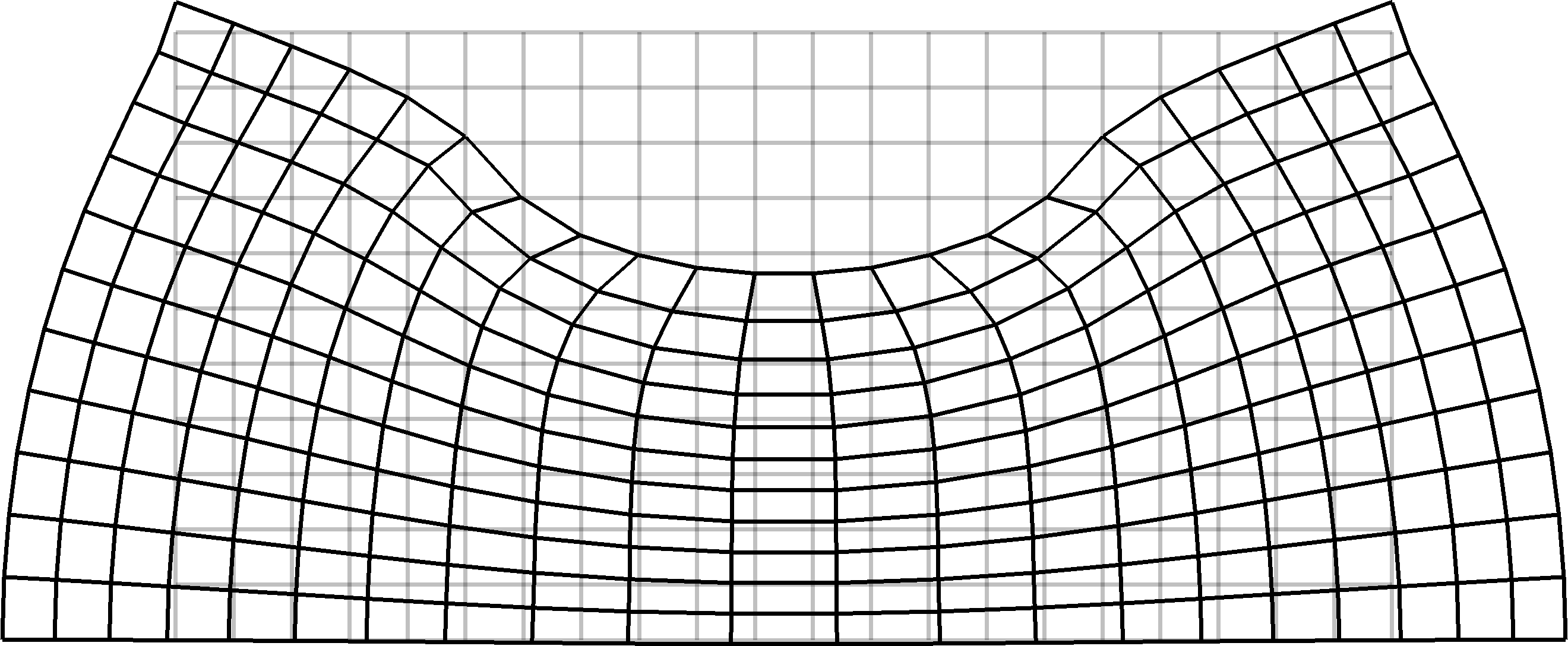}
\end{center}
\caption{The displaced mesh after compression with the discontinuous loading pressure.  
These results are with the sharp interface method, the steady state formulation for $\varphi$, and $\nu = 0$.  
The initial configuration of the block is in the background.}
\label{fig:block6}
\end{figure}

\begin{figure}[h!]
\begin{center}
\includegraphics[scale=0.4,trim=0 -40 0 0]{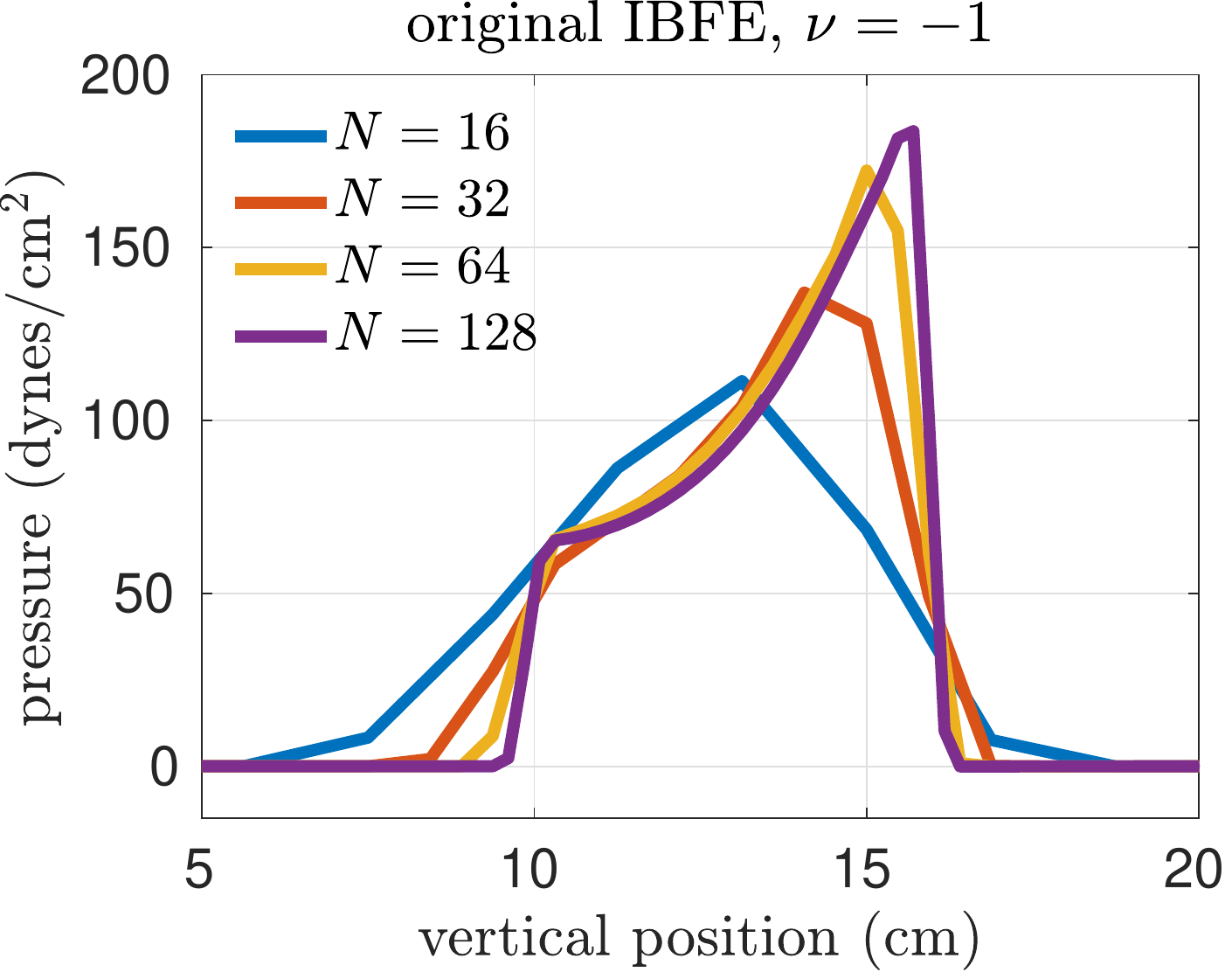}
\includegraphics[scale=0.4,trim=0 -62 0 0]{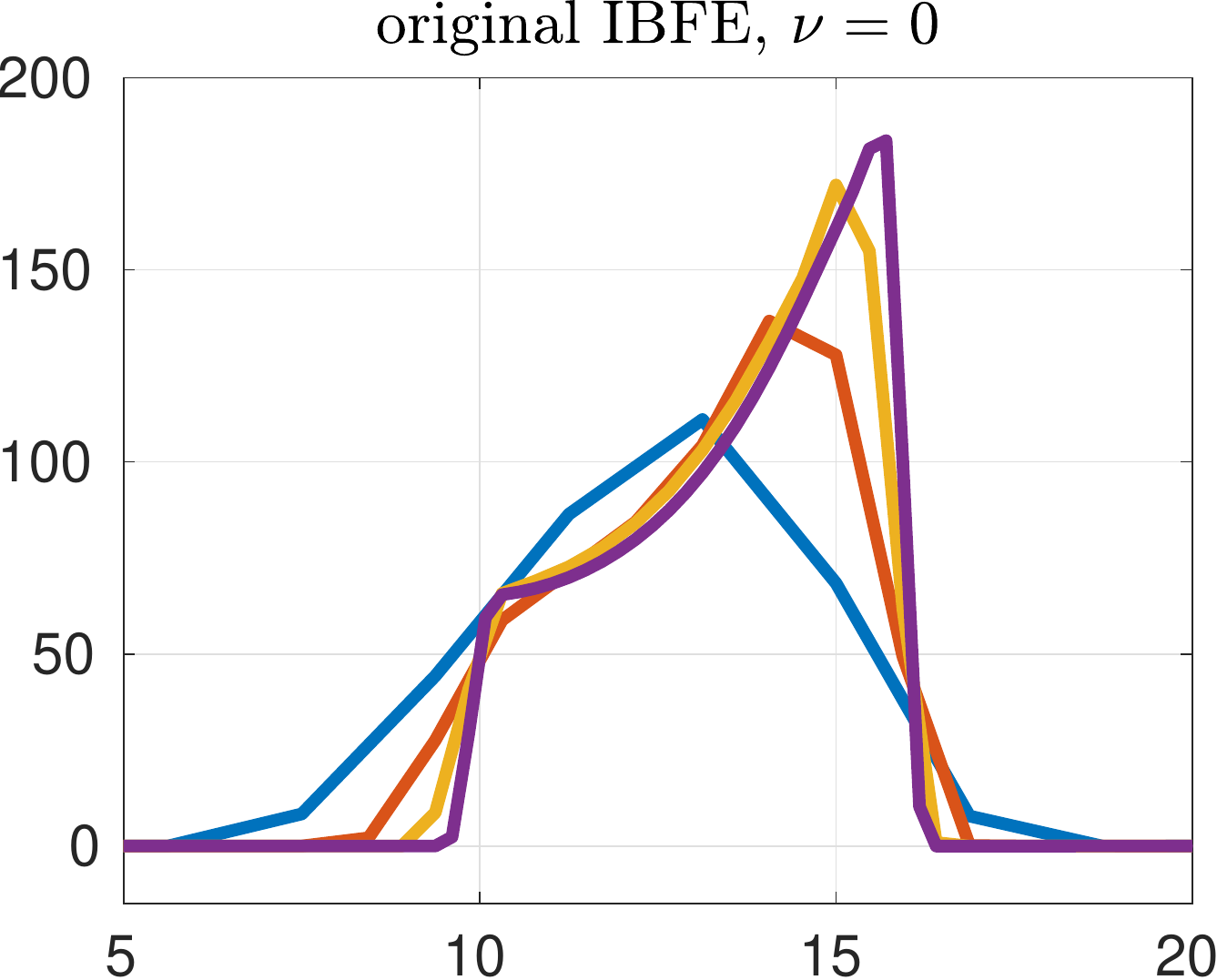}
\includegraphics[scale=0.4,trim=0 -62 0 0]{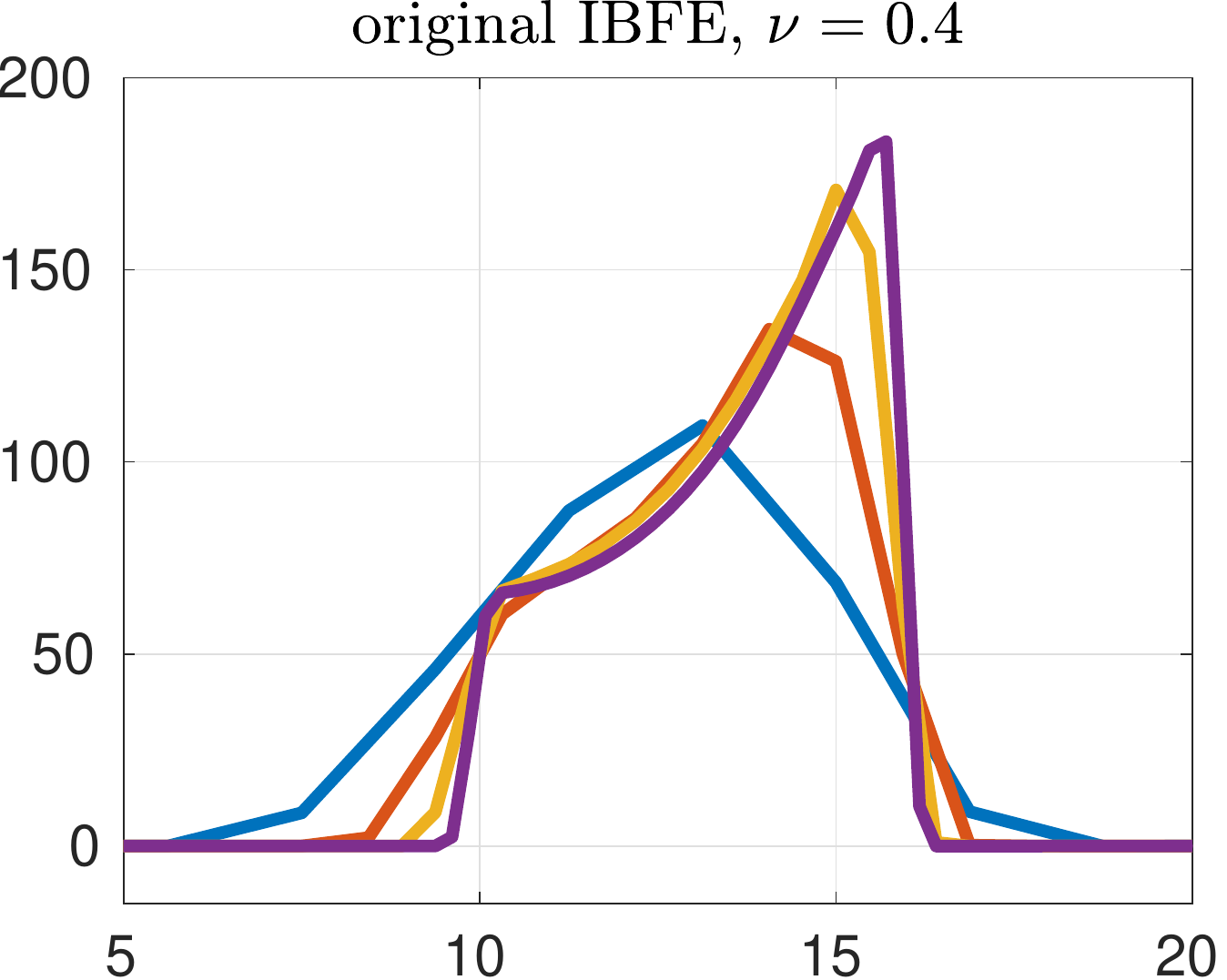} \\ 
\includegraphics[scale=0.4,trim=0 0 0 0]{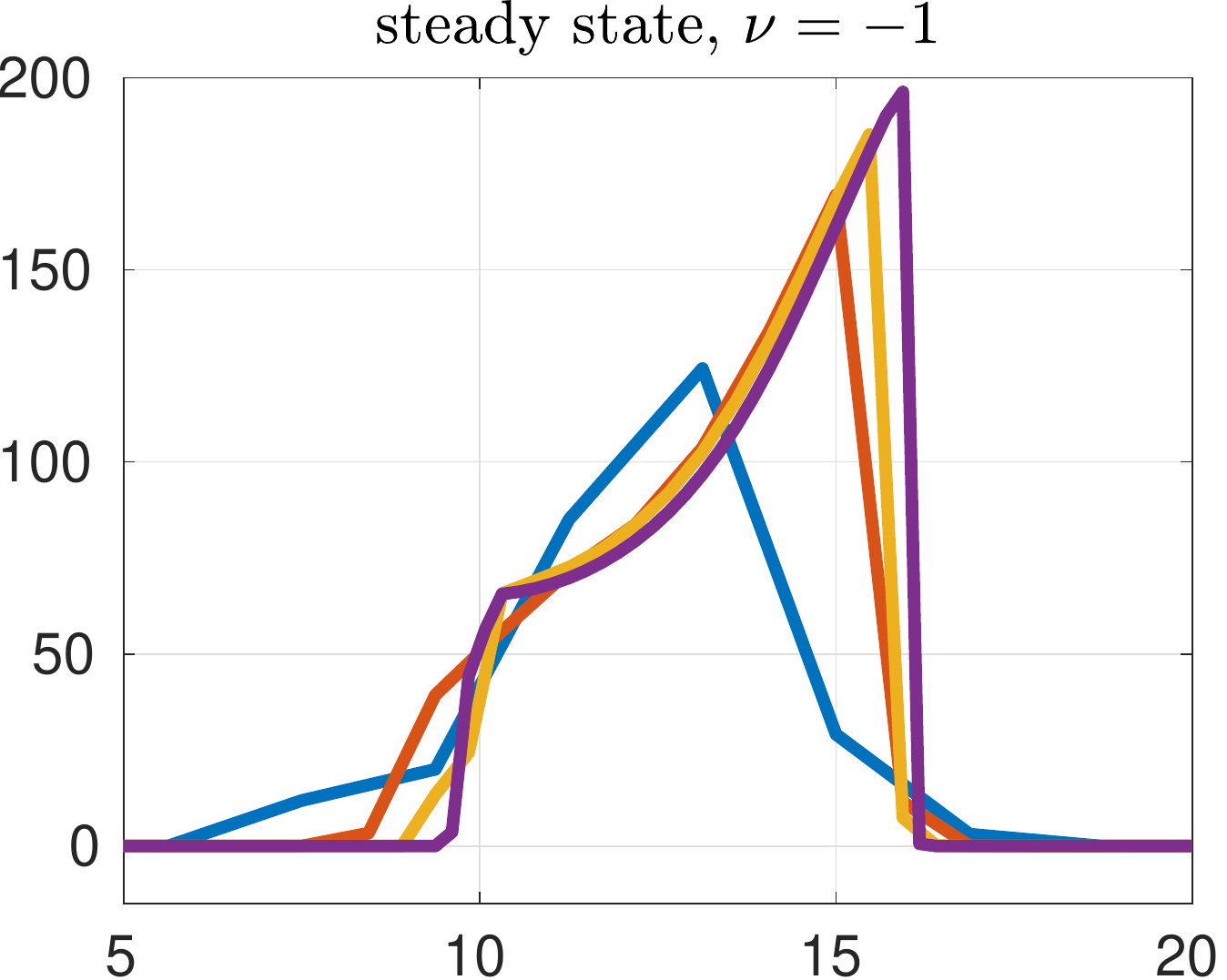}
\includegraphics[scale=0.4,trim=0 0 0 0]{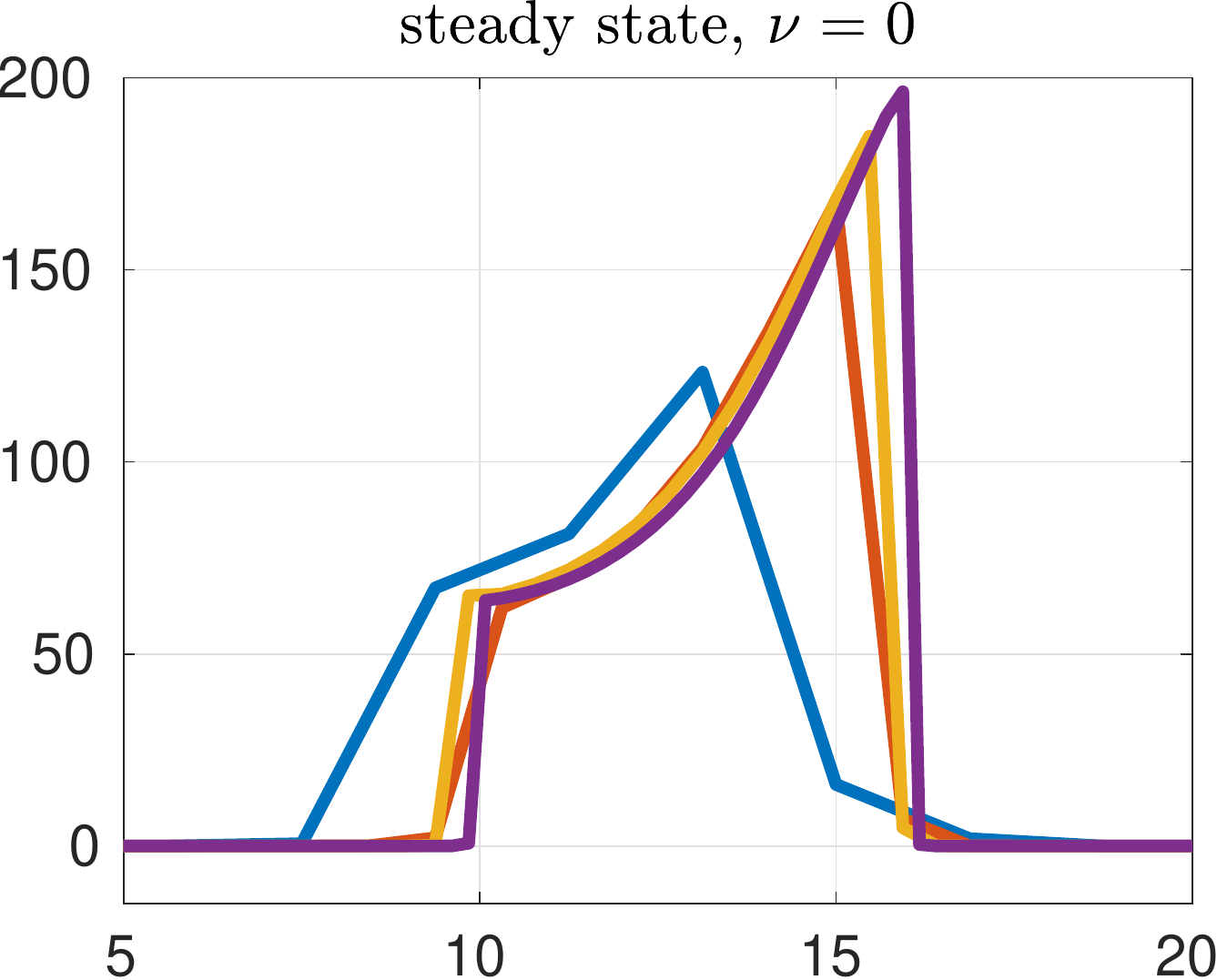}
\includegraphics[scale=0.4,trim=0 0 0 0]{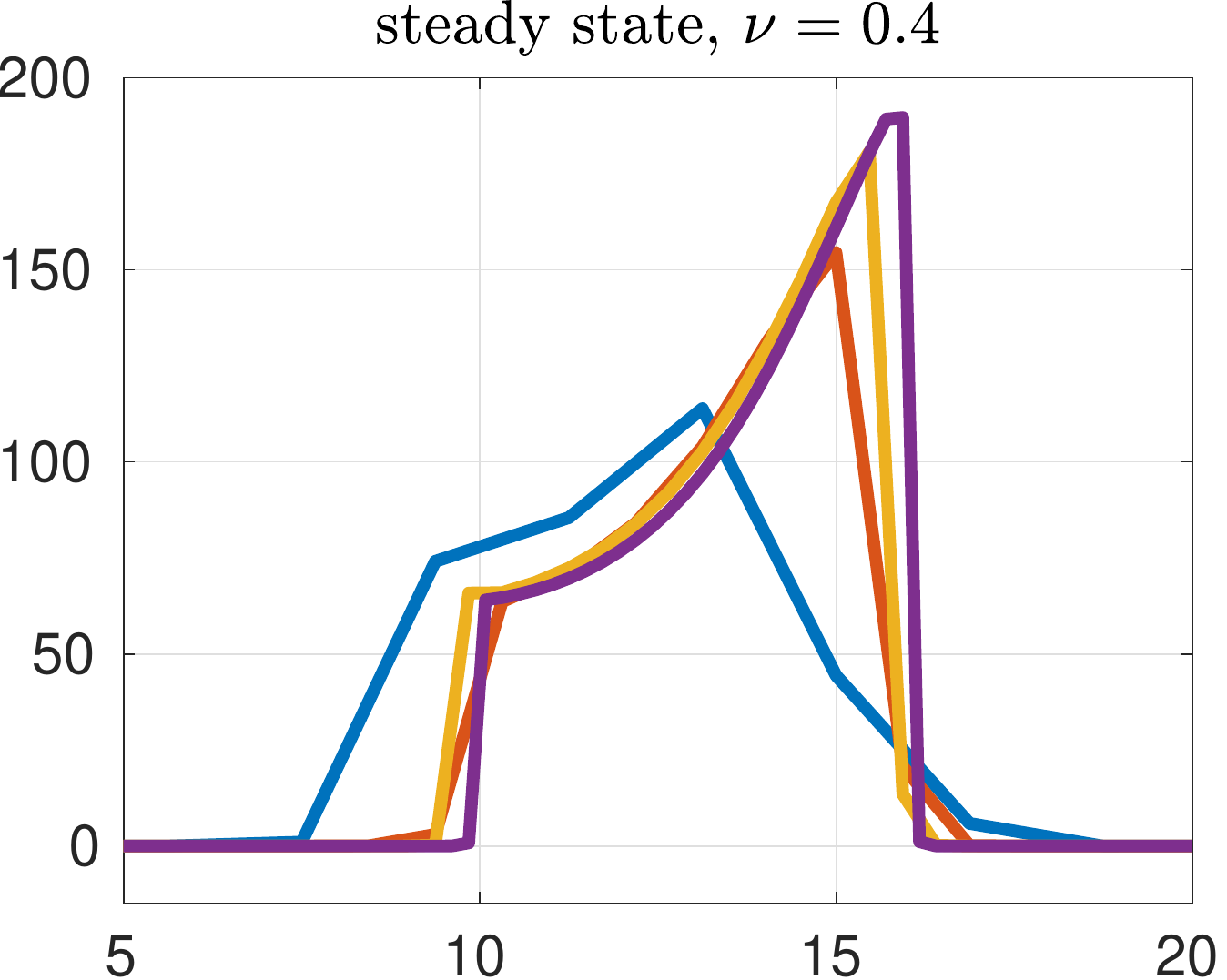} 
\caption{Slices of the pressure field vertically down the center of the block, at $X_1 = 15$, for the case with the discontinuous loading pressure.  
Results from the original IBFE method are on the top row, and results from the sharp interface method with the steady state formulation for $\varphi$ are on the bottom row. The sharp interface method yields crisper pressure discontinuities and faster convergence.}
\label{fig:block5}
\end{center}
\end{figure}

Figure \ref{fig:block6} shows the displaced mesh for the compressed block with the discontinuous loading pressure.  
For this result, the steady state formulation for $\varphi$ is used and $\nu = 0$. 
Figure \ref{fig:block5} examines slices of the pressure field down the center of the block at $X_2 = 15$ mm, with the original IBFE method on the top row and the sharp interface method on the bottom row.  
As $N$ increases, the pressure slices converge to a sharp profile in both methods. 
The sharp interface method converges faster and more cleanly resolves pressure discontinuities on the top and bottom faces of the block.  
For the case $\nu = 0.4$, the discontinuities in the pressure field are less clear with the sharp interface method.   
This observation indicates some subtle interplay between the volumetric energy and the splitting of the pressure field, which requires additional investigation.  
We consider this issue in the next section with the continuous loading pressure. 

\begin{figure}[h!]
\begin{center}
\includegraphics[scale=0.11,trim=0 0 -250 0]{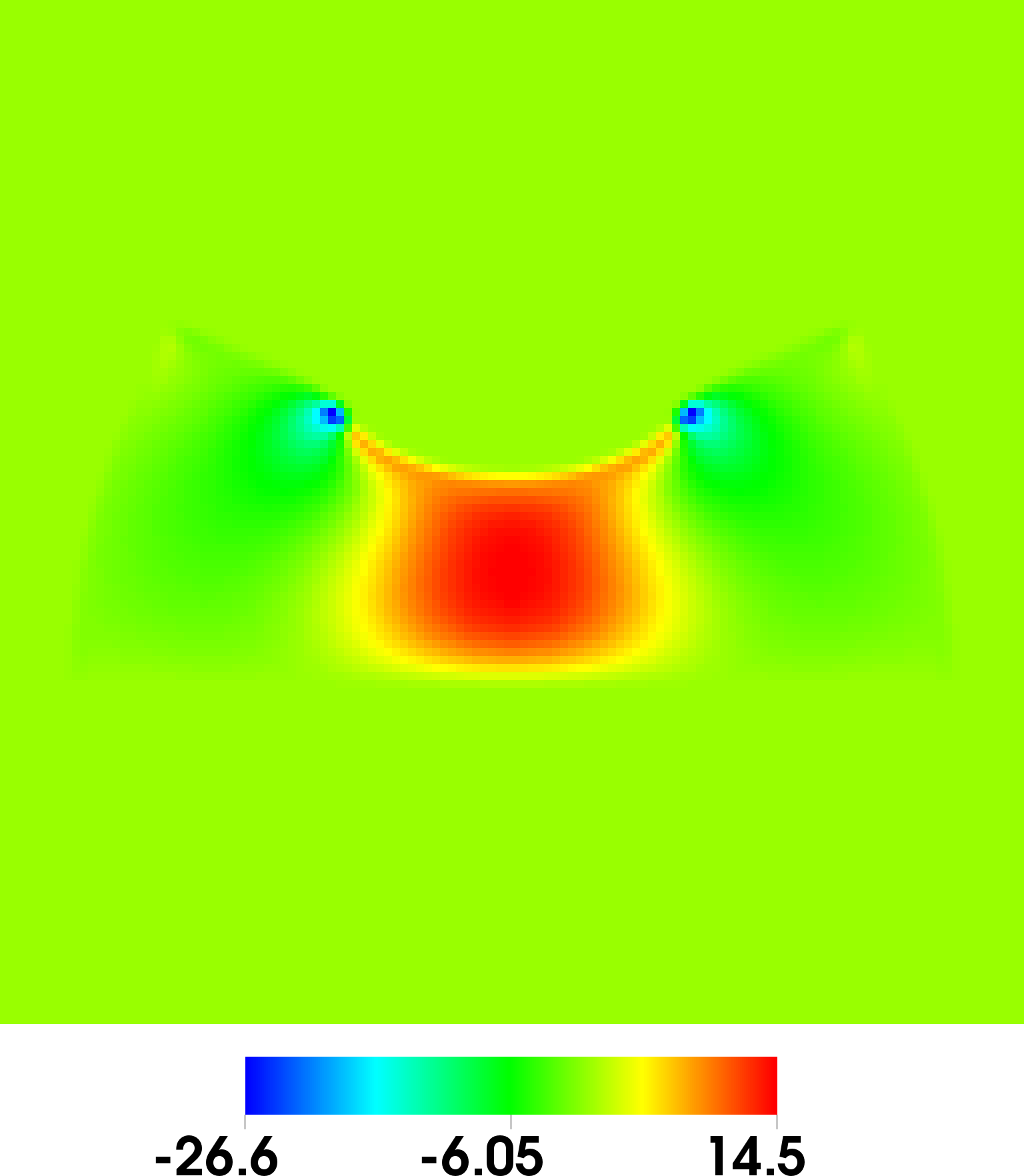}
\includegraphics[scale=0.11,trim=-250 0 0 0]{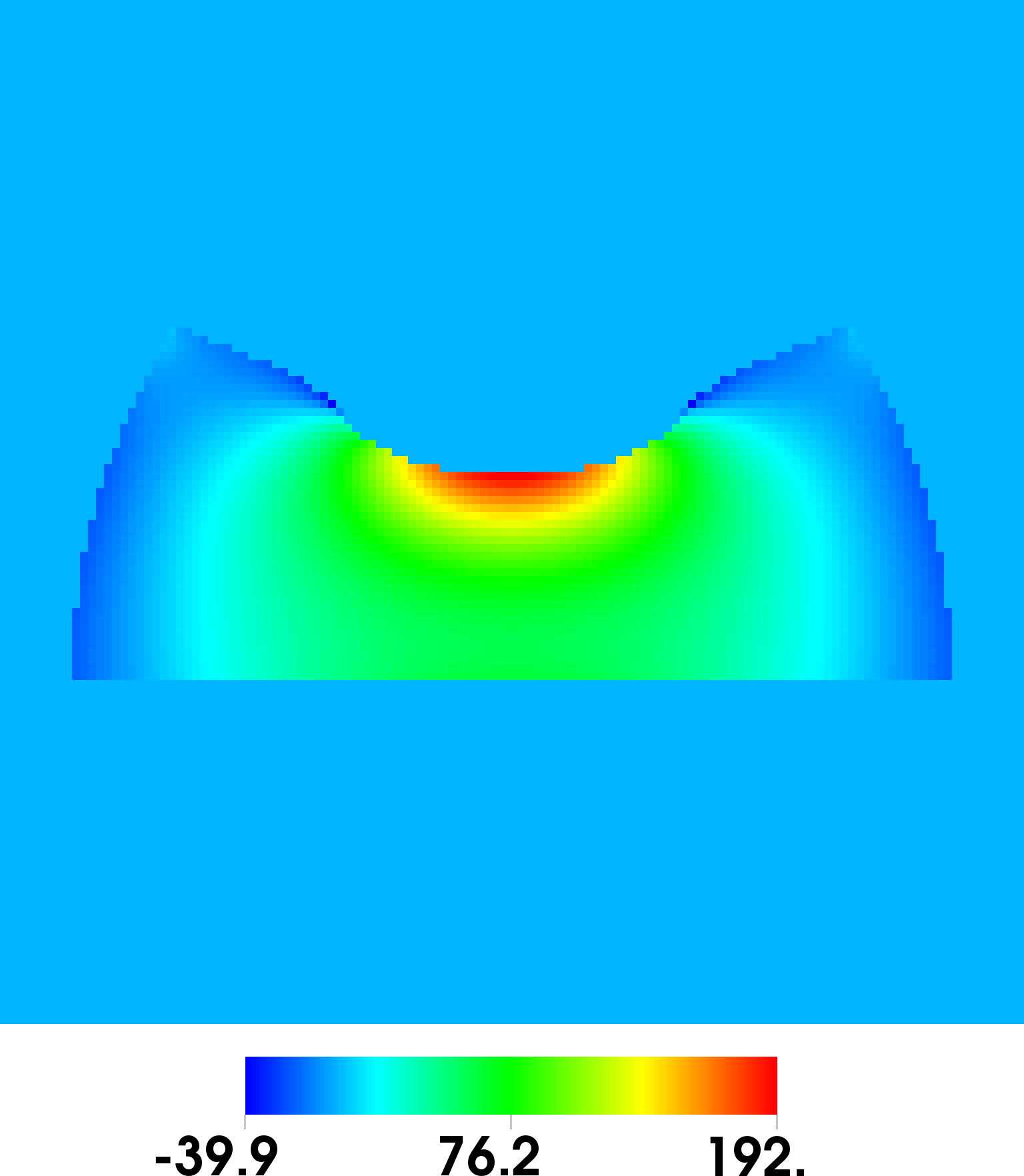}
\caption{The pressure--like fields for the compressed block with the discontinuous loading pressure. 
On the left is the $\pi$ field and on the right is the $\varphi$ field, computed with the steady state formulation, with $\nu = 0$.  
The Cartesian grid composed of $128^2$ cells.  The $\pi$ field appears to have discontinuities at the locations of singularities of the force used to compress the block.}
\label{fig:block7}
\end{center}
\end{figure}

\begin{figure}[h!]
\begin{center}
\includegraphics[scale=0.11,trim=0 0 -250 0]{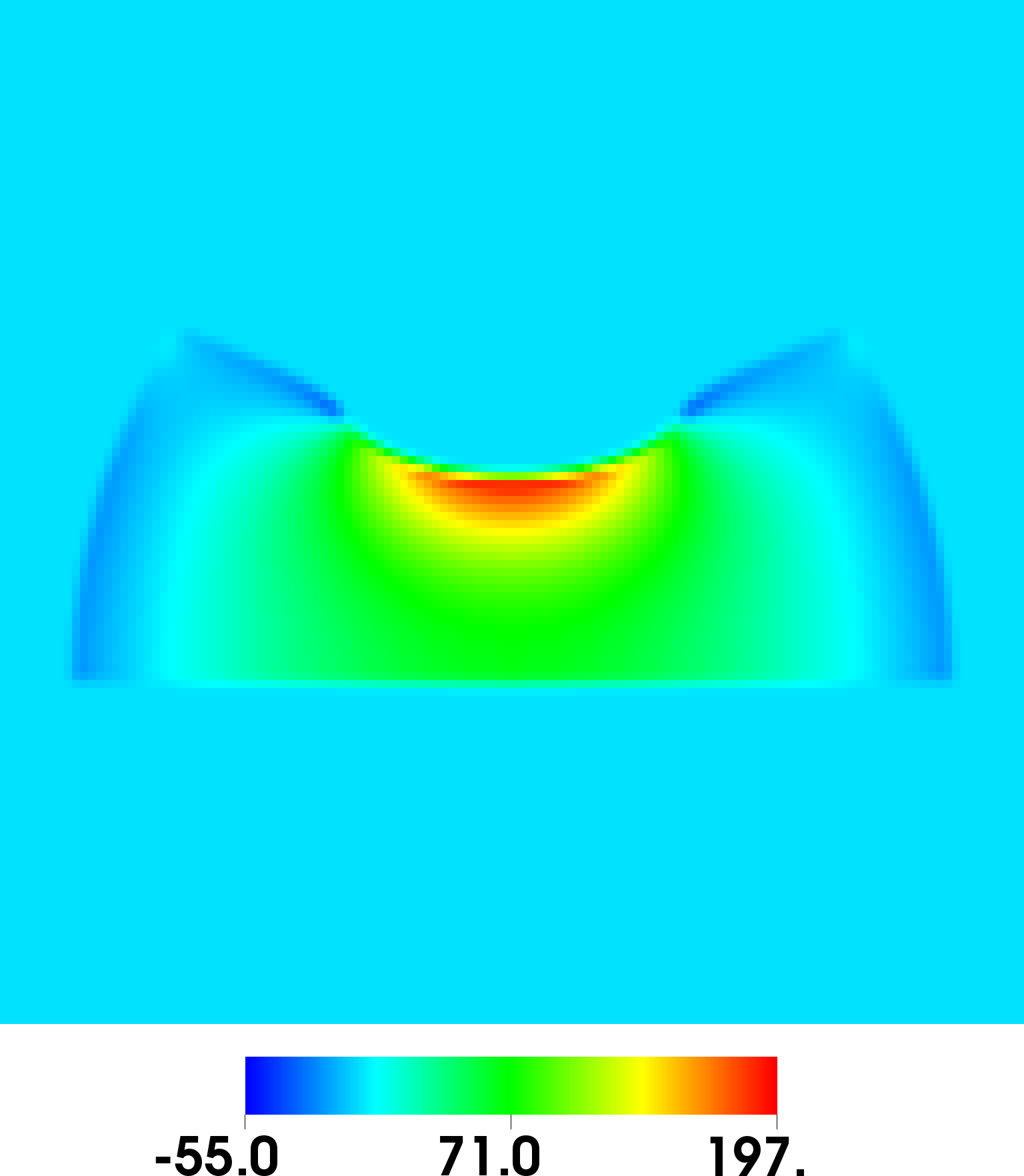}
\includegraphics[scale=0.11,trim=-250 0 0 0]{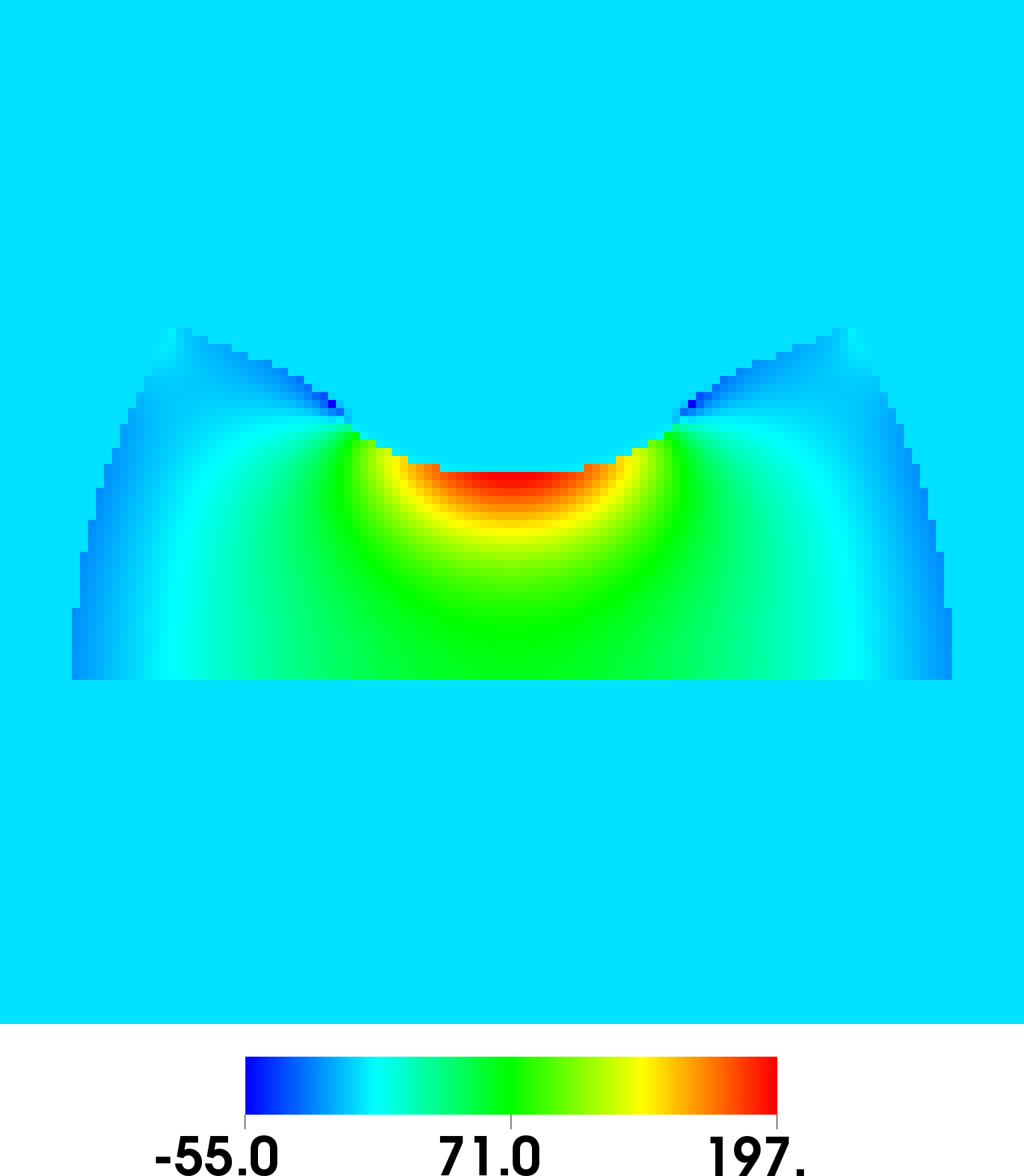}
\caption{The pressure fields for the compressed block with the discontinuous loading pressure.  
On the left is the pressure field computed with the original IBFE method, and on the right is the pressure field computed with the steady state solver for $\varphi$.  
The parameter $\nu = 0$, and the Cartesian grid contains $128^2$ cells.}
\label{fig:block8}
\end{center}
\end{figure}

In Figure \ref{fig:block7}, we plot $\pi$ and $\varphi$ for $\nu = 0$.  
The $\pi$ field appears to contain discontinuities at the approximate location of the jumps in the pressure loading force.  
The sharp interface method presented here assumes that $\pi$ is continuous, so this example presents some challenges for our approach.  
Even so, the physical pressure field from the sharp interface method, in Figure \ref{fig:block8} on the right, more sharply captures pressure discontinuities compared to the original IBFE method. 

\subsubsection{Continuous loading pressure}

The smoothed version of the loading pressure we use takes the form:

\begin{align*}
P_\text{load}(X_1,t) = 
\begin{cases}
P_\text{ramp}(t) \times \exp\left(\frac{(b-a)^2}{(2X_1 - a - b)^2 - (b-a)^2} + 1\right) & \text{if } a < X_1 < b \\
0 & \text{otherwise},
\end{cases}
\end{align*}
The parameters $a = 4$ mm, $b = 16$ mm, and all other numerical and physical parameters are the same as in the previous section.

\begin{figure}[h!]
\begin{center}
\includegraphics[scale=0.1,trim=0 0 0 0]{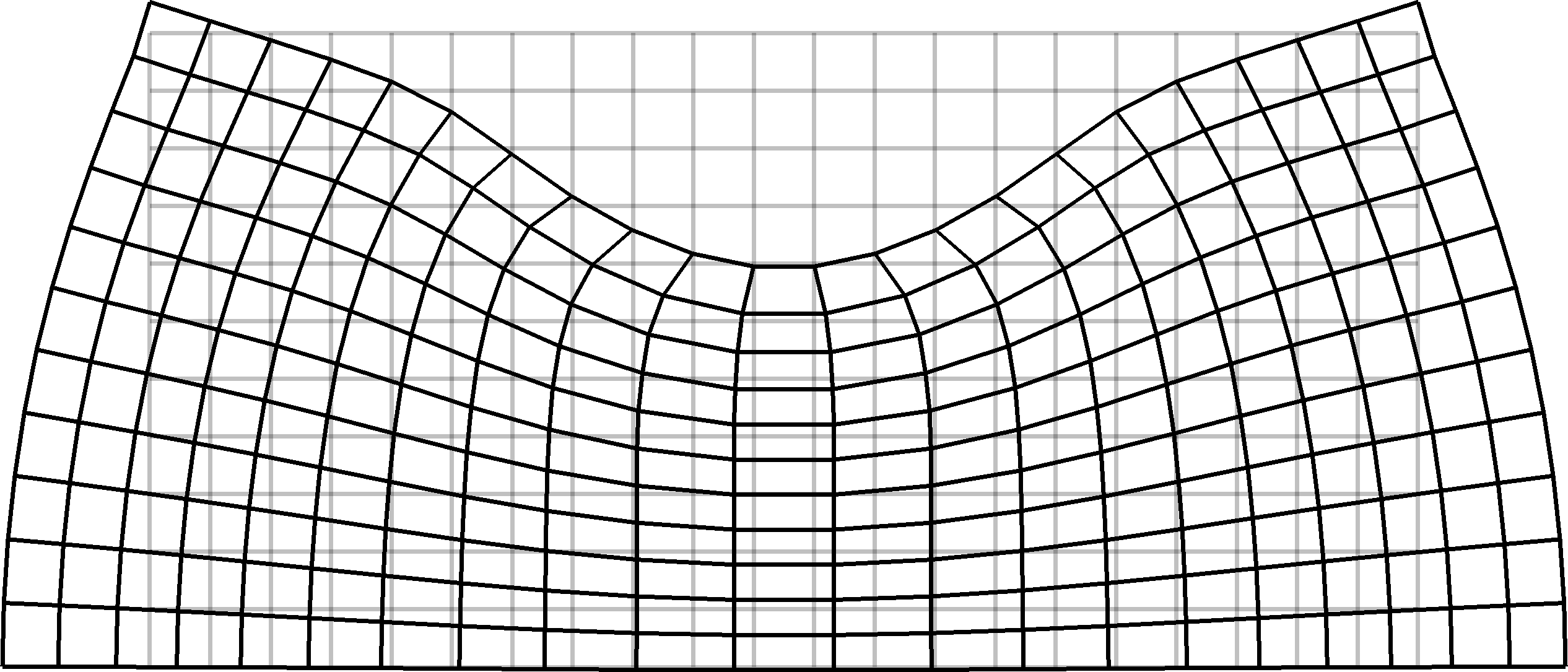}
\end{center}
\caption{The displaced mesh after compression with the continuous loading pressure.  
These results are with the sharp interface method, the steady state formulation for $\varphi$, and $\nu = 0$.  
The initial configuration of the block is in the background.}
\label{fig:block1}
\end{figure}

\begin{figure}[h!]
\begin{center}
\includegraphics[scale=0.09,trim=0 0 -100 0]{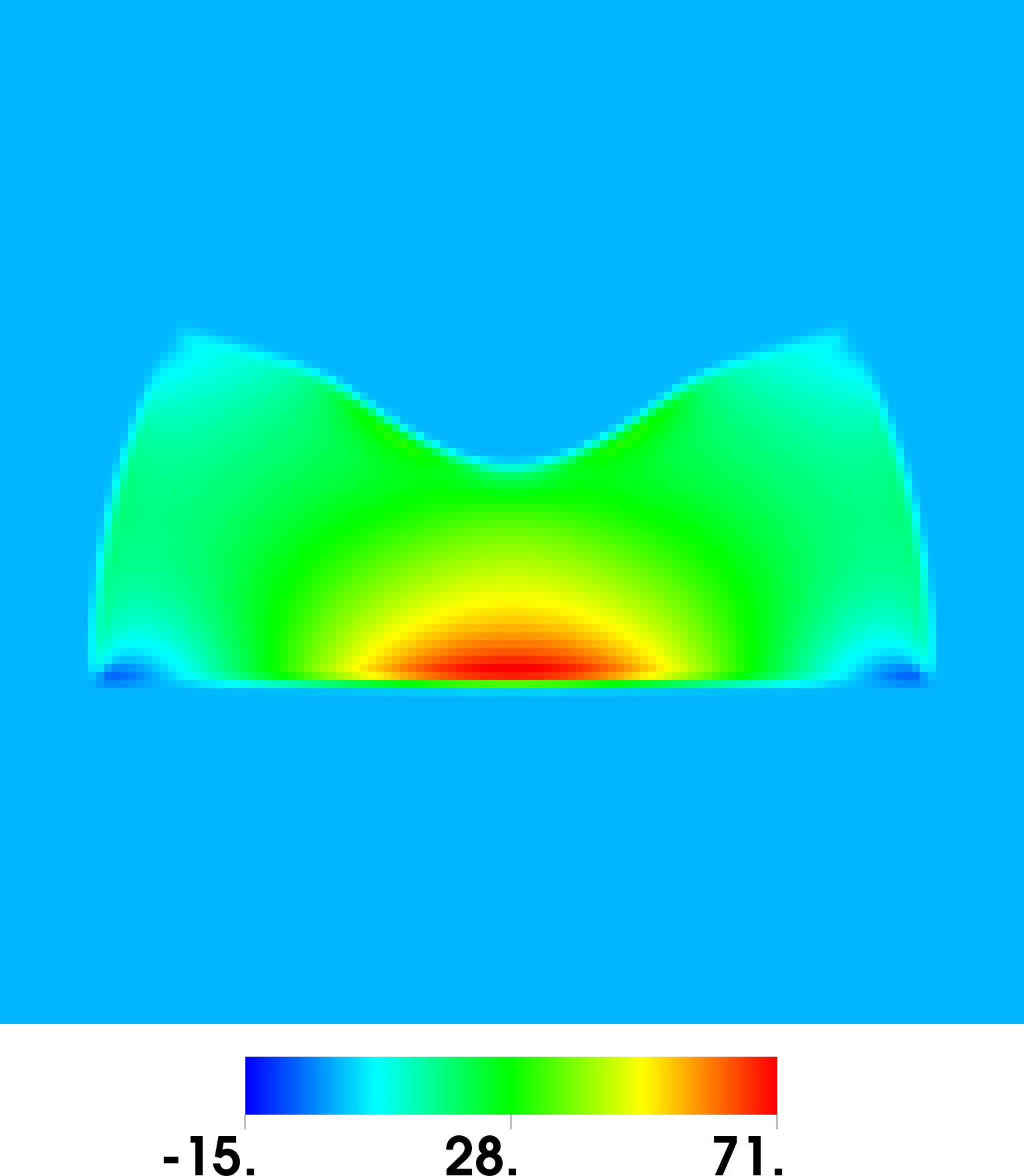}
\includegraphics[scale=0.09,trim=-100 0 -100 0]{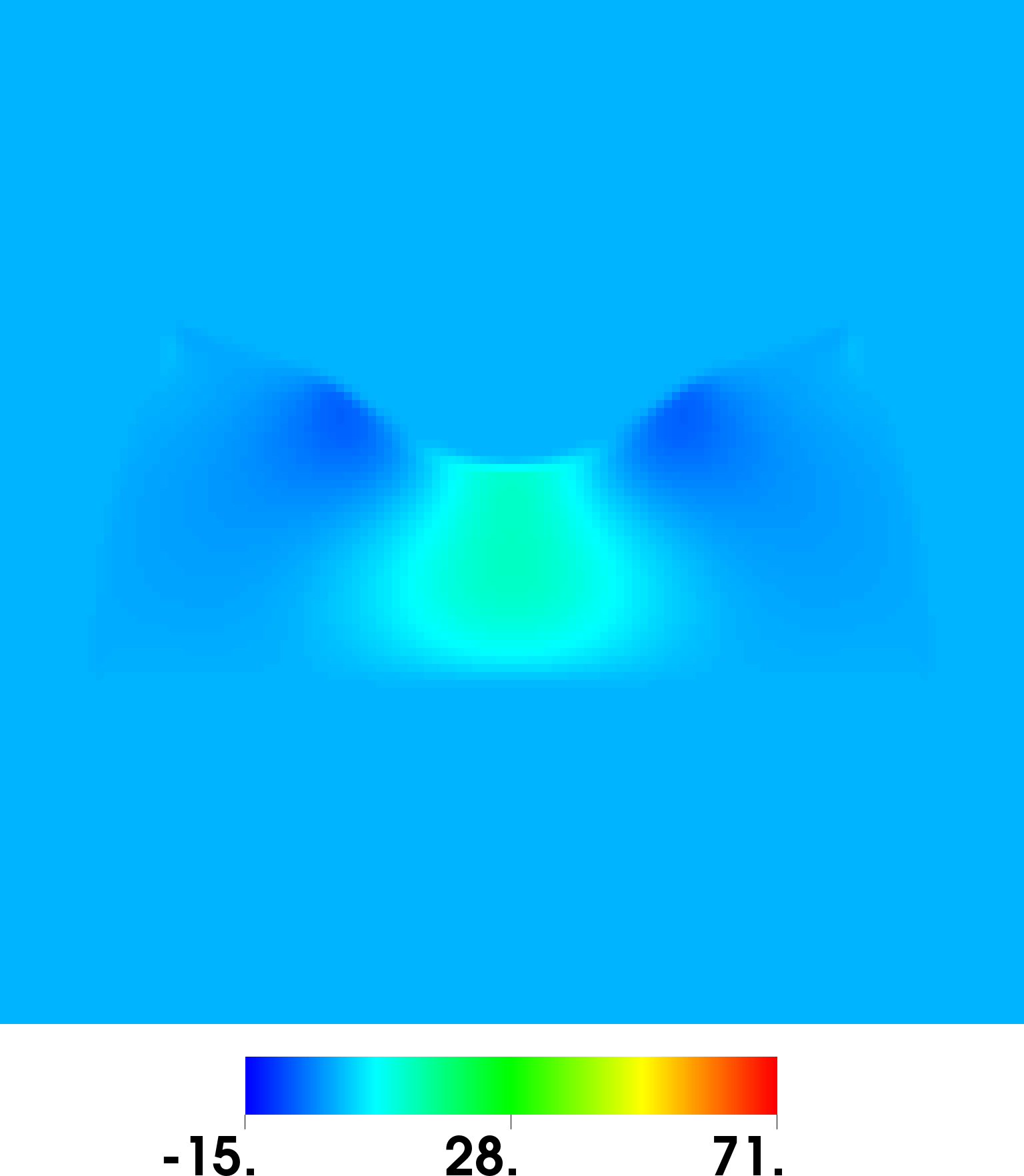}
\includegraphics[scale=0.09,trim=-100 0 0 0]{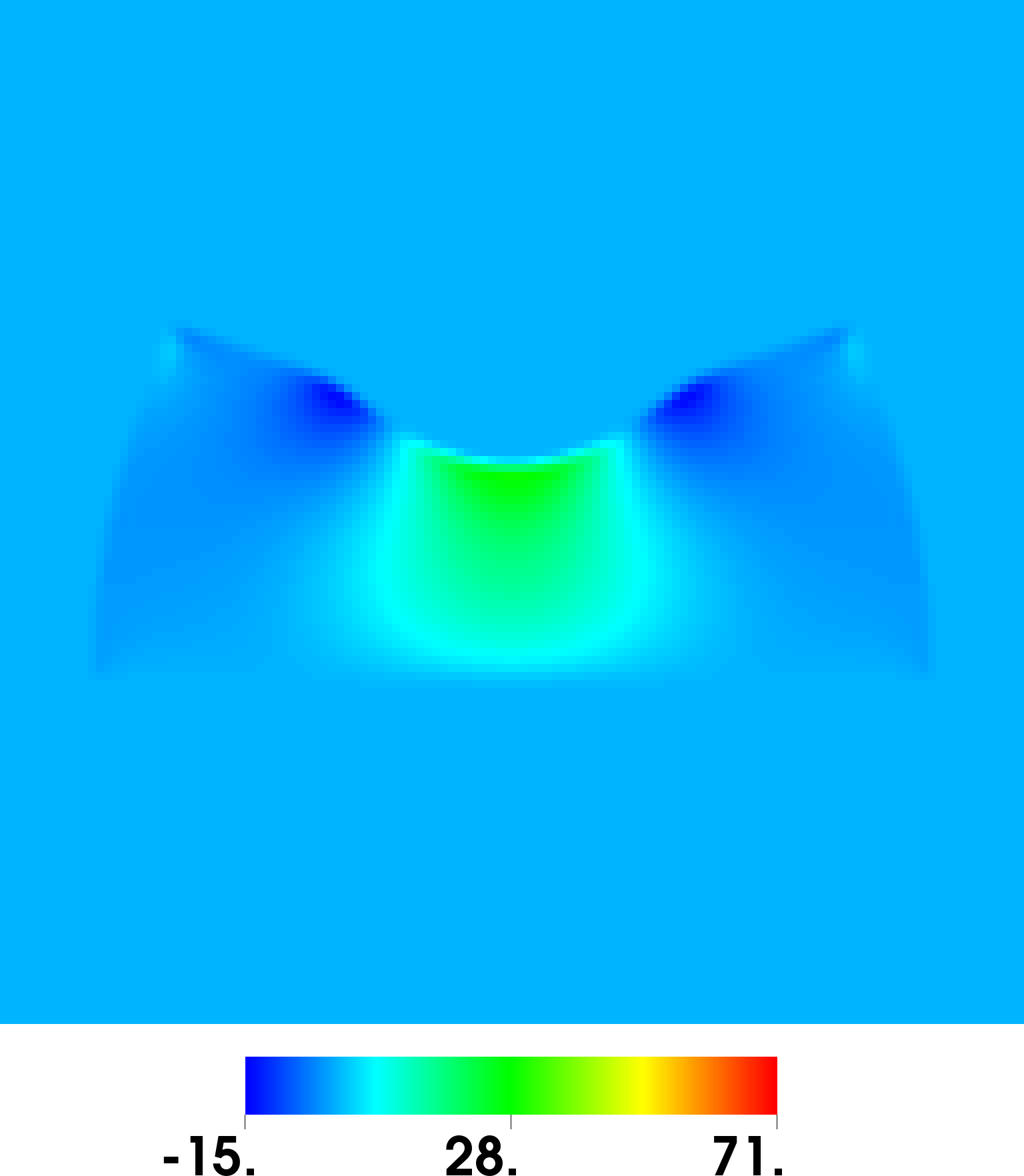}
\put(-410,165){$\nu = -1$}
\put(-245,165){$\nu = 0$}
\put(-87,165){$\nu = 0.4$}
\caption{Visualization of the $\pi$ field for the compressed block test, with the continuous loading pressure, for $\nu = -1, 0, 0.4$.  
The grid contains $128^2$ cells, and the steady state formulation is used to solve for $\varphi$.}
\label{fig:block_pi_fields}
\end{center}
\end{figure}

\begin{figure}[h!]
\begin{center}
\includegraphics[scale=0.09,trim=0 0 -100 0]{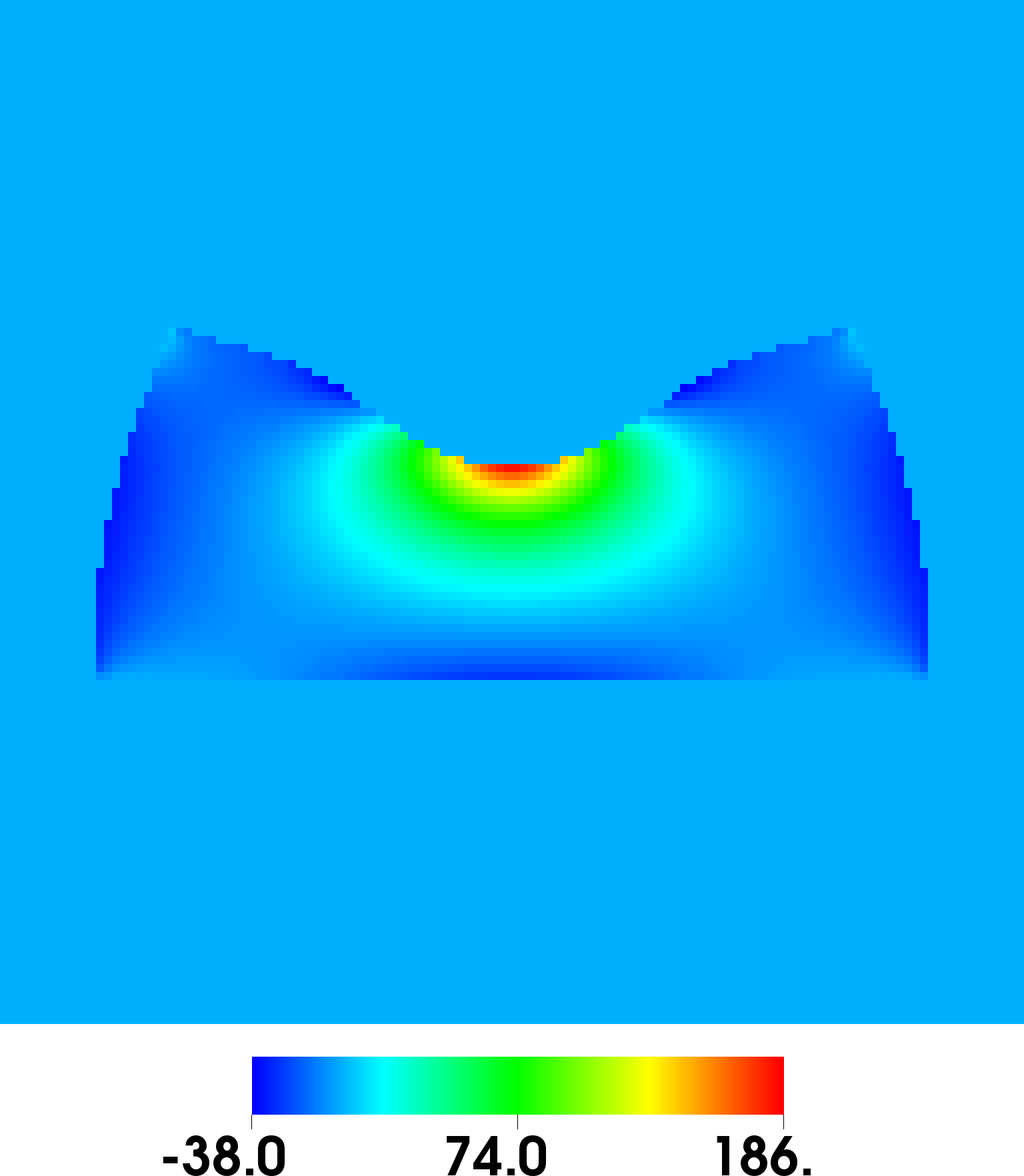}
\includegraphics[scale=0.09,trim=-100 0 -100 0]{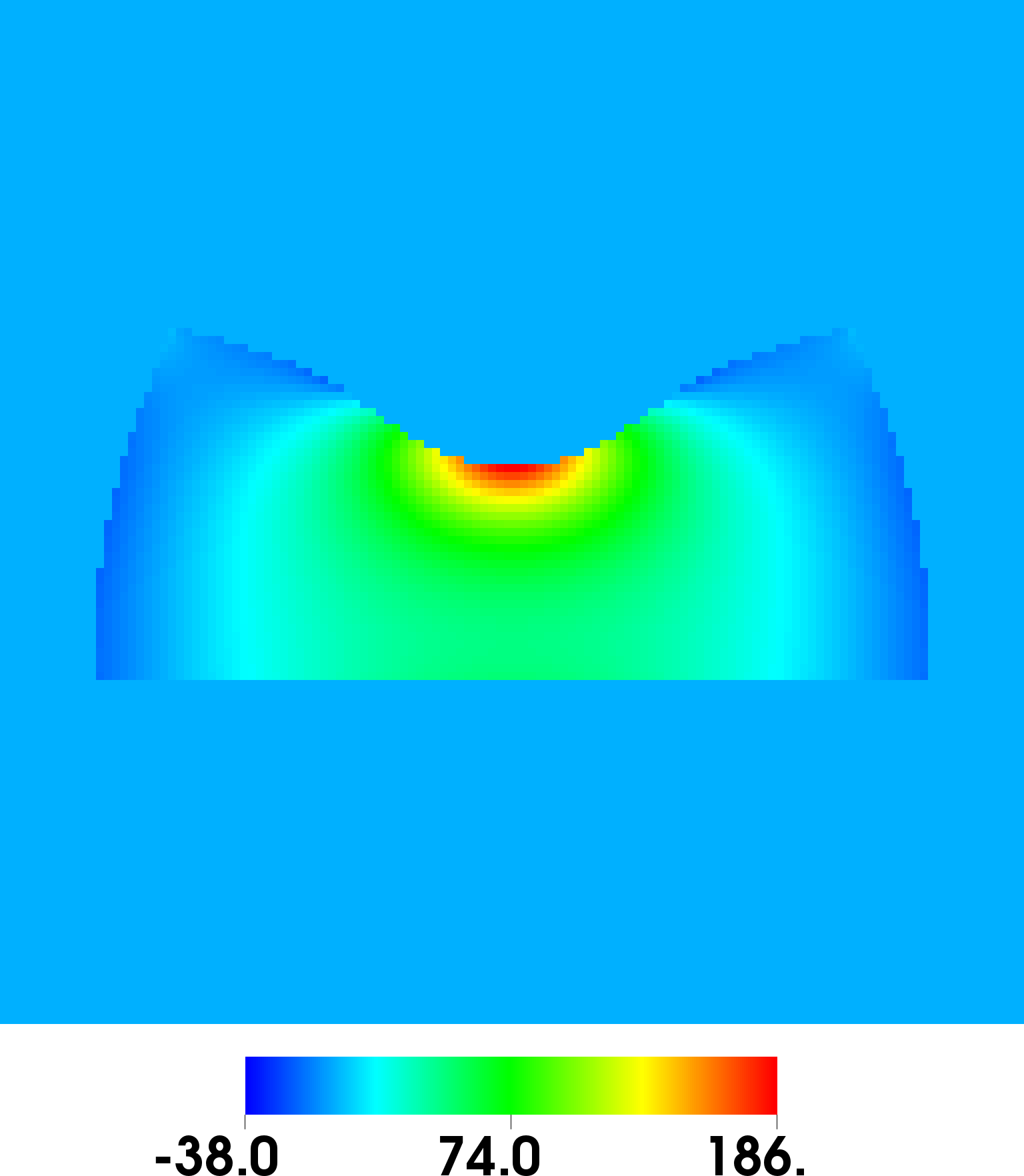}
\includegraphics[scale=0.09,trim=-100 0 0 0]{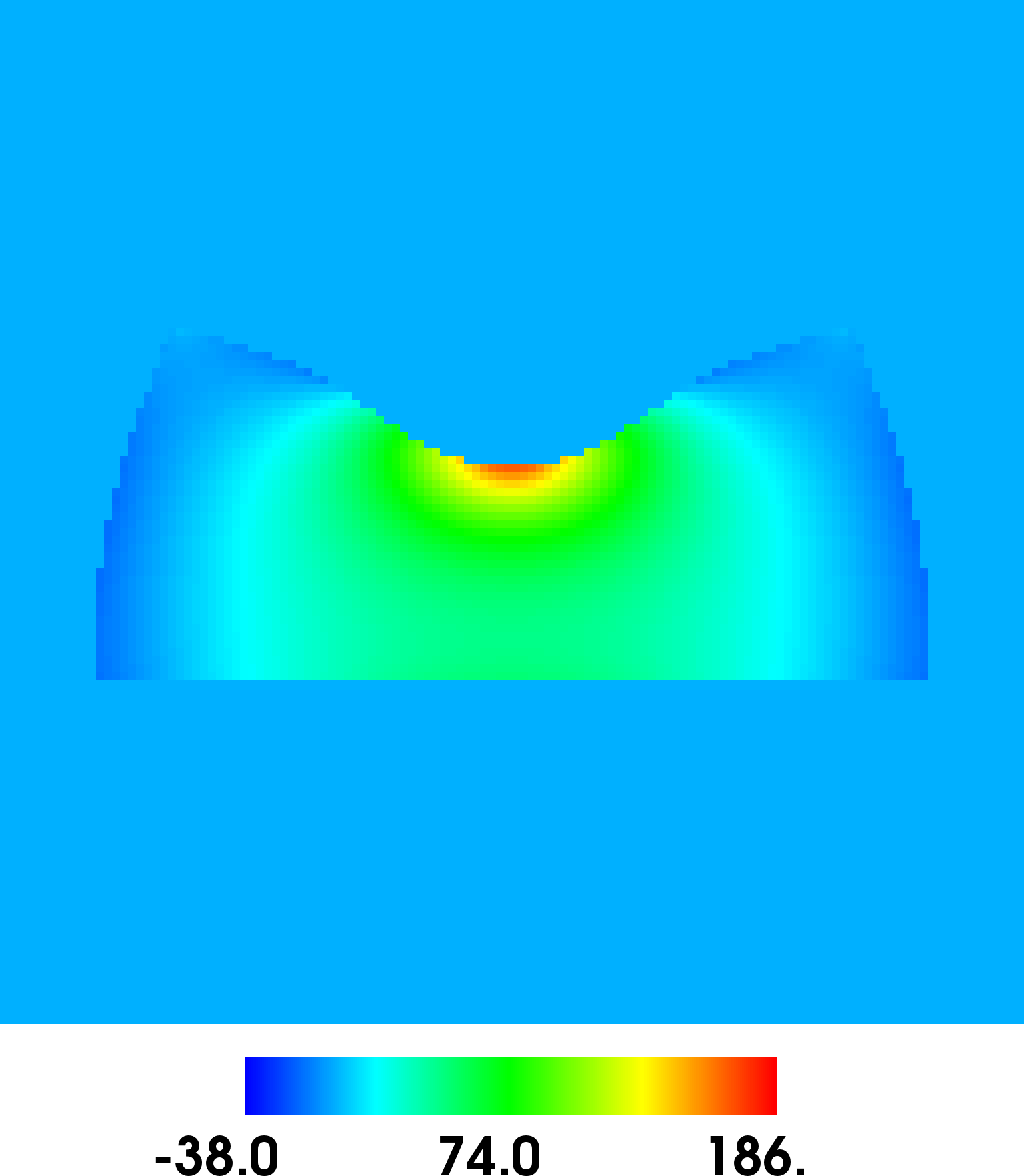}
\put(-410,165){$\nu = -1$}
\put(-245,165){$\nu = 0$}
\put(-87,165){$\nu = 0.4$}
\caption{Visualization of the $\varphi$ field for the compressed block test, with the continuous loading pressure, for $\nu = -1, 0, 0.4$.  
The grid contains $128^2$ cells, and the steady state formulation is used to solve for $\varphi$.}\label{fig:block_phi_fields}
\end{center}
\end{figure}

Figure \ref{fig:block1} depicts the displaced mesh after the continuous loading pressure is applied, compared to the initial mesh in the background. 
First, we more closely examine the impact of the volumetric energy on the sharp interface method by displaying the $\pi$ and $\varphi$ fields separately in Figures \ref{fig:block_pi_fields} and \ref{fig:block_phi_fields}.  
When the volumetric energy is zero, i.e. $\nu = -1$, the $\pi$ field has a large gradient at the bottom of the block.  
When the volumetric energy is included in the strain energy density, corresponding to the cases $\nu = 0$ and $\nu = 0.4$, the gradients in the $\pi$ field are more mild, with larger values toward the top of the block where the pressure load is applied.  
The $\varphi$ field in all cases takes it largest values at the top of the block, and in the cases with nonzero volumetric energy has smaller discontinuities on the block boundary.

\begin{figure}[h!]
\begin{center}
\includegraphics[scale=0.11,trim=0 0 -250 0]{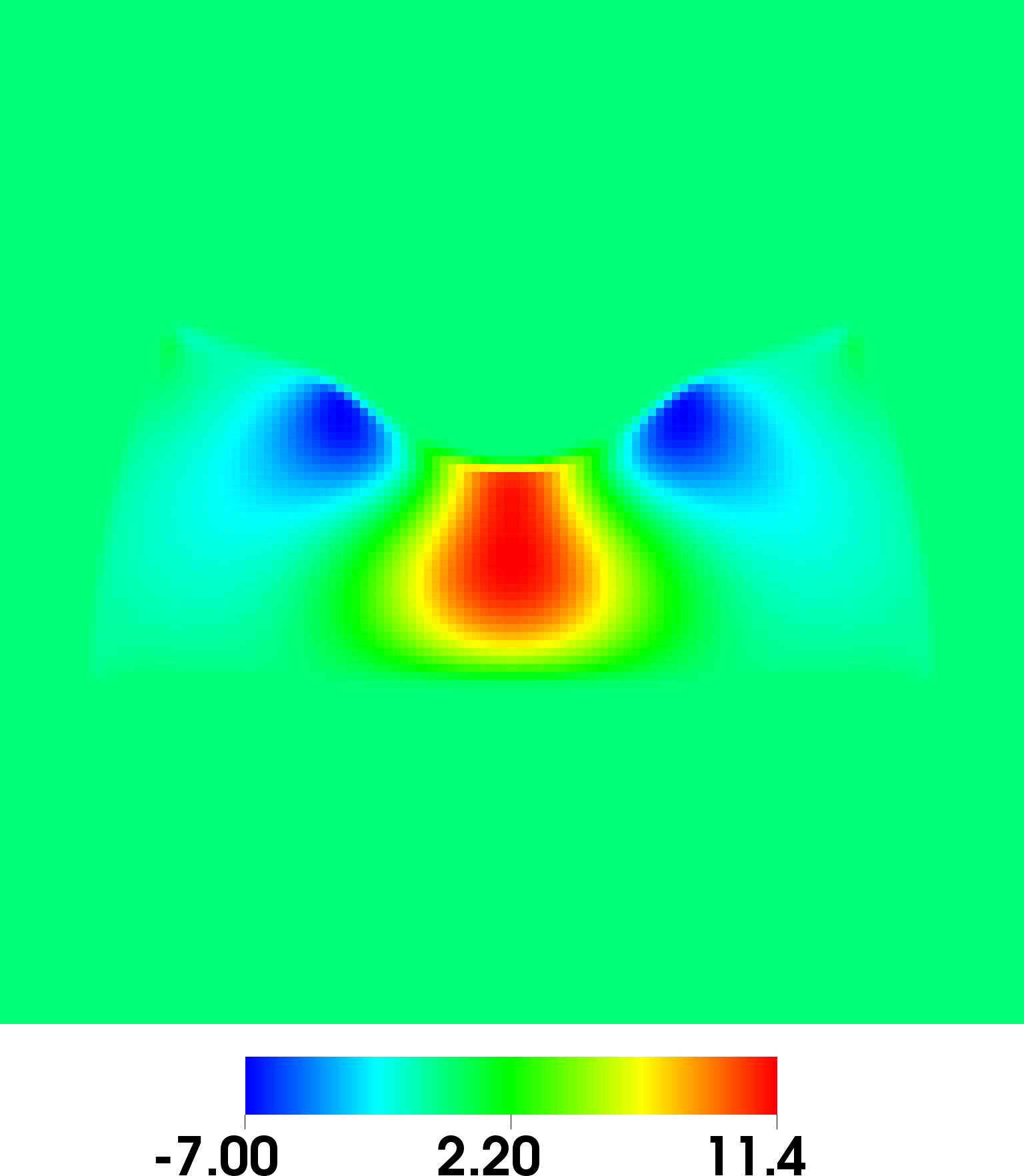}
\includegraphics[scale=0.11,trim=-250 0 0 0]{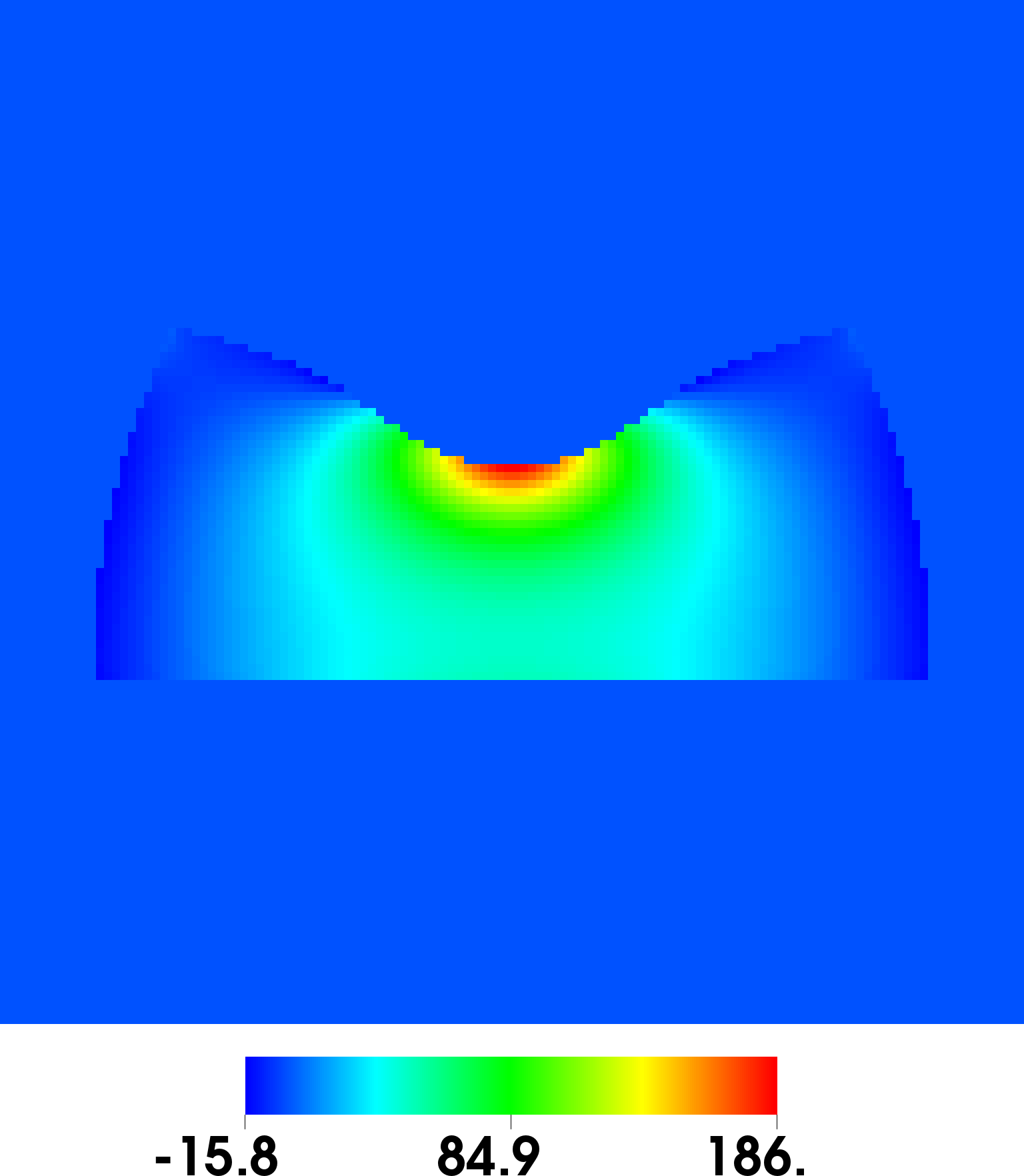}
\caption{The pressure--like fields for the compressed block with the continuous loading pressure.  
On the left is the $\pi$ field and on the right is the $\varphi$ field, computed with the steady state formulation, with $\nu = 0$.  
The Cartesian grid composed of $128^2$ cells.  The $\pi$ field appears to be much smoother than in the case with the discontinuous loading pressure.}
\label{fig:block2}
\end{center}
\end{figure}

\begin{figure}[h!]
\begin{center}
\includegraphics[scale=0.11,trim=0 0 -250 0]{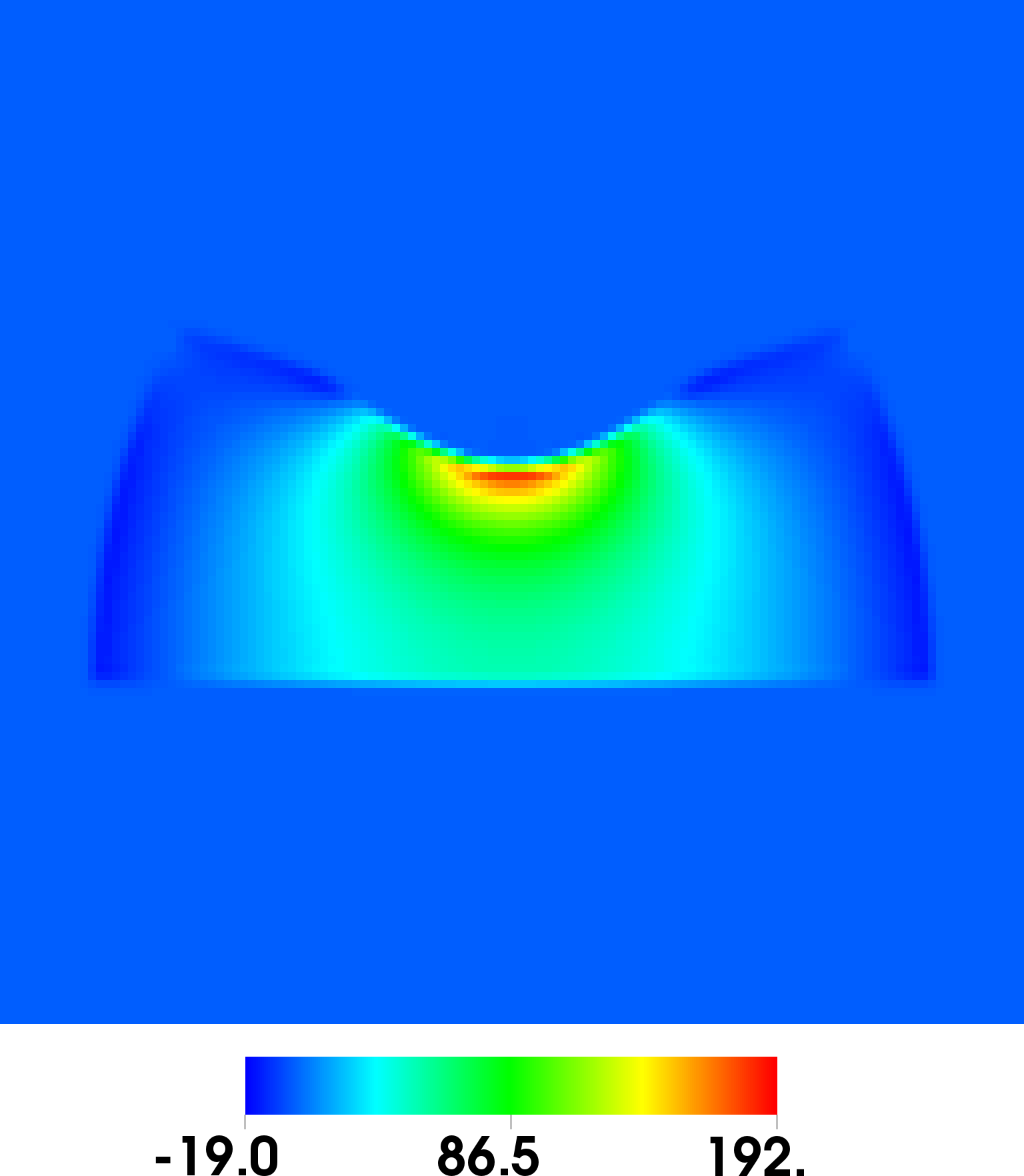}
\includegraphics[scale=0.11,trim=-250 0 0 0]{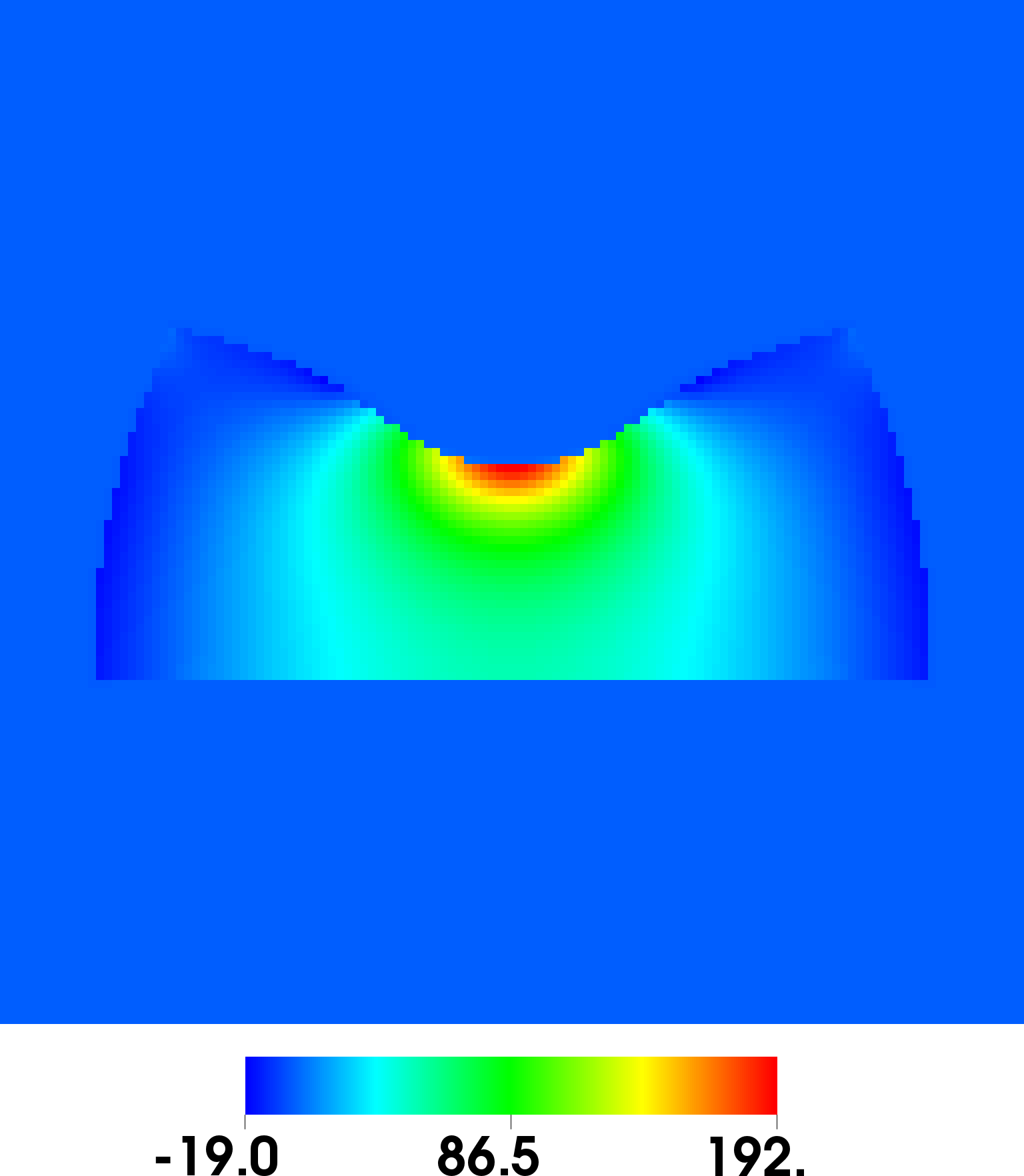}
\caption{The pressure fields for the compressed block with the continuous loading pressure. 
On the left is the pressure field computed with the original IBFE method, and on the right is the pressure field computed with the steady state solver for $\varphi$.  
The parameter $\nu = 0$, and the Cartesian grid contains $128^2$ cells.}
\label{fig:block3}
\end{center}
\end{figure}

In Figure \ref{fig:block2}, we examine the $\pi$ and $\varphi$ fields for the case $\nu = 0$.  
Not surprisingly, the largest pressure discontinuities in this test appear on the top of the block where the loading pressure is applied. Further, the computed $\pi$ field is more smooth than in the case with the discontinuous loading pressure.  
Figure \ref{fig:block3} shows the physical pressure field; on the left are the results from the original IBFE method and on the right is the sum $\pi + \varphi$ from the fields in Figure \ref{fig:block2}.  
This sharp interface method better resolves the pressure discontuities, especially in the center of the top of the block where the loading pressure is applied. 

\begin{figure}[h!]
\begin{center}
\includegraphics[scale=0.4,trim=0 -40 0 0]{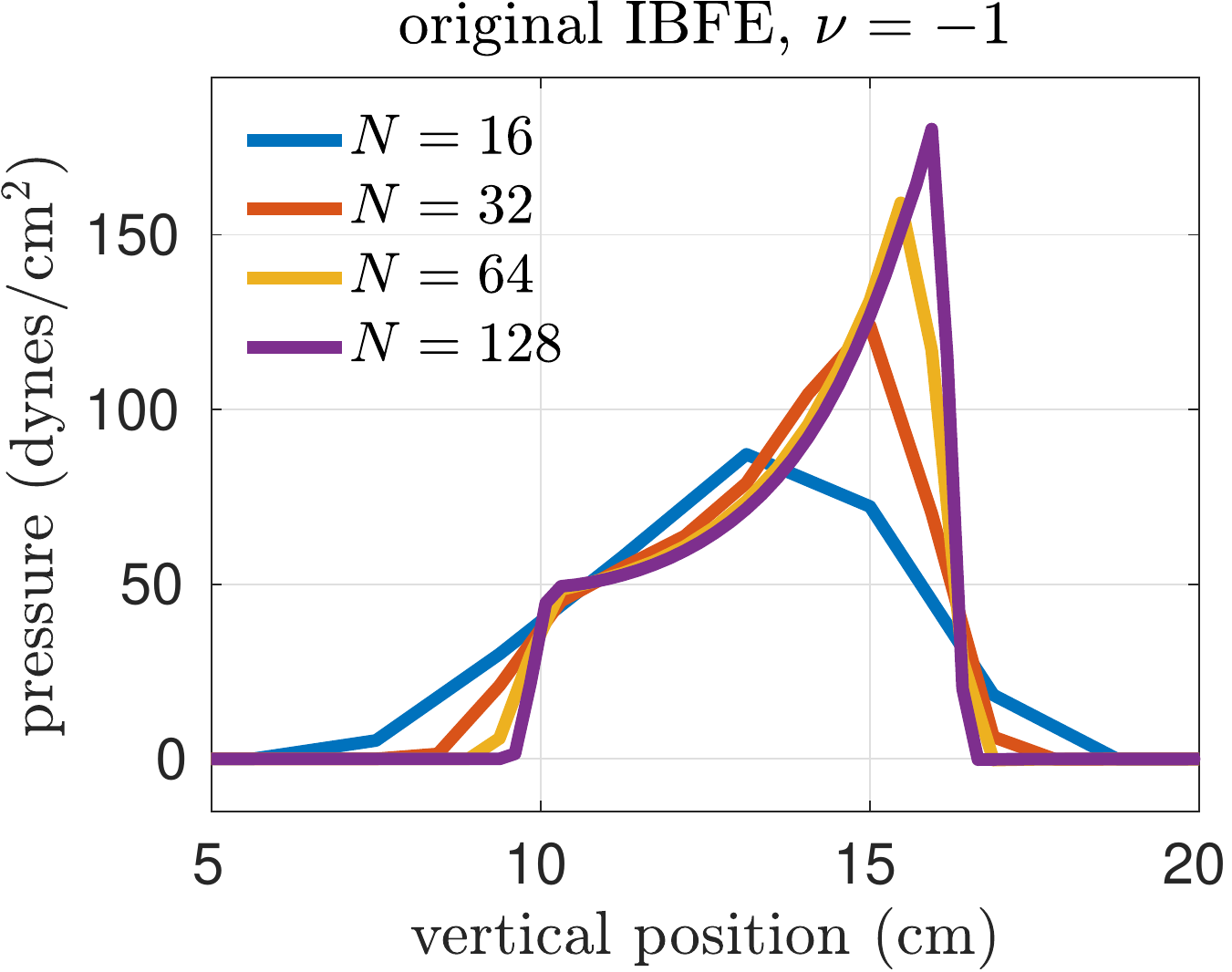}
\includegraphics[scale=0.4,trim=0 -62 0 0]{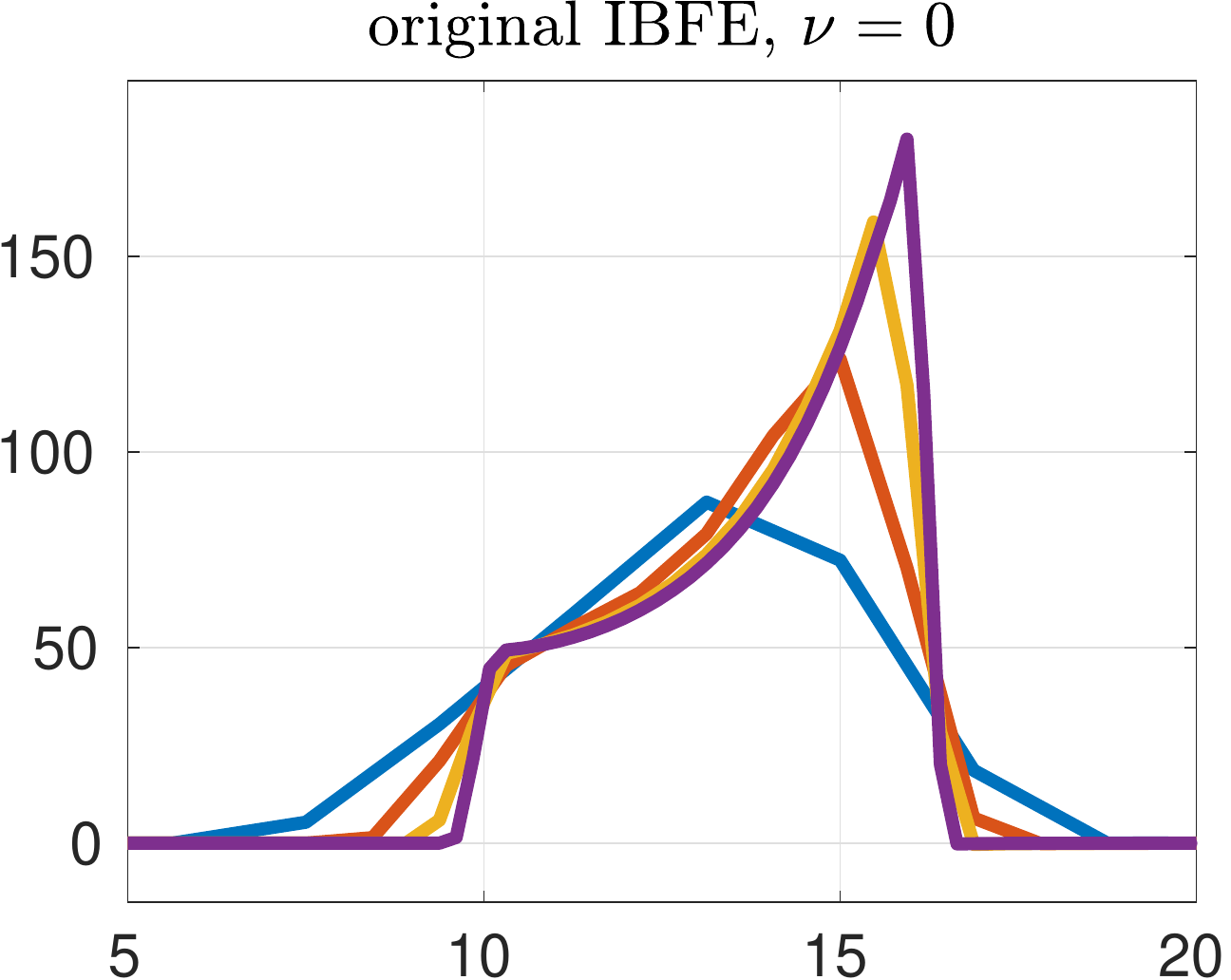}
\includegraphics[scale=0.4,trim=0 -62 0 0]{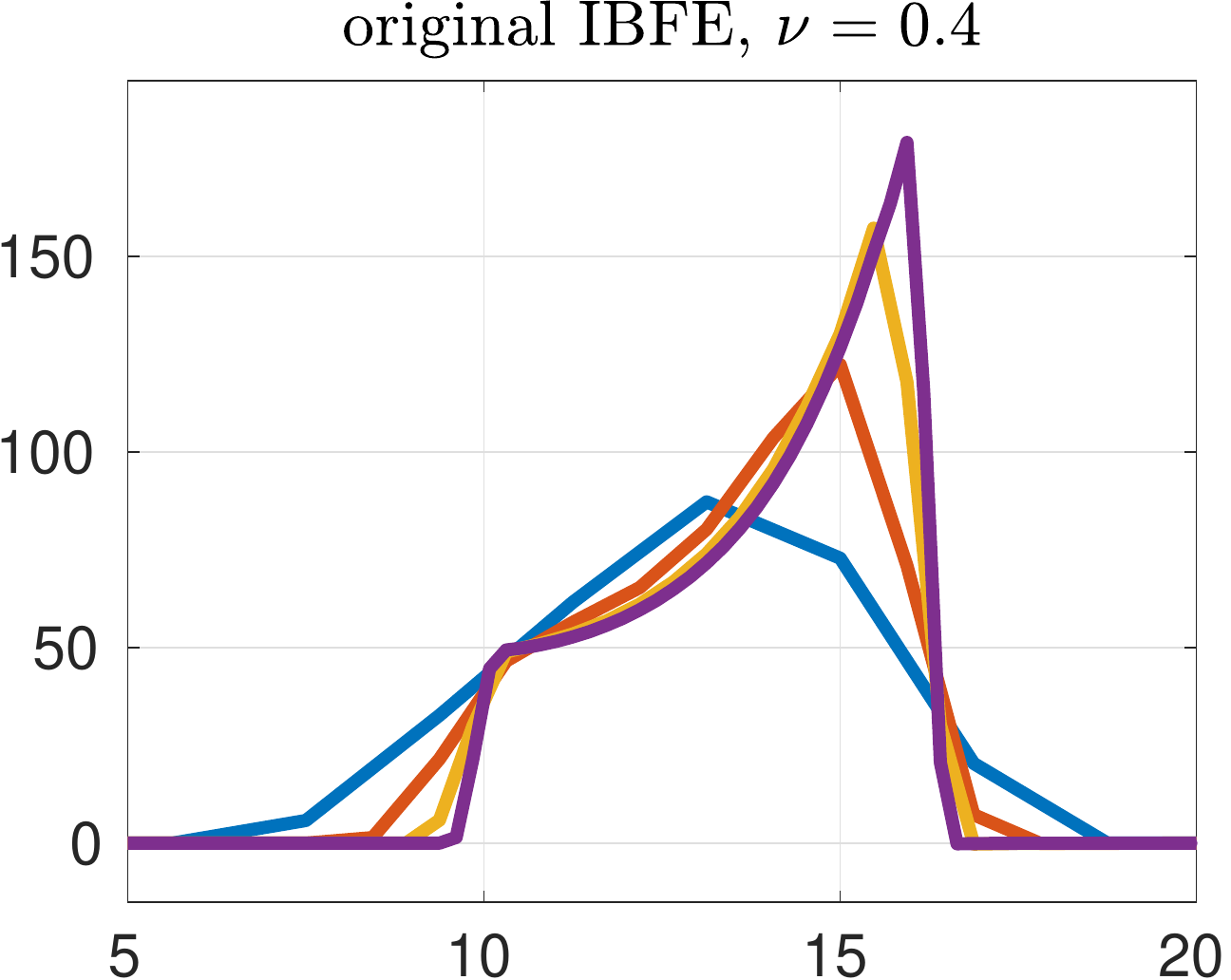} \\ 
\includegraphics[scale=0.4,trim=0 0 0 0]{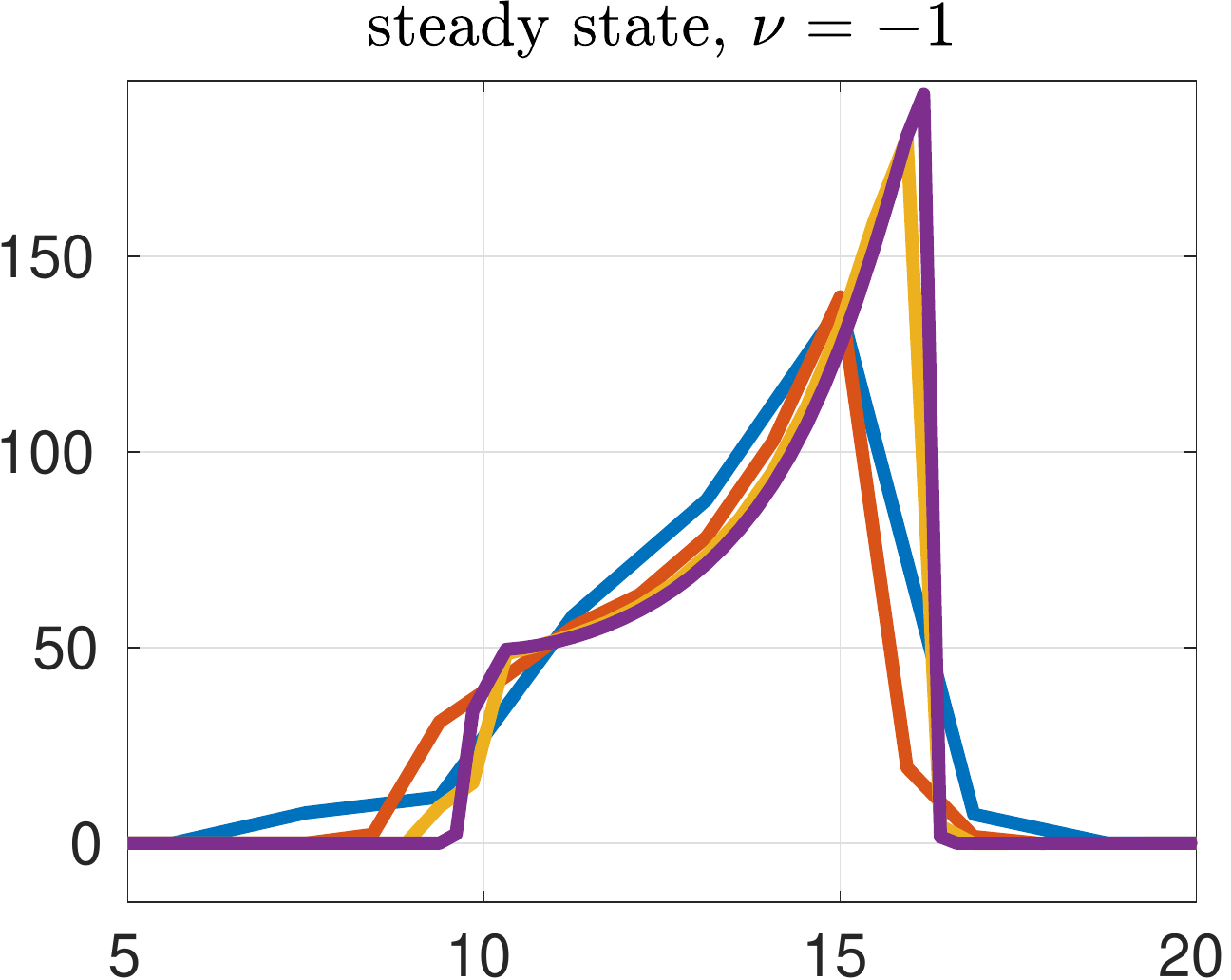}
\includegraphics[scale=0.4,trim=0 0 0 0]{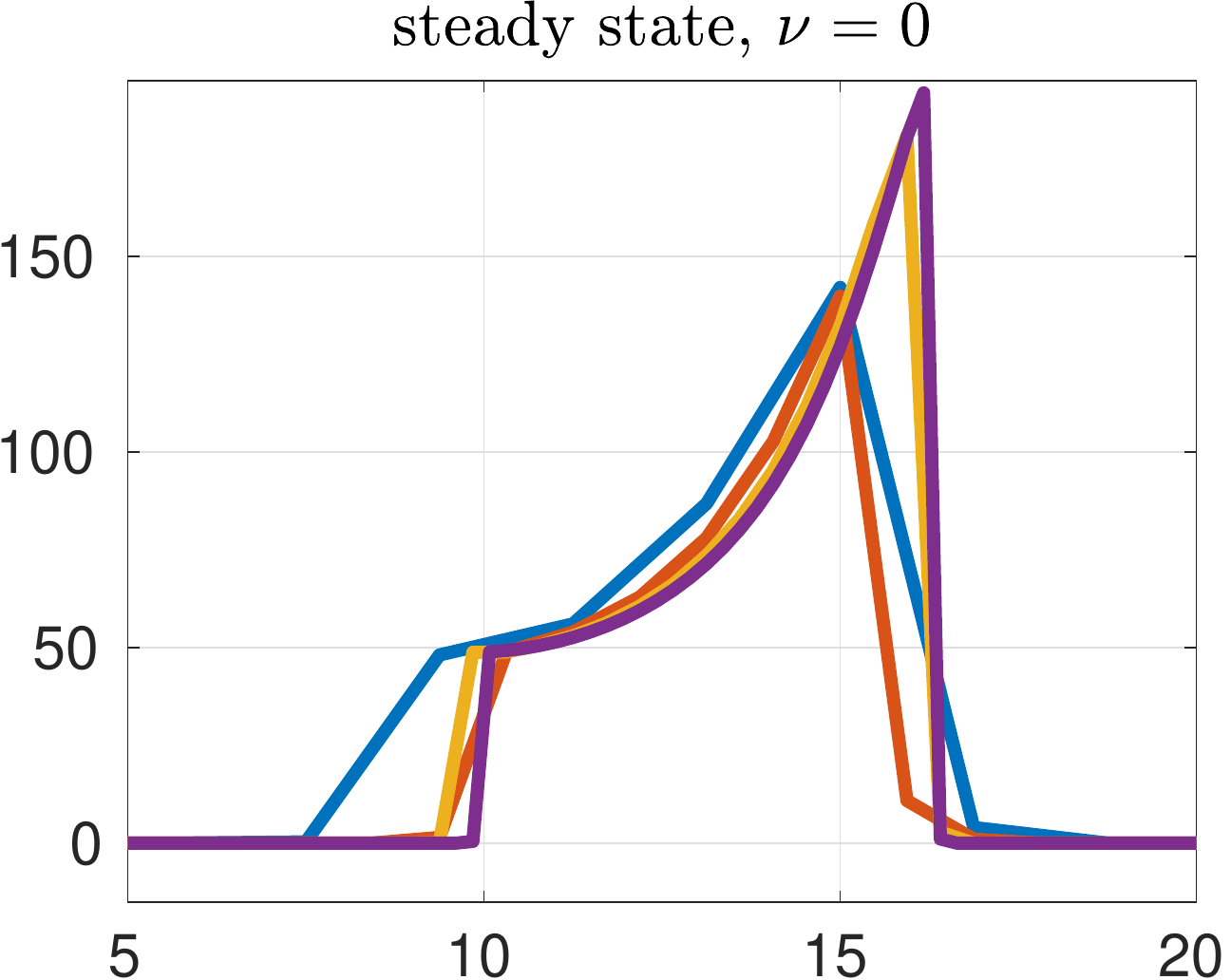}
\includegraphics[scale=0.4,trim=0 0 0 0]{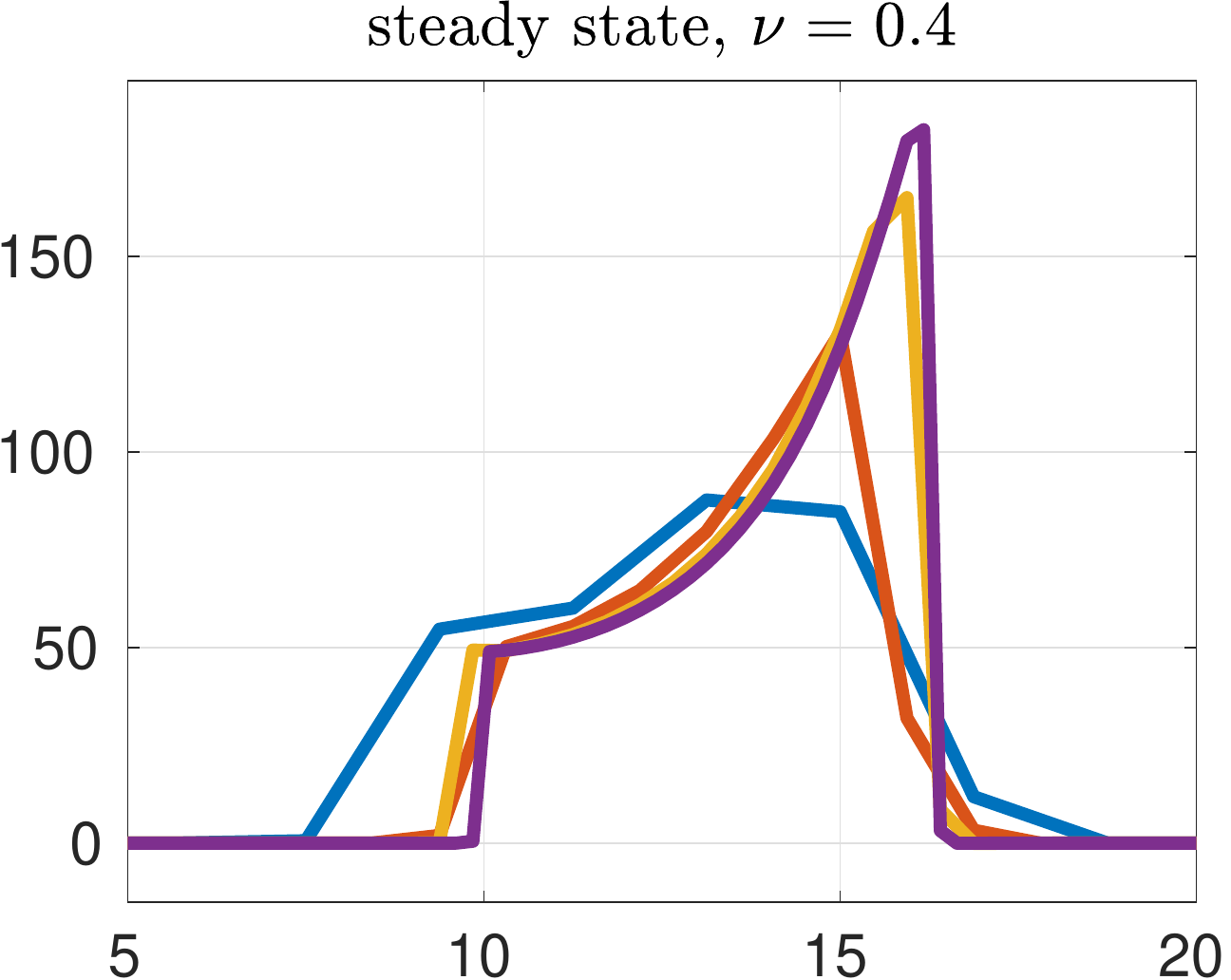} 
\caption{Slices of the pressure field vertically down the center of the block, at $X_1 = 15$. 
Results are displayed for different values of $\nu$.  
Results from the original IBFE method are on the top row, and results from the sharp interface method with the steady state formulation for $\varphi$ are on the bottom row.}
\label{fig:block4}
\end{center}
\end{figure}

Figure \ref{fig:block4} examines a slice of the pressure field taken vertically at $X_2 = 15$ mm, for different mesh sizes and different values for $\nu$.  
We only display slices for the sharp interface method with the steady state formulation for $\varphi$ since results for the diffusion formulation look similar.

\begin{figure}[h!]
\begin{center}
\includegraphics[scale=0.5,trim=0 0 -20 0]{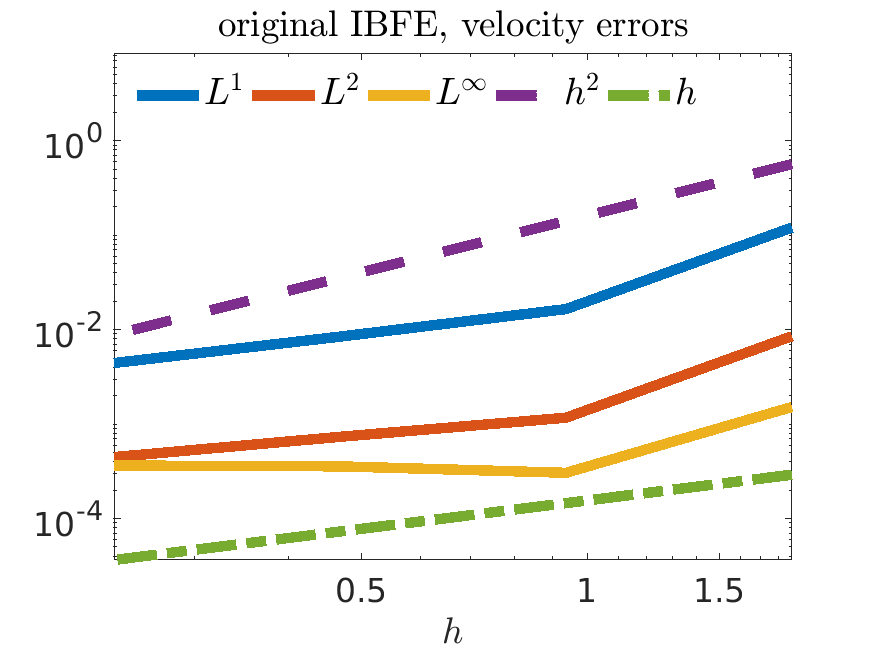}
\includegraphics[scale=0.5,trim=0 0 0 0]{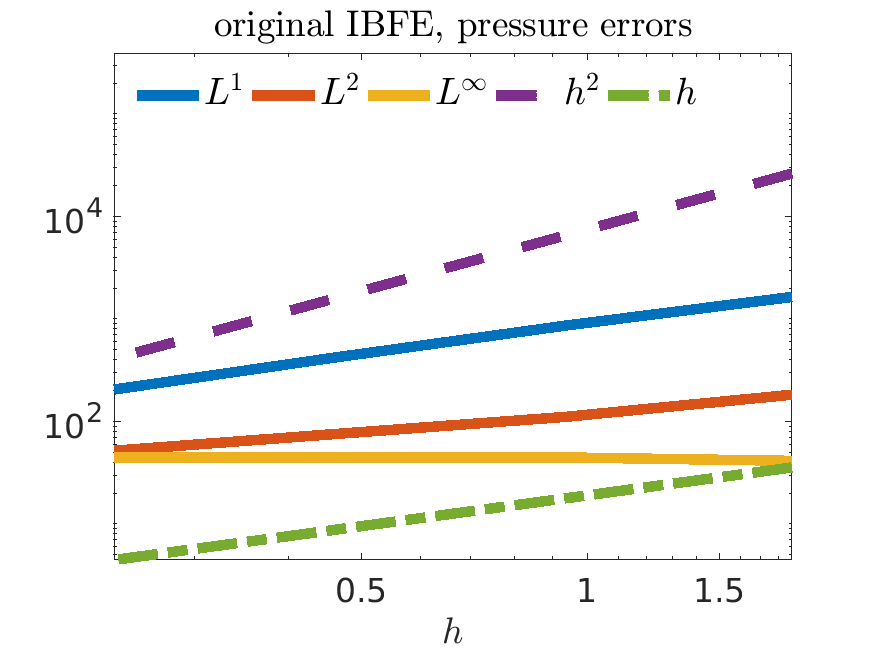}
\caption{Approximate absolute velocity and pressure errors for the compressed block test with the original IBFE method, the continuous loading pressure, $\nu = 0$. There is no pointwise convergence of the velocity and pressure.}
\label{fig:block_errors1}
\end{center}
\end{figure}

\begin{figure}[h!]
\begin{center}
\includegraphics[scale=0.5,trim=0 0 -20 0]{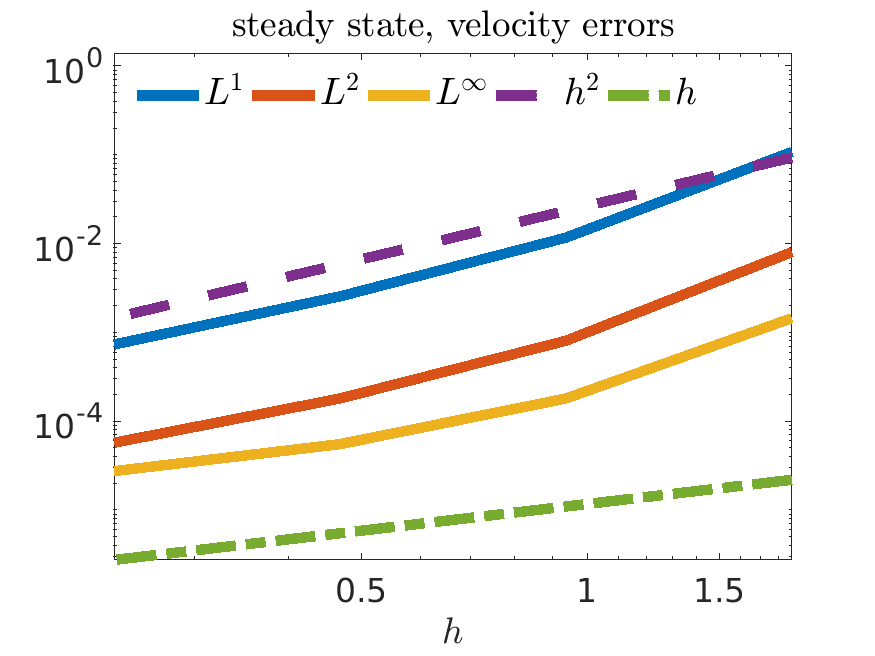}
\includegraphics[scale=0.5,trim=0 0 0 0]{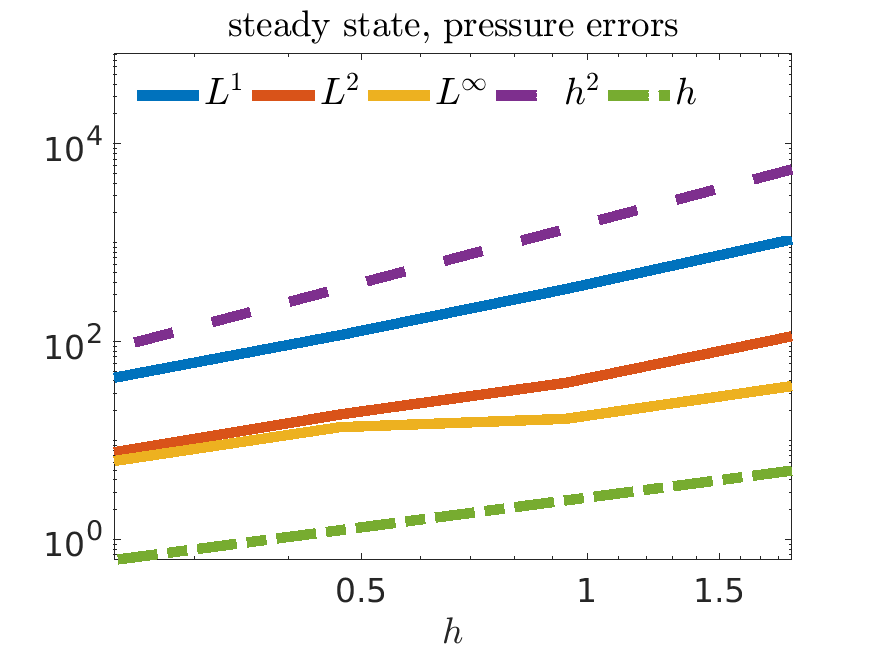}
\caption{Approximate absolute velocity and pressure errors for the compressed block test with the sharp interface method, the continuous loading pressure, and $\nu = 0$.   The steady state formulation is used to determine $\varphi$.  The pressure and velocity appear to converge in all norms.}
\label{fig:block_errors2}
\end{center}
\end{figure}

\revtwo{Figures \ref{fig:block_errors1} and \ref{fig:block_errors2}} display approximate absolute pressure and velocity errors, computed on the finite element mesh.  
The errors are computed as the difference between coarse and fine solutions at consecutive levels of refinement.  
For the original IBFE method, Figure \ref{fig:block_errors1}, the $L^1$ error for the velocity appears to be converging with a rate of 1 and the $L^\infty$ error is not converging.  
The $L^2$ converges with a rate less than 1.  The pressure errors show a similar trend.  
\revtwo{The approximate errors for the sharp interface method, shown in Figure \ref{fig:block_errors2}, are improved. 
All velocity and pressure errors converge with approximately a rate of 1.  Results are similar when using the diffusion equation for $\varphi$ with either $\gamma = 1$ or $\gamma = h$.}

%
%\begin{figure}[h!]
%\begin{center}
%\includegraphics[scale=0.5,trim=0 0 -20 0]{IMAGES/diff_block_velocity-errors}
%\includegraphics[scale=0.5,trim=0 0 0 0]{IMAGES/diff_block_pressure-errors}
%\caption{Approximate absolute velocity and pressure errors for the compressed block test with the sharp interface method, the continuous loading pressure, and $\nu = 0$.  
%The diffusion formulation with $\gamma = 1$ is used to determine $\varphi$.  The pressure and velocity appear to converge in all norms.}
%\label{fig:block_errors3}
%\end{center}
%\end{figure}
%
%\begin{figure}[h!]
%\begin{center}
%\includegraphics[scale=0.5,trim=0 0 -20 0]{IMAGES/diff_dx_block_velocity-errors}
%\includegraphics[scale=0.5,trim=0 0 0 0]{IMAGES/diff_dx_block_pressure-errors}
%\caption{Approximate absolute velocity and pressure errors for the compressed block test with the sharp interface method, the continuous loading pressure, and $\nu = 0$.  
%The diffusion formulation with $\gamma = h$ is used to determine $\varphi$.  The pressure and velocity appear to converge in all norms.}
%\label{fig:block_errors4}
%\end{center}
%\end{figure}

\clearpage
\subsection{\revtwo{Actively contracting thick torus}}

In our final example, we apply the method to a thick, actively contracting torus which creates a nonzero velocity field. This problem is inspired by the work of McQueen and Peskin \cite{McQueen89}.  The constitutive model contains an active stress component describing a wave of contraction; this wave pushes fluid within the torus, creating a peristaltic pump.

The geometrical parameters are $R_1$, $R_2$, and $w$, and the thick toroidal geometry is parametrized with respect to the coordinates ${\bs s} =  (s_1, s_2, s_3) = [-\pi R_1, \pi R_1) \times [0,w] \times [-\pi R_2, \pi R_2)$.  The reference coordinates ${\bs X} = (X_1, X_2, X_3)$ are defined as:
\begin{align*}
X_1 &= R_2 \sin (s_3/R_2) + (R_1+s_2) \cos(s_1/R_1) \sin(s_3/R_2), \\
X_2 &= (R_1+s_2) \sin(s_1/R_1), \\
X_3 &= R_2 \cos (s_3/R_2) + (R_1+s_2) \cos(s_1/R_1) \cos(s_3/R_2).
\end{align*}
We identify this parametrization by a function $\mathcal{G}$, i.e. ${\bs X} = \mathcal{G}({\bs s})$.  The parameter $w = 0.125$ mm is the thickness of the wall, built from inner and outer toroidal surfaces.  The parameters $R_1 = 0.25$ mm and $R_2 = 1$ mm are the minor and major radii of the inner toroidal surface respectively.

The elastic stress is a sum of two components, a passive stress and an active stress which depends on a fiber vector field ${\bs f}_0$ specified in the reference configuration:
\begin{align*}
\mathbb{P}^\text{e}({\bs X},t) = \mu_e \left(\mathbb{F} - \mathbb{F}^{-T} \right) + T(\mathcal{G}^{-1}({\bs X}),t)\, \mathbb{F}\,{\bs f}_0 \otimes {\bs f}_0.
\end{align*}

\begin{figure}[h!]
\begin{center}
\includegraphics[scale=0.225,trim=0 0 0 0]{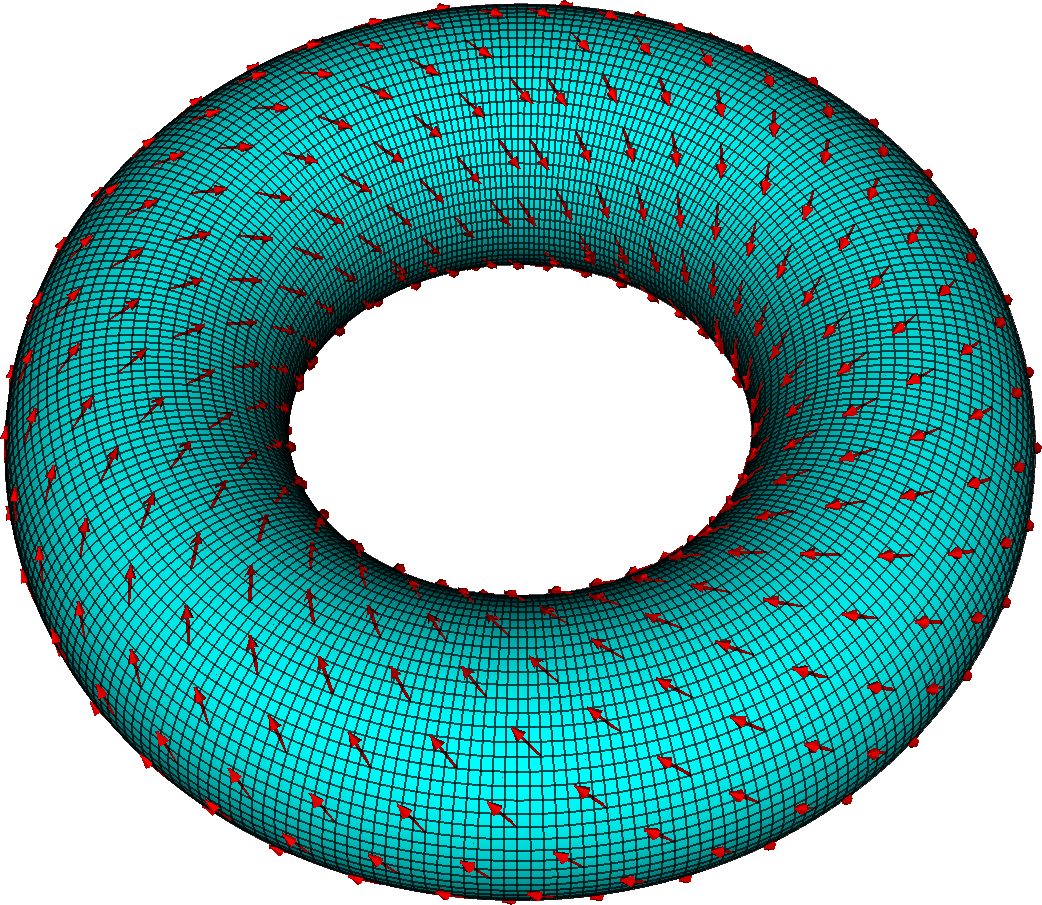}
\caption{The mesh for the torus, with the fiber vector field ${\bs f}_0$ superimposed.}
\label{fig:torus1}
\end{center}
\end{figure} 
\noindent
The size, strength, and location of the contraction is determined by a tension function $T = T({\bs s}, t)$.  The support of this function has measure $\ell_\text{c}$ in the third component $s_3$ and translates with velocity $\nu_\text{c}$.  Explicitly, the tension function is defined as
\begin{align*}
T({\bs s},t) = 
\begin{cases}
T_\text{ramp}(t) \times \exp\left(\frac{(b-a)^2}{(2(s_3 + \nu_\text{c} t) - a - b)^2 - (b-a)^2} + 1\right) & \text{if } a < s_3 + \nu_\text{c} t < b \\
0 & \text{otherwise},
\end{cases}
\end{align*}
with $a = -\frac{1}{2} \ell_\text{c}$ and $b = \frac{1}{2}\ell_\text{c}$.  Notice the tension function smoothly decays to zero.  The maximum of the tension is ramped up linearly for an amount of time $t_\text{ramp}$ to a value $T_\text{max}$:
\begin{align*}
T_\text{ramp}(t) = 
\begin{cases}
\frac{t}{t_\text{ramp}} \times T_\text{max} & \text{if } t < t_\text{ramp},\\
T_\text{max} & \text{otherwise}.
\end{cases}
\end{align*}
The finite element mesh for the torus, with the fiber field ${\bs f}_0$ superimposed, is shown in Figure \ref{fig:torus1}.  The angle of the fiber field rotates from $60^\text{o}$ to $-60^\text{o}$ from the outer surface to the inner surface.  The mesh is centered in the computational domain defined to be $\Omega = [-2L, 2L] \times [-L, L] \times [-2L, 2L]$ with $L = 0.75$.  This domain is discretized with $2N$ cells in the $x$ and $z$ directions and $N$ cells in the $y$ direction, with $N = 64$.  The finite element mesh is defined so $M_\text{fac}$ is approximately 2.  Boundary conditions for the computational domain are set to hold the pressure at zero.  The timestep size is $\Delta t = 0.0025 \times h$, with $h = L/N$. 

\begin{figure}
\begin{center}
\includegraphics[scale=0.125,trim=0 0 -20 0]{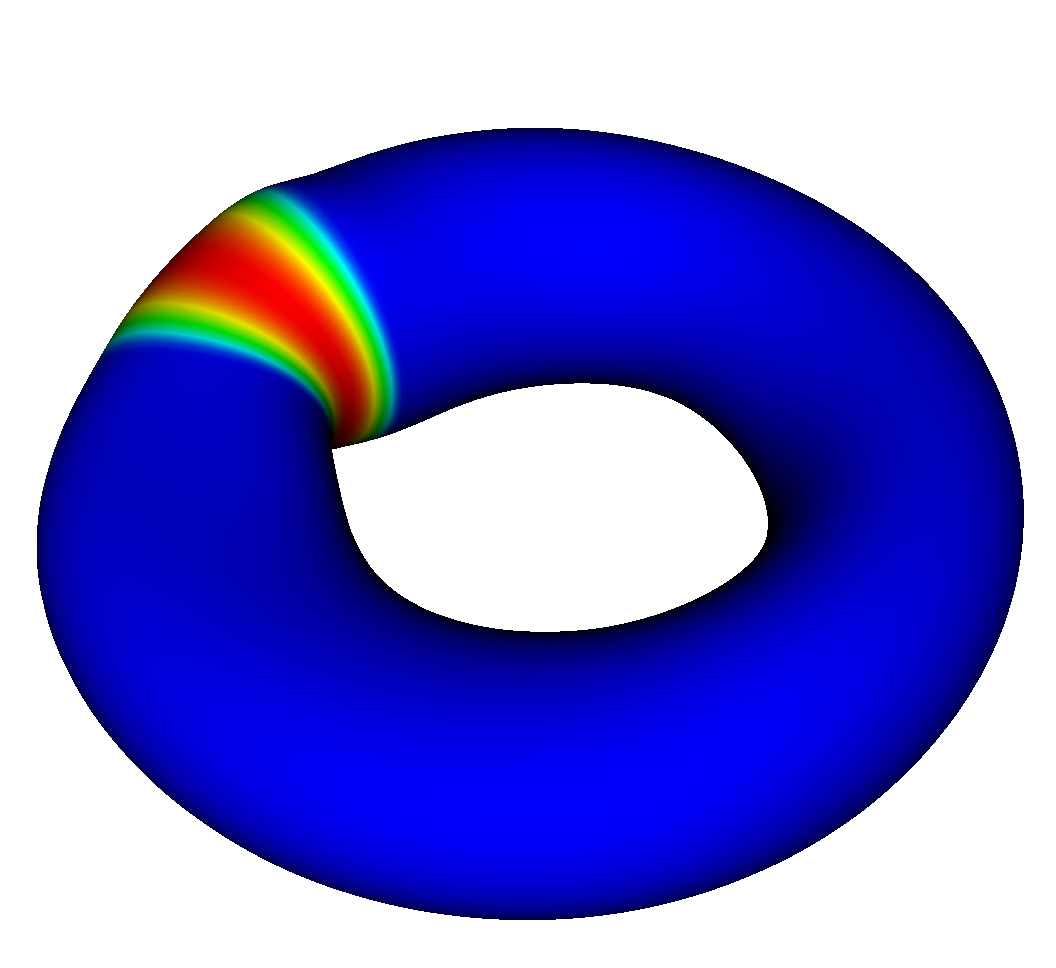}
\includegraphics[scale=0.125,trim=0 0 -20 0]{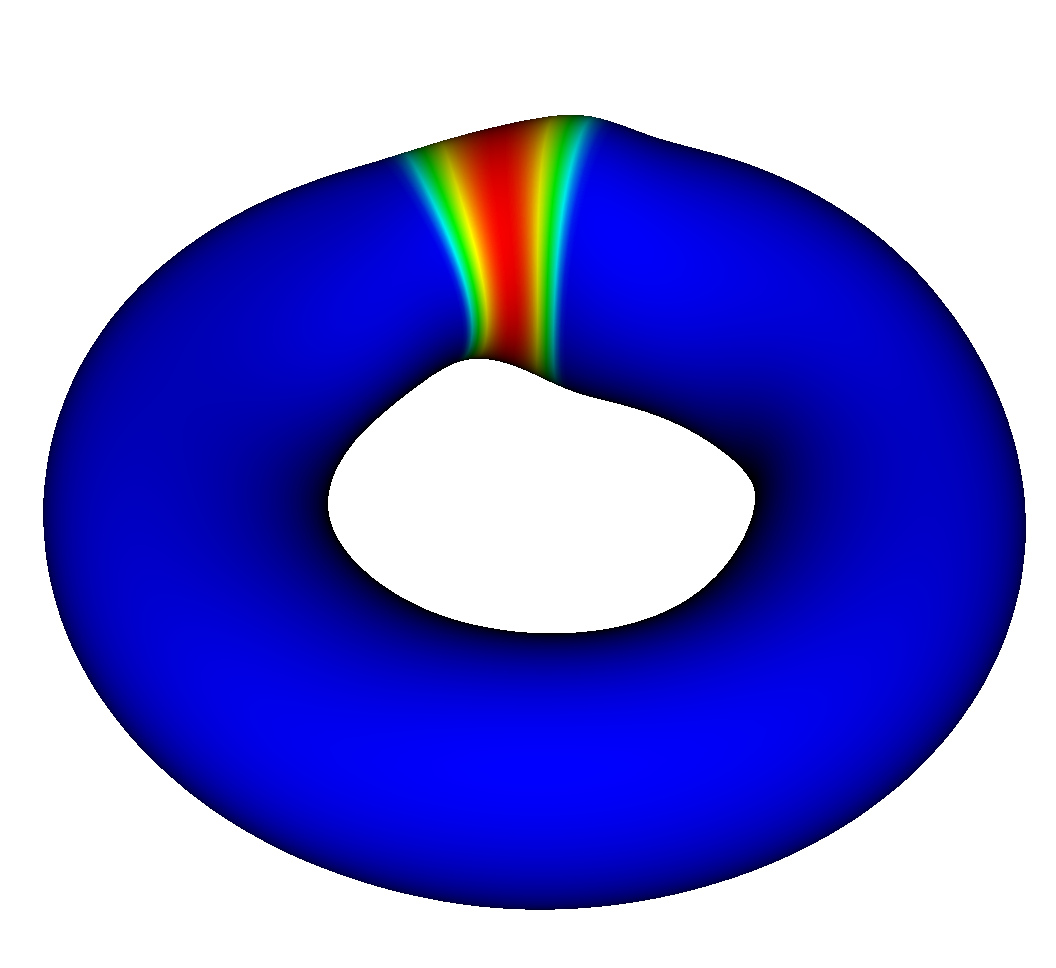}
\includegraphics[scale=0.125,trim=0 0 -20 0]{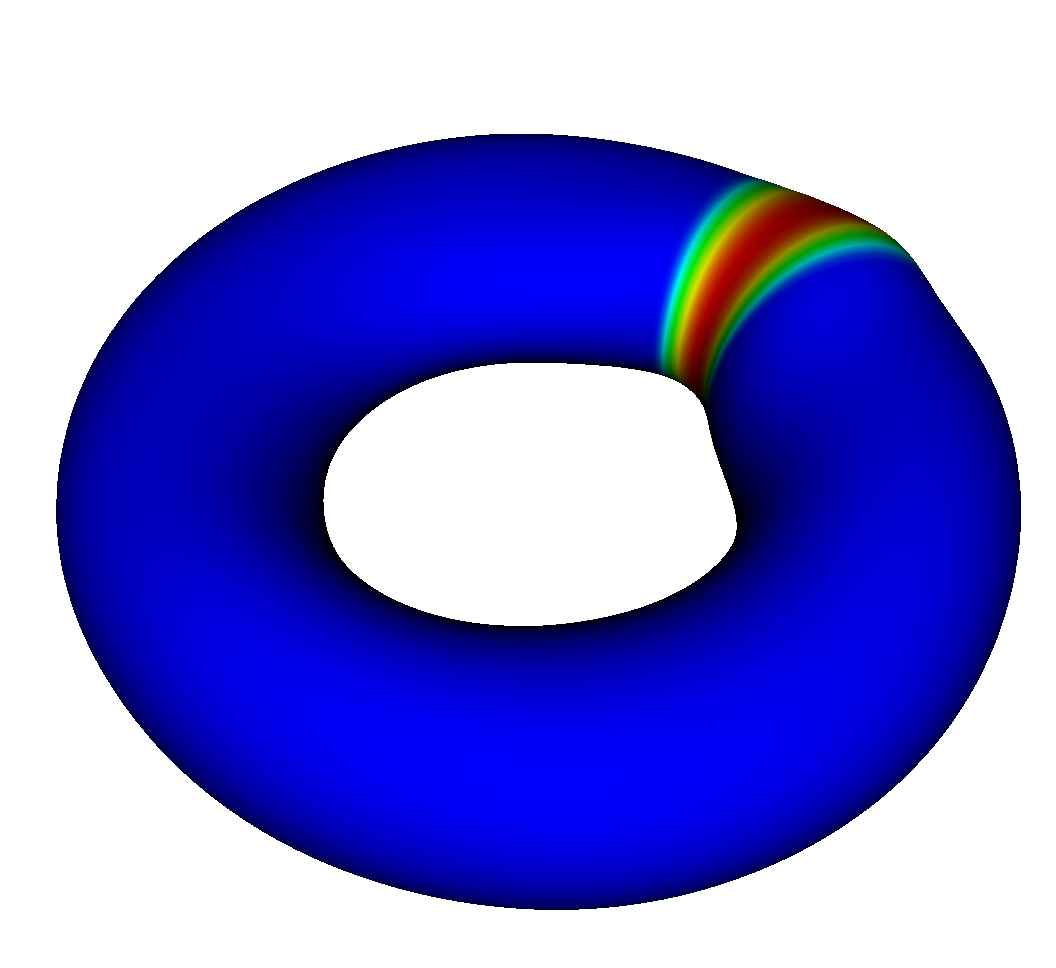} \\
\includegraphics[scale=0.125,trim=0 0 -20 0]{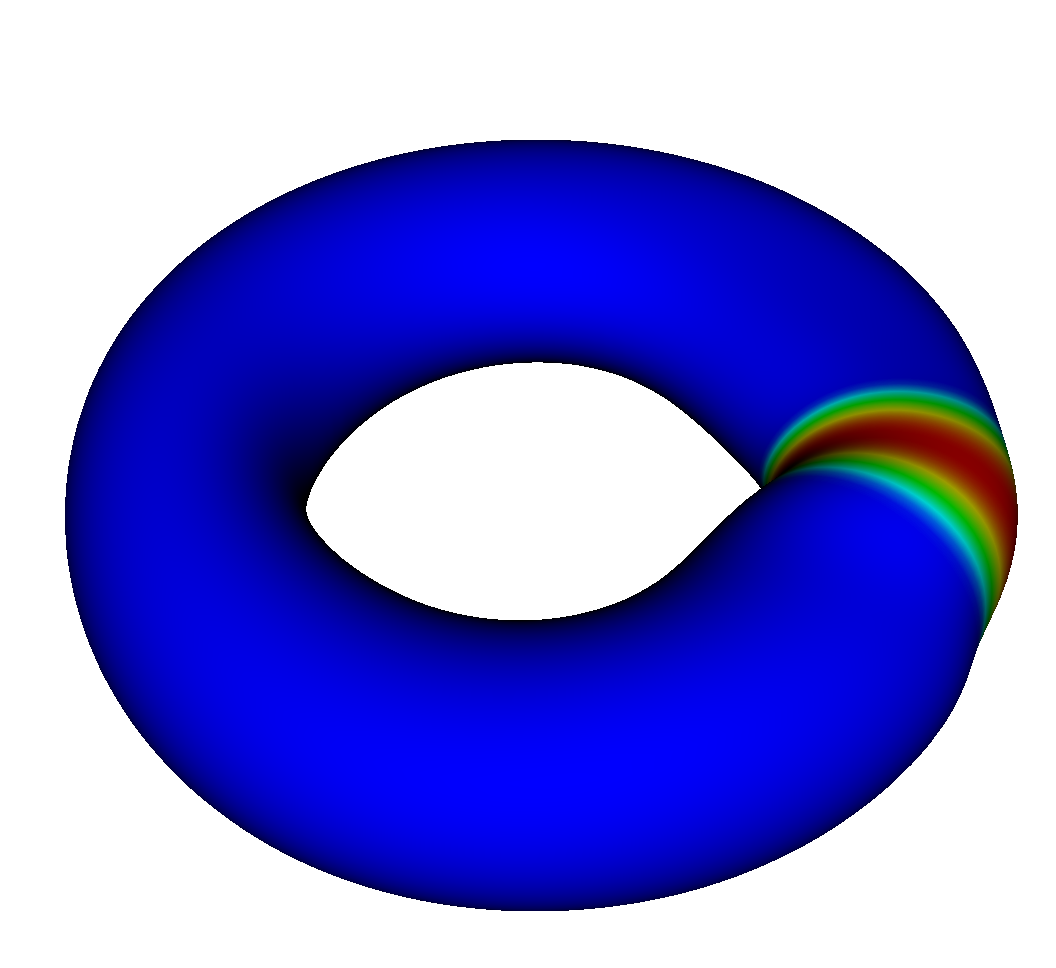}
\includegraphics[scale=0.125,trim=0 0 -20 0]{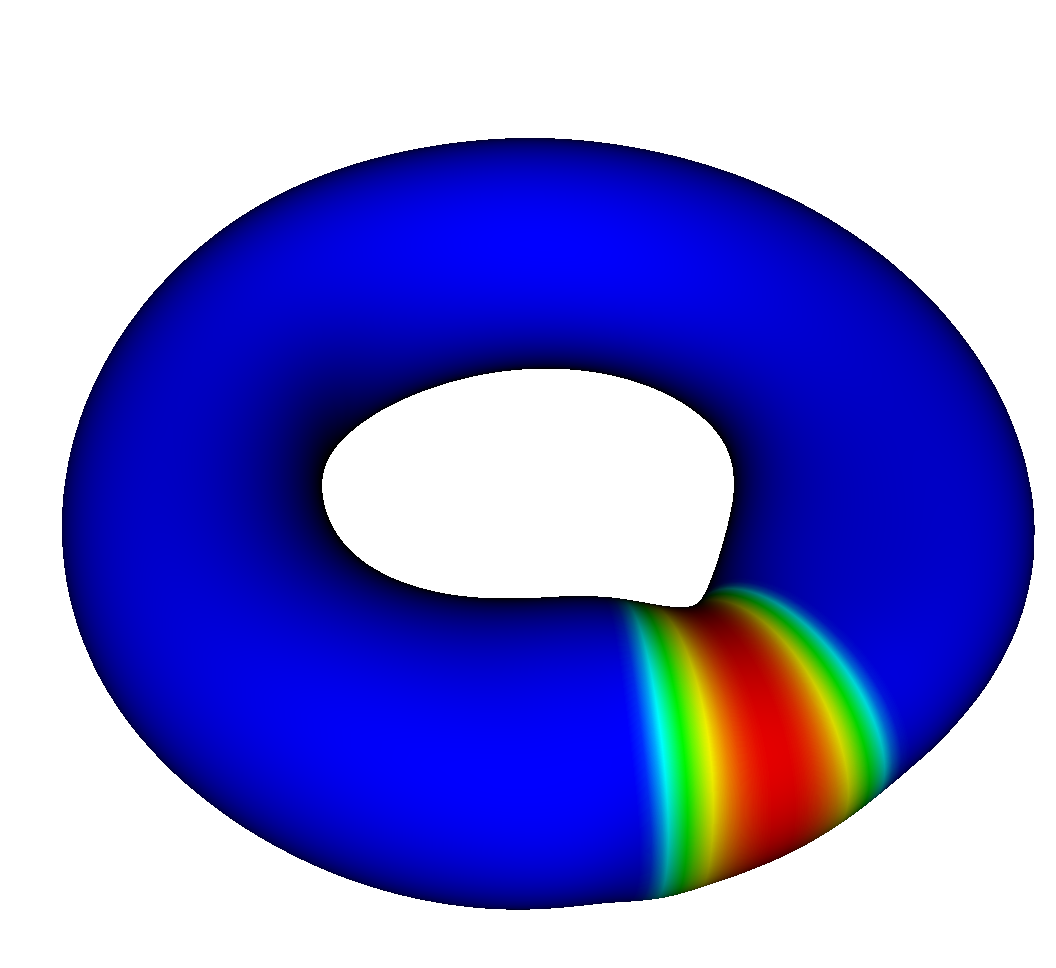}
\includegraphics[scale=0.125,trim=0 0 -20 0]{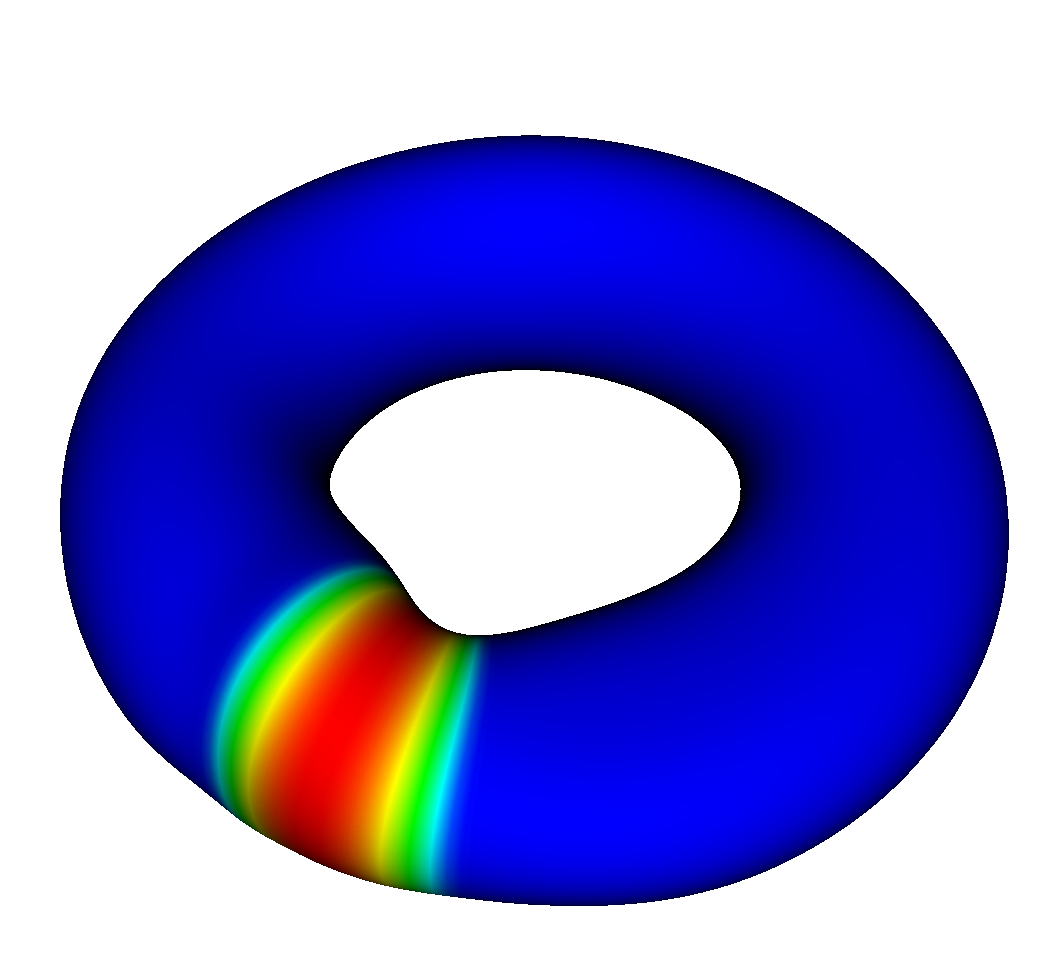} \\\caption{A visualization of contraction.  Time increases from left--to--right and top--to--bottom.  The color indicates the value of the tension function.}
\label{fig:torus2}
\end{center}
\end{figure}

\begin{figure}
\begin{center}
\includegraphics[scale=0.25,trim=0 0 -250 0]{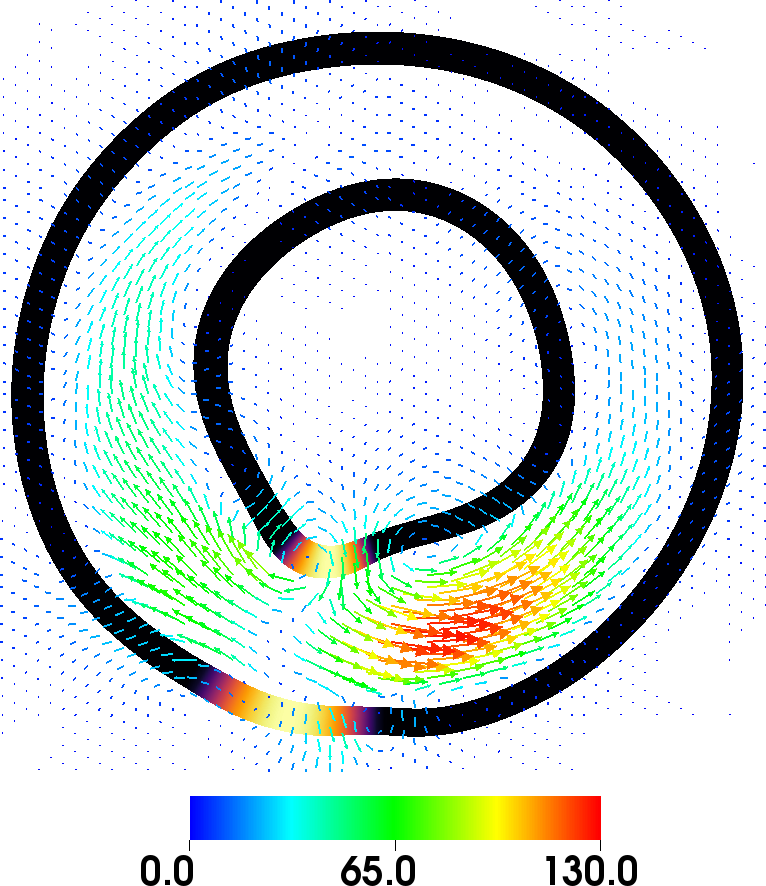}
\includegraphics[scale=0.25]{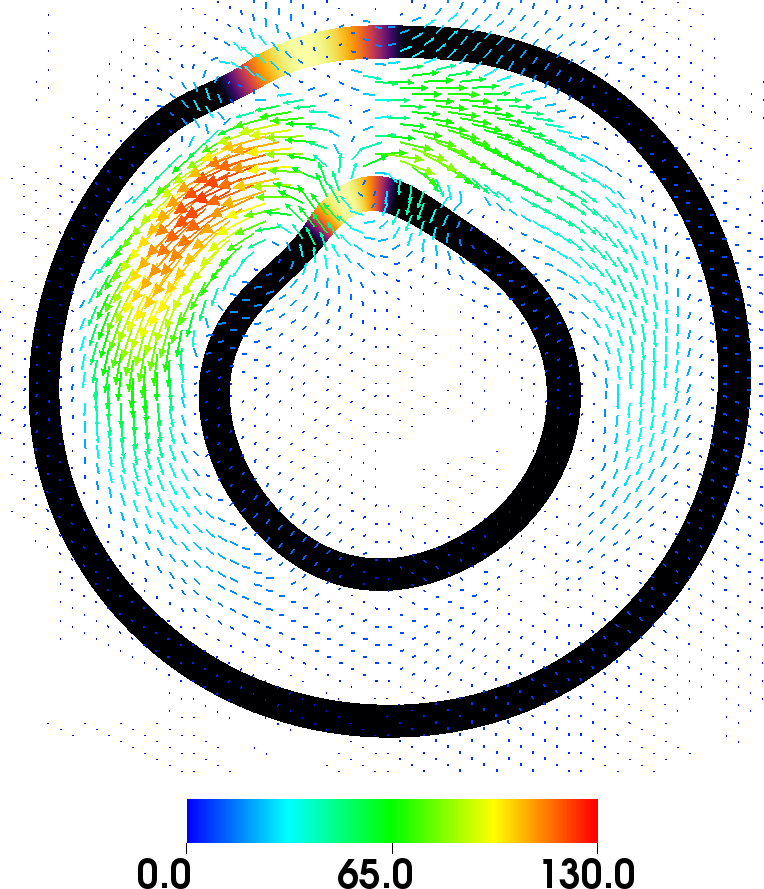}
\caption{A slice of the velocity field for the original IBFE method shown at two different snapshots in time.  The torus is also shown, with the color on the torus indicating the value of the active tension function $T$.}
\label{fig:torus3}
\end{center}
\end{figure}

\begin{figure}
\begin{center}
\includegraphics[scale=0.25,trim=0 0 -250 0]{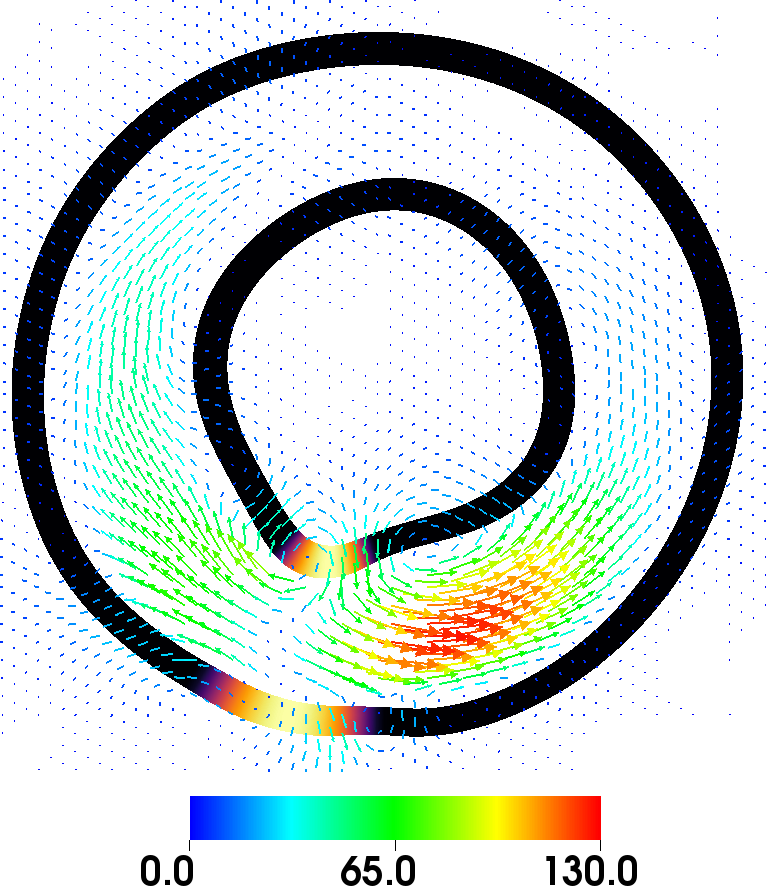}
\includegraphics[scale=0.25]{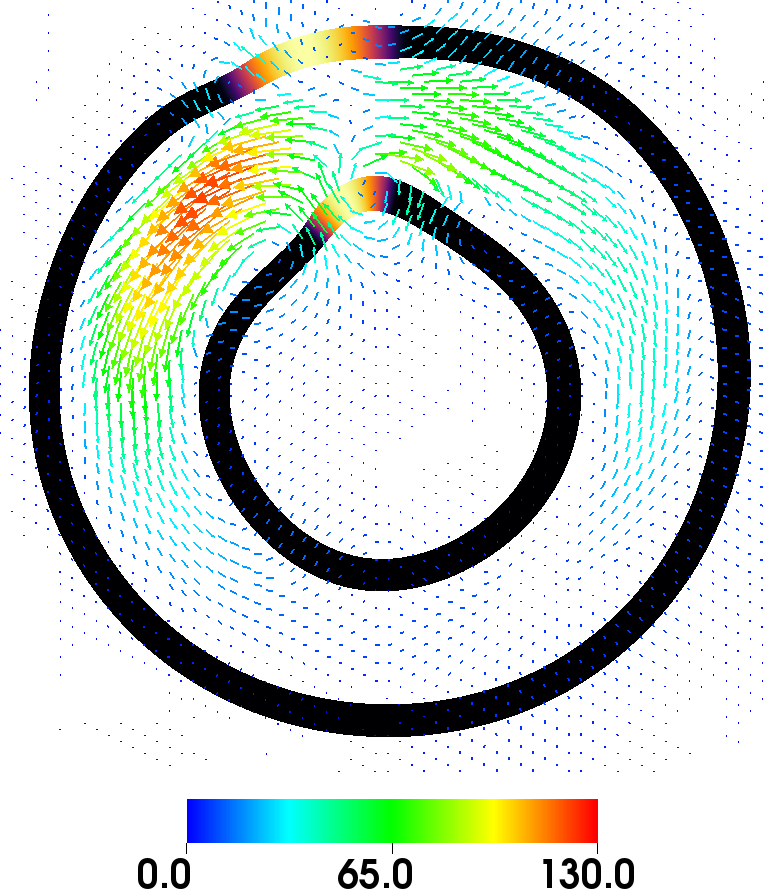}
\caption{A slice of the velocity field for the sharp interface method shown at two different snapshots in time.  The torus is also shown, with the color on the torus indicating the value of the active tension function $T$.}
\label{fig:torus4}
\end{center}
\end{figure}

The displacement of the torus during contraction is shown in Figure \ref{fig:torus2}, with time increasing from left--to--right, top--to--bottom.  The color indicates the value of the tension function.  The torus is not pre--stressed, but the active stress term results in large pressure discontinuities at the fluid--structure interface and induces a nonzero velocity field. Figures \ref{fig:torus3} and \ref{fig:torus4} show the velocity field on a slice through the center of the displaced torus for the original IBFE method and the sharp interface method respectively. The Reynolds number is approximately 195, and both methods qualitatively produce the same velocity field.  Notice the region of forward flow directly in front of the contracted part of the torus.

\begin{figure}
\begin{center}
\includegraphics[scale=0.25,trim=0 0 -250 0]{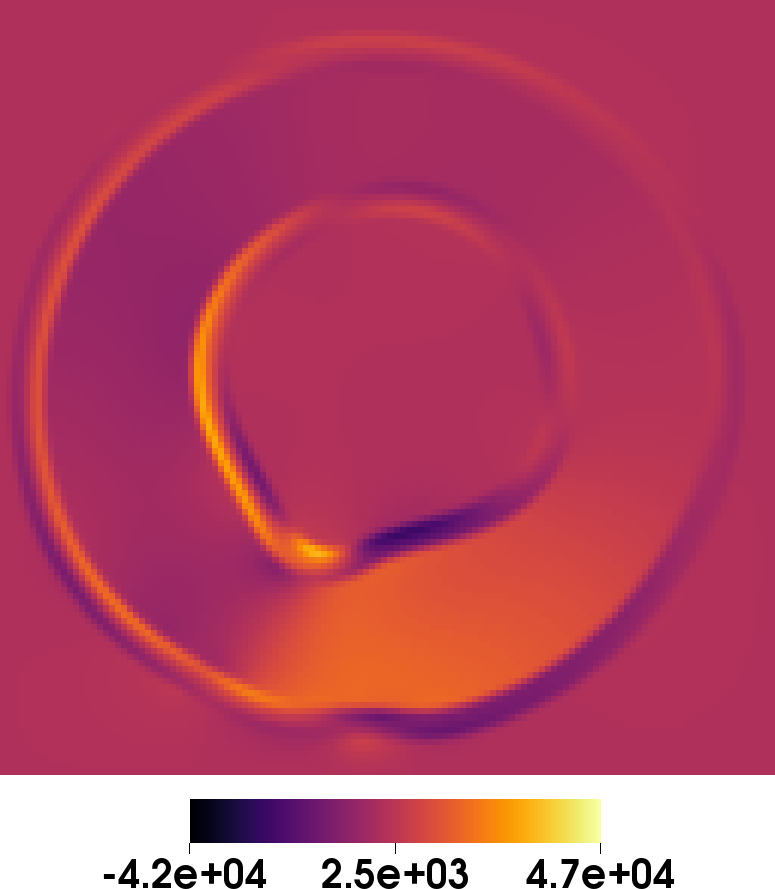}
\includegraphics[scale=0.25]{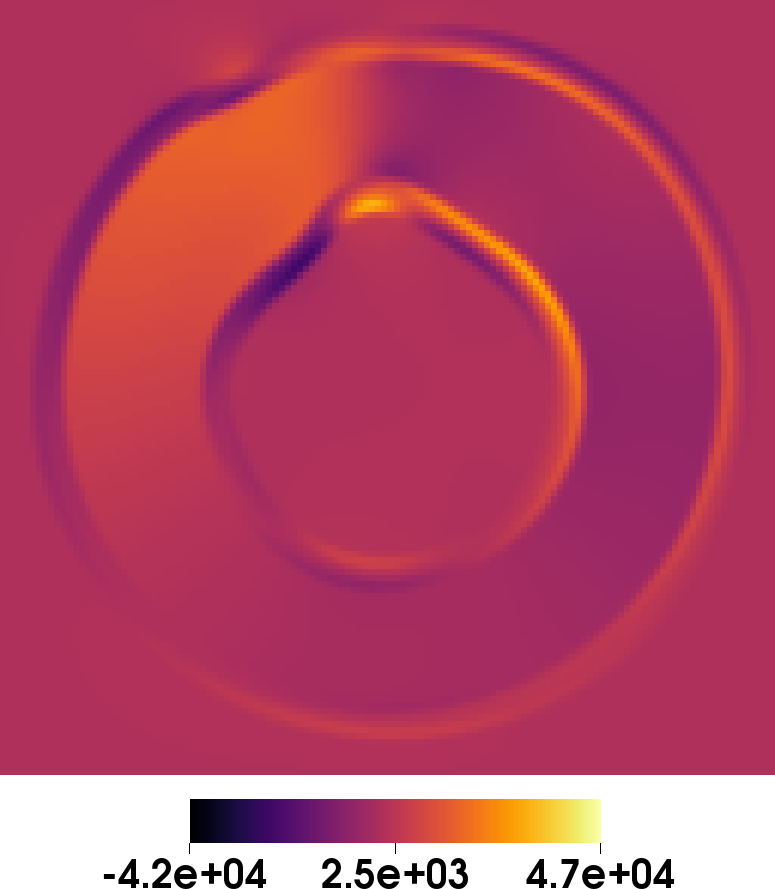}
\caption{A slice of the pressure for the original IBFE method at two difference snapshots in time.  Note that the pressure is regularized at the fluid--structure interface.}
\label{fig:torus5}
\end{center}
\end{figure}

\begin{figure}
\begin{center}
\includegraphics[scale=0.25,trim=0 0 -250 0]{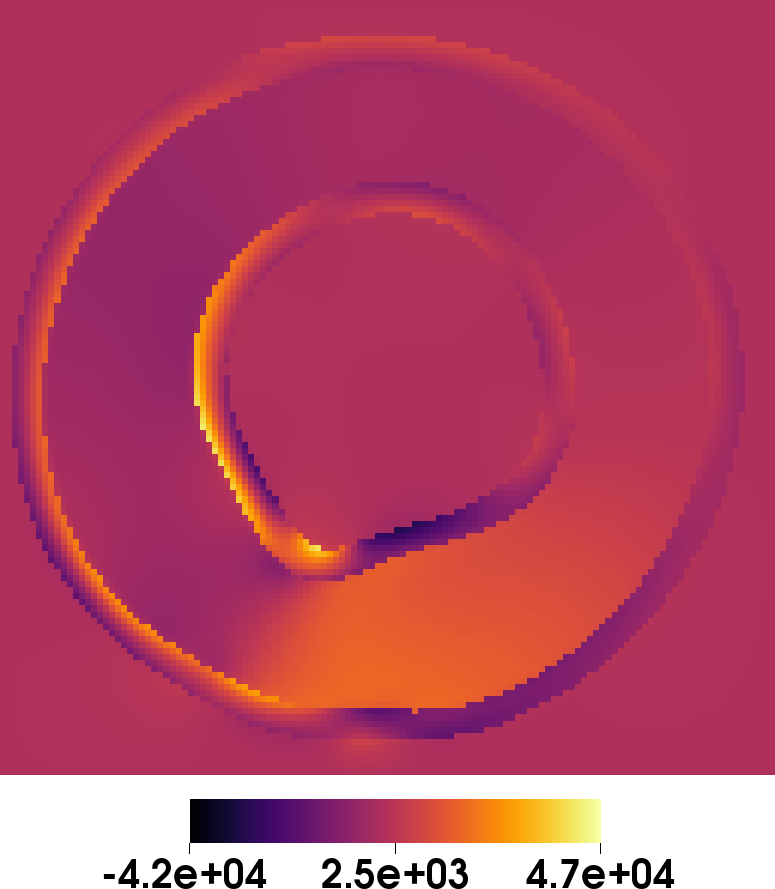}
\includegraphics[scale=0.25]{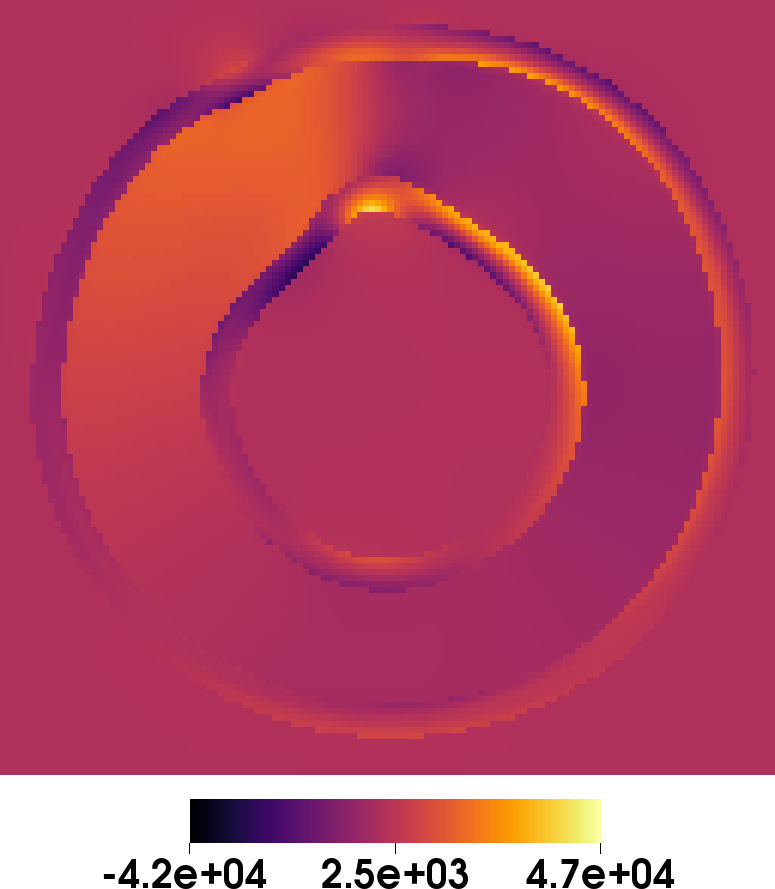}
\caption{A slice of the pressure for the sharp interface method at difference snapshots in time.  Our method resolves the pressure discontinuities at the fluid-structure interface which are generated by the contraction.}
\label{fig:torus6}
\end{center}
\end{figure}

Figures \ref{fig:torus5} and \ref{fig:torus6} show a slice of the pressure field from the original IBFE method and the sharp interface method respectively.  The active stress term creates large pressures in the wall of the torus, leading to pressure discontinuities at the fluid--structure interface.  These pressure jumps are smoothed out with the original approach.  The sharp interface method resolves these discontinuities.

\clearpage
\section{Conclusion}

This paper describes a numerical method for resolving pressure discontinuities in an immersed boundary finite element formulation.  
This method applies to immersed solids of codimension 0 which may undergo large deformations.
The pressure is split into two pressure--like fields, one of which is continuous and defined on the entire domain, and the other is defined only on the finite element mesh.  
Calculation of the pressure--like field defined on the finite element mesh is done by solving either a steady state harmonic problem or a diffusion equation, with boundary conditions that depend on the normal component of the elastic traction. 

The method is tested on \revtwo{four} examples.
One example statically contains pressure discontinuities, and the other \revtwo{three} involve large deformations.  
In each test, the method resolves the pressure discontinuities on the fluid--structure interface.  
In the thick ring examples with analytic solutions, we see improved convergence for the pressure, velocity, displacement, and elastic stress, and a substantial decrease in the magnitude of errors for modest numbers of degrees of freedom.  
Further, the inflating ring example demonstrates that the method helps to enforce solid incompressibility, at least for coarse finite element meshes.
In the more complex compressed block example, the method  more sharply resolves the pressure field for both smooth and discontinuous loading forces.  Results with the discontinuous loading force suggest that care must be taken in managing discontinuous forcing on the fluid--structure interface.  
In the case of a smooth loading force, extrapolated errors for both the pressure and velocity appear to be converging.  This example highlights the method's flexibility in dealing with additional surface forces and volumetric energies with modestly large stabilization parameters \cite{Vadala18}. \revtwo{In the final example, we demonstrated the method works in three dimensions with transient flow dynamics and a stress function which depends on a fiber vector field.  The active stress term creates large pressures in the toroidal wall, leading to pressure discontuities at the fluid--structure interface which are easily resolved with our approach.} 

As shown in the compressed block example, the proposed method appears to struggle with large volumetric stabilization parameters.  A related question is if this approach can serve as a replacement for the volumetric stabilization technique, given its performance in enforcing incompressible solid deformations for the inflating ring.  These questions are topics of future work.  

Additional future work includes the resolution of other jump conditions, including the normal derivative of the velocity at the fluid--structure interface, possibly by solving higher order harmonic problems on the solid finite element mesh.  Further, \revtwo{a modification of the method may work for codimension 1 structures which correspond to the boundary of a codimension 0 structure.  A similar technique    may be possible if a volumetric mesh for the corresponding codimension 0 structure is available; in this case, a harmonic problem could be solved on the associated volumetric mesh to impose known jump conditions across the thin interface.}

\clearpage
\section{Acknowledgments}

Charles Puelz was supported in part by the Research Training Group in Modeling and Simulation funded by the National Science Foundation via grant RTG/DMS--1646339. Boyce Griffith was supported in part by  NSF DMS 1664645, NSF OAC 1652541, NSF OAC 1450327. The authors thank Charles Peskin for many insightful conversations.

\clearpage
\bibliographystyle{plain}
\bibliography{IB_refs}

\begin{thebibliography}{10}

\bibitem{Bedrossian10}
Jacob Bedrossian, James~H Von~Brecht, Siwei Zhu, Eftychios Sifakis, and
  Joseph~M Teran.
\newblock A second order virtual node method for elliptic problems with
  interfaces and irregular domains.
\newblock {\em {Journal of Computational Physics}}, 229(18):6405--6426, 2010.

\bibitem{Boffi08}
Daniele Boffi, Lucia Gastaldi, Luca Heltai, and Charles~S Peskin.
\newblock On the hyper-elastic formulation of the immersed boundary method.
\newblock {\em {Computer Methods in Applied Mechanics and Engineering}},
  197(25-28):2210--2231, 2008.

\bibitem{Griffith17}
Boyce E.~Griffith and Xiaoyu Luo.
\newblock Hybrid finite difference/finite element immersed boundary method.
\newblock {\em {International Journal for Numerical Methods in Biomedical
  Engineering}}, 33(12):e2888, 2017.

\bibitem{Lai01}
Ming-Chih Lai and Zhilin Li.
\newblock {A remark on jump conditions for the three-dimensional Navier-Stokes
  equations involving an immersed moving membrane}.
\newblock {\em {Applied Mathematics Letters}}, 14(2):149--154, 2001.

\bibitem{Lee03}
Long Lee and Randall~J LeVeque.
\newblock {An immersed interface method for incompressible Navier--Stokes
  equations}.
\newblock {\em {SIAM Journal on Scientific Computing}}, 25(3):832--856, 2003.

\bibitem{Leveque94}
Randall~J Leveque and Zhilin Li.
\newblock The immersed interface method for elliptic equations with
  discontinuous coefficients and singular sources.
\newblock {\em SIAM {J}ournal on {N}umerical {A}nalysis}, 31(4):1019--1044,
  1994.

\bibitem{Li06}
Zhilin Li and Kazufumi Ito.
\newblock {\em {The immersed interface method: numerical solutions of PDEs
  involving interfaces and irregular domains}}, volume~33.
\newblock {SIAM}, 2006.

\bibitem{Liu06}
Wing~Kam Liu, Yaling Liu, David Farrell, Lucy Zhang, X~Sheldon Wang, Yoshio
  Fukui, Neelesh Patankar, Yongjie Zhang, Chandrajit Bajaj, Junghoon Lee,
  et~al.
\newblock Immersed finite element method and its applications to biological
  systems.
\newblock {\em Computer {M}ethods in {A}pplied {M}echanics and {E}ngineering},
  195(13-16):1722--1749, 2006.

\bibitem{McQueen89}
David~M McQueen and Charles~S Peskin.
\newblock {A three--dimensional computational method for blood flow in the
  heart. II. Contractile fibers}.
\newblock {\em {Journal of Computational Physics}}, 82(2):289--297, 1989.

\bibitem{Peskin72}
Charles~S Peskin.
\newblock Flow patterns around heart valves: a numerical method.
\newblock {\em {Journal of Computational Physics}}, 10(2):252--271, 1972.

\bibitem{Peskin77}
Charles~S Peskin.
\newblock Numerical analysis of blood flow in the heart.
\newblock {\em {Journal of Computational Physics}}, 25(3):220--252, 1977.

\bibitem{Peskin93}
Charles~S Peskin and Beth~Feller Printz.
\newblock Improved volume conservation in the computation of flows with
  immersed elastic boundaries.
\newblock {\em {Journal of Computational Physics}}, 105(1):33--46, 1993.

\bibitem{Reese99}
Stefanie Reese, Martin K{\"u}ssner, and Batmanathan~Dayanand Reddy.
\newblock A new stabilization technique for finite elements in non-linear
  elasticity.
\newblock {\em International {J}ournal for {N}umerical {M}ethods in
  {E}ngineering}, 44(11):1617--1652, 1999.

\bibitem{Seo11}
Jung~Hee Seo and Rajat Mittal.
\newblock A sharp-interface immersed boundary method with improved mass
  conservation and reduced spurious pressure oscillations.
\newblock {\em {Journal of Computational Physics}}, 230(19):7347--7363, 2011.

\bibitem{Stein16}
David~B Stein, Robert~D Guy, and Becca Thomases.
\newblock Immersed boundary smooth extension: a high-order method for solving
  {PDE} on arbitrary smooth domains using {F}ourier spectral methods.
\newblock {\em Journal of {C}omputational {P}hysics}, 304:252--274, 2016.

\bibitem{Udaykumar01}
HS~Udaykumar, R~Mittal, P~Rampunggoon, and A~Khanna.
\newblock {A sharp interface Cartesian grid method for simulating flows with
  complex moving boundaries}.
\newblock {\em {Journal of Computational Physics}}, 174(1):345--380, 2001.

\bibitem{Vadala18}
Ben Vadala-Roth, Simone Rossi, and Boyce~E Griffith.
\newblock Stabilization approaches for the hyperelastic immersed boundary
  method for problems of large-deformation incompressible elasticity.
\newblock {\em arXiv preprint arXiv:1811.06620}, 2018.

\bibitem{Wang04}
Xiaodong Wang and Wing~Kam Liu.
\newblock Extended immersed boundary method using {FEM} and {RKPM}.
\newblock {\em {Computer Methods in Applied Mechanics and Engineering}},
  193(12-14):1305--1321, 2004.

\bibitem{Ye99}
Tao Ye, Rajat Mittal, HS~Udaykumar, and Wei Shyy.
\newblock {An accurate Cartesian grid method for viscous incompressible flows
  with complex immersed boundaries}.
\newblock {\em {Journal of Computational Physics}}, 156(2):209--240, 1999.

\bibitem{Zhang04}
Lucy Zhang, Axel Gerstenberger, Xiaodong Wang, and Wing~Kam Liu.
\newblock Immersed finite element method.
\newblock {\em Computer {M}ethods in {A}pplied {M}echanics and {E}ngineering},
  193(21-22):2051--2067, 2004.

\end{thebibliography}

\clearpage
\appendix

\section{Exact solution for thick inflating ring}

Consider a two--dimensional thick ring with with inner radius $R_\text{in}$ and outer radius $R_\text{out}$. 
We inject fluid on the interior of the ring with volume $A_\text{add}$.  By incompressibility, the deformed radius within the ring is
\begin{align*}
r(R) = \left(R^2 + \frac{A_\text{add}}{\pi} \right)^{1/2}, \quad R_\text{in} \leq R \leq R_\text{out}.
\end{align*} 
This uses a motion map in polar coordinates, which may be expressed as:
\begin{align*}
{\bs \chi}_\text{polar}(R, \Theta) = \left( r(R), \Theta \right)
\end{align*}
The deformation gradient in polar coordinates is:
\begin{align*}
\mathbb{F}_\text{polar} = 
\begin{bmatrix}
\frac{R}{\left(R^2 + \frac{A_\text{add}}{\pi} \right)^{1/2}} & 0 \\
0 & 1
\end{bmatrix} = 
\begin{bmatrix}
\frac{\left(r^2 - \frac{A_\text{add}}{\pi} \right)^{1/2}}{r} & 0 \\
0 & 1
\end{bmatrix}.
\end{align*}
Let us consider the following material model:
\begin{align*}
\mathbb{P}_\text{polar}^\text{e} = \mu_e \left(\mathbb{F}_\text{polar} - \mathbb{F}_\text{polar}^{-T} \right)
\end{align*}
Assuming incompressibility, $J = 1$, the elastic Cauchy stress is
\begin{align*}
{\bs \sigma}^\text{e,polar} = \mu_e \left(\mathbb{F}_\text{polar}\, \mathbb{F}_\text{polar}^{T} - {\bs I}\right) =
\begin{bmatrix}
-\mu_e\frac{A_\text{add}}{\pi r^2} & 0 \\
0 & 0
\end{bmatrix}.
\end{align*}
By force balance, at steady state, we have:
\begin{align*}
\frac{1}{r} \frac{\partial}{\partial r}\left(r {\bs \sigma}^\text{e,polar}_{rr}\right) = \frac{\partial p}{\partial r},
\end{align*} 
This implies that 
\begin{align*}
p(r) =  -\frac{\mu_e}{2}\frac{A_\text{add}}{\pi r^2} + \text{constant}, \quad r_\text{in} \leq r \leq r_\text{out} .
\end{align*}
If we hold the pressure on the outside of the ring equal to zero, we can determine this constant from the jump condition on the outer radius. 
More precisely:
\begin{align*}
[p]|_{r_\text{out}} =p^+ - p^- = -{\bs n} \cdot {\bs \sigma}_\text{polar}^\text{e}\, {\bs n} = \mu_e \frac{A_\text{add}}{\pi r_\text{out}^2}.
\end{align*}
Choosing $p^+ = 0$ implies that $p^- = -\mu_e \frac{A_\text{add}}{\pi r_\text{out}^2}$.  But since
\begin{align*}
p^- = p(r_\text{out}) = -\frac{\mu_e}{2}\frac{A_\text{add}}{\pi r^2} + \text{constant}
\end{align*}
we have determined the constant to be $-\frac{\mu_e}{2}\frac{A_\text{add}}{\pi r_\text{out}^2}$.  Thus, we obtain
\begin{align*}
p(r) =  -\frac{\mu_e A_\text{add}}{2\pi}\left( \frac{1}{r^2} + \frac{1}{r_\text{out}^2} \right) , \quad r_\text{in} \leq r \leq r_\text{out}.
\end{align*}
The jump condition for the pressure on the interior radius allows us to derive a formula for the pressure of the fluid interior to the ring.  Similarly, we have
\begin{align*}
[p]|_{r_\text{in}} =p^+ - p^- = -{\bs n} \cdot {\bs \sigma}_\text{polar}^\text{e}\, {\bs n} = \mu_e \frac{A_\text{add}}{\pi r_\text{in}^2}.
\end{align*}
Solving for $p^+$ gives us a formula for the pressure of the fluid interior to the ring.  The pressure is:
\begin{align*}
p(r) = 
\begin{cases}
 \frac{\mu_e A_\text{add}}{2\pi}\left( \frac{1}{r_\text{in}^2} - \frac{1}{r_\text{out}^2} \right) , \quad  r \leq r_\text{in} \\
  -\frac{\mu_e A_\text{add}}{2\pi}\left( \frac{1}{r^2} + \frac{1}{r_\text{out}^2} \right) , \quad r_\text{in} \leq r \leq r_\text{out} \\
0, \quad r > r_\text{out} 
\end{cases}
\end{align*}
As a remark, the implementation of this model uses Cartesian coordinates, and care is needed when converting between Cartesian and polar coordinate systems. 
Define the reference polar coordinates to be ${\bs P} = (R, \Theta)$ and the current polar coordinates to be ${\bs p} = (r,\theta)$.  
Converting the polar Cauchy stress ${\bs \sigma}_\text{polar}^\text{e}$ to the Cartesian version ${\bs \sigma}_\text{cart}^\text{e}$, with $x = r\cos\theta$ and $y = r\sin\theta$, uses the Givens rotation:
\begin{align*}
{\bs \sigma}_\text{cart}^\text{e} = 
\begin{bmatrix}
\cos \theta & -\sin \theta \\
\sin \theta & \cos \theta
\end{bmatrix}
{\bs \sigma}_\text{polar}^\text{e}
\begin{bmatrix}
\cos \theta & \sin \theta \\
-\sin \theta & \cos \theta
\end{bmatrix}.
\end{align*}
Also, the chain rule can be used to convert the Cartesian deformation gradient $\mathbb{F}_\text{cart}$ to the polar deformation gradient $\mathbb{F}_\text{polar}$:
\begin{align*}
\mathbb{F}_\text{polar} = \frac{\partial {\bs \chi}_\text{polar}}{\partial {\bs P}} =  \frac{\partial {\bs \chi}_\text{polar}}{\partial {\bs \chi}_\text{cart}} \frac{\partial {\bs \chi}_\text{cart}}{\partial {\bs X}} \frac{\partial {\bs X}}{\partial {\bs P}} = \frac{\partial {\bs p}}{\partial {\bs x}} \mathbb{F}_\text{cart} \frac{\partial {\bs X}}{\partial {\bs P}}.
\end{align*}

\end{document}